\documentclass[11pt]{article}

\usepackage{epsfig,amsmath,latexsym,amssymb, amsthm, enumerate, mathtools}
\usepackage{graphicx}
\usepackage{lscape}
\usepackage{picture, eso-pic, tikz} 
\usepackage{dsfont}
\usepackage{listings}
\usepackage{changes}
\usepackage{hyperref}
\usepackage{bbding}
\usepackage{booktabs}
\usepackage{comment}
\usepackage{multicol}

\renewcommand{\baselinestretch}{1.1}
\def\R{{\mathbb R}}  
\def\N{{\mathbb N}}  
\def\p{{\mathbb P}}  
\def\Q{{\mathbb Q}}  
\def\E{{\mathbb E}}  %
\def\X{{\mathcal X}}  
\def\Beweis{\footnotesize}

\newcommand{\Remm}[1]{}
\newtheorem{theo}{Theorem}[section]
\newtheorem{lemma}[theo]{Lemma}
\newtheorem{prop}[theo]{Proposition}

\newtheorem{defi}[theo]{Definition}
\newtheorem{model ass}[theo]{Model Assumptions}

\newtheorem{ass}[theo]{Assumption}

\newcommand*\interior[1]{\mathring{#1}}

\newtheorem{rem}[theo]{Remark}
\def\EndProof{\hfill {\scriptsize $\Box$}}

\numberwithin{equation}{section}

\definecolor{MyGray}{rgb}{0.92,0.92,0.92}


\definecolor{British racing}{rgb}{0.0, 0.5, 0.0}
\def\bx{\boldsymbol{x}}

\def\bY{\boldsymbol{Y}}
\def\bX{\boldsymbol{X}}

\def\b0{\boldsymbol{0}}
\def\bone{\boldsymbol{1}}

\def\bmu{\boldsymbol{\mu}}

\def\b0{\boldsymbol{0}}

\def\bv{\boldsymbol{v}}

\begin{document}
\author{{\L}ukasz Delong\footnote{University of Warsaw, Faculty of Economic Sciences,
l.delong@uw.edu.pl} \and 
Selim Gatti\footnote{RiskLab, Department of Mathematics, ETH Zurich,
selim.gatti@math.ethz.ch} \and 
  Mario V.~W\"uthrich\footnote{RiskLab, Department of Mathematics, ETH Zurich,
mario.wuethrich@math.ethz.ch}}

\date{Version of October 8, 2025}
\title{Calibration Bands for Mean Estimates within the Exponential Dispersion Family}
\maketitle

\begin{abstract} A statistical model is said to be calibrated if the resulting mean estimates perfectly match the true means of the underlying responses. Aiming for calibration is often not achievable in practice as one has to deal with finite samples of noisy observations. A weaker notion of calibration is auto-calibration. An auto-calibrated model satisfies that the expected value of the responses for a given mean estimate matches this estimate. Testing for auto-calibration has only been considered recently in the literature and we propose a new approach based on calibration bands. Calibration bands denote a set of lower and upper bounds such that the probability that the true means lie simultaneously inside those bounds exceeds some given confidence level. Such bands were constructed by Yang--Barber (2019) for sub-Gaussian distributions. Dimitriadis et al.~(2023) then introduced narrower bands for the Bernoulli distribution. We use the same idea in order to extend the construction to the entire exponential dispersion family that contains for example the binomial, Poisson, negative binomial, gamma and normal distributions. Moreover, we show that the obtained calibration bands allow us to construct various tests for calibration and auto-calibration, respectively. As the construction of the bands does not rely on asymptotic results, we emphasize that our tests can be used for any sample size.
\noindent

~

\end{abstract}

\section{Introduction}


Various statistical methods can be used to derive mean estimates from available observations, and it is important to understand whether these mean estimates are reliable for decision making. 
A statistical model is said to be {\it calibrated} if the resulting mean estimates perfectly match the true means of the underlying responses. In practice, calibration is often not achievable, as estimates are obtained from finite samples of noisy observations. A desirable property for a statistical model is \textit{auto-calibration}, which is a related and weaker notion of calibration; see Krüger--Ziegel \cite{Kruger} and Gneiting--Resin \cite{Gneiting}. This property means that when the responses are partitioned according to their mean estimates, i.e., responses with equal mean estimates are grouped, the estimated mean within each of these groups should match the expected value of the responses of that group. Pohle \cite{Pohle}, Gneiting--Resin \cite{Gneiting}, Krüger--Ziegel \cite{Kruger}, Denuit et al.~\cite{Denuit}, Fissler et al.~\cite{Fissler} as well as Wüthrich--Merz \cite{WM2023} emphasize the importance of auto-calibration when assessing a fitted model, especially for insurance pricing, because an auto-calibrated pricing system avoids systematic cross-subsidy.

Testing for auto-calibration has only been considered recently in the literature. Denuit et al.~\cite{Denuit2} propose a test using Lorenz and concentration curves that requires the evaluation of a non-explicit asymptotic distribution using Monte-Carlo simulations. Simpler versions of this test are provided in Wüthrich \cite{Wüthrich} for discrete and finite regression functions. Additionally, Delong--Wüthrich \cite{Delong--Wüthrich} consider the use of bootstrap techniques to assess the auto-calibration of a model. In the special case of binary observations, Hosmer--Lemeshow \cite{Hosmer--Lemeshow} derive a $\chi^2$-test by binning observations over disjoint intervals, whereas Gneiting--Resin \cite{Gneiting} propose a bootstrap approach to test for auto-calibration in this binary setup.

We take a different approach in this paper. Our goal is to construct \textit{calibration bands} for mean estimates within the exponential dispersion family (EDF). A calibration band denotes a set of lower and upper bounds for each mean estimate such that the probability of having simultaneously all the true means lying inside these bounds exceeds a given confidence level. This allows us to assess the calibration of a model by evaluating whether the mean estimates fall inside these bounds, and as the construction of the band does not rely on any asymptotical results, one can use it to construct statistical tests for calibration and auto-calibration for all sample sizes, in contrast to the approaches mentioned above.


Calibration bands were first constructed by Yang--Barber \cite{Yang_Barber} for mean estimates of sub-Gaussian distributions, which are distributions having similar or lighter tails than a Gaussian distribution as, for example, the binomial and the normal distributions. Dimitriadis et al.~\cite{Dimitriadis} then provided another construction in the binary case and showed that the resulting calibration bands are narrower than Yang--Barber's bands for the same case. Our construction of the calibration bands is similar to the construction of Dimitriadis et al.~\cite{Dimitriadis}. We extend their results to the entire EDF by exploiting stochastic ordering results and convolution formulas within the EDF.

The EDF is a broad class of distributions commonly used in statistical modeling and, particularly, in generalized linear models (GLMs); we refer the reader to McCullagh--Nelder \cite{Nelder--Wedderburn}, J\o rgensen \cite{Jorgenssen1, Jorgenssen2} and Barndorff-Nielsen \cite{Barndorff-Nielsen}.  
We provide a general construction of the calibration bands for the EDF, and we show that these bands can be expressed explicitly for the binomial, Poisson, negative binomial, gamma and normal cases.
Our bands are identical to the ones derived by Dimitriadis et al.~\cite{Dimitriadis} in the binary case, and we show that they are narrower than the calibration bands of Yang--Barber \cite{Yang_Barber} in the normal case.

We then extend the above construction to regression modeling, where the mean estimation task consists in approximating the conditional mean of a response given an observed set of features. 
In this framework, we construct two opposite statistical tests for calibration using calibration bands, i.e., statistical tests where the calibration property once lies in the null-hypothesis and once in the alternative. Moreover, we show that one can construct two opposite tests for auto-calibration.


\medskip
{\bf Organization.} This manuscript is organized as follows. In the next section, we introduce the EDF and state the framework under which we aim at constructing calibration bands on the mean. In Section \ref{section properties}, we outline the necessary assumptions and derive stochastic ordering results within the EDF that allow us, along with convolution formulas, to derive the bands. Then, we exploit these results by using a union bound argument in order to construct the calibration bands in Section \ref{section no regression}. In Section \ref{sec explicit bands},
we show that these bands can actually be expressed explicitly for the binomial, Poisson, negative binomial, gamma and normal distributions. In Section \ref{sec regression}, we extend the construction of the calibration bands to regression modeling, and in Section \ref{sec tests}, we introduce the auto-calibration property and provide conditions under which this property is equivalent to calibration. Moreover, we derive in the same section statistical tests that enable to test for calibration and auto-calibration of a given regression function. Finally, in Section \ref{numerical examples}, we highlight the impact of various factors on the resulting calibration bands through a series of numerical examples. The last section concludes this work. All mathematical
 proofs are provided in the appendix.

\section{Calibration bands within the exponential dispersion family}

\label{sec definitions}
A random variable $Y$ belongs to the EDF if its density can be written as
\begin{equation}
    \label{density EDF}
        f_Y(y; \theta, v, \varphi, \kappa(\cdot)) = \exp \left\{ \frac{y \theta - \kappa(\theta)}{\varphi/v} + a(y; v/\varphi) \right\}, 
    \end{equation}
where $\theta \in \Theta$ is called the canonical parameter, $\Theta$ is the effective domain, $\kappa : \Theta \to \R$ is the cumulant function, $v > 0$ is the volume, $\varphi > 0$ is the dispersion parameter and $a(y; v/\varphi)$ is a normalizing function depending only on $y$ and $v/\varphi$ such that the density integrates to one. We write $Y \sim \textrm{EDF}(\theta, v, \varphi, \kappa(\cdot))$ to denote a member of the EDF and emphasize that the density in \eqref{density EDF} is understood w.r.t.~a $\sigma$-finite measure $\nu$ on $\R$ that is independent of the specific choice of the canonical parameter $\theta \in \Theta$. In particular, the random variable $Y$ can, for example, be absolutely continuous or discrete. 

In this paper, we construct a calibration band on the mean of responses belonging to the EDF. To this end, we consider $n$ independent responses $Y_i~\sim~\textrm{EDF}(\theta_i, v_i, \varphi, \kappa(\cdot))$ for a fixed and known cumulant function $\kappa$ and a given dispersion parameter $\varphi > 0$. Let $\bY = (Y_1, \dots, Y_n)$ and denote the mean of each response by $\mu_i = \E[Y_i]$. Under the assumption that the responses are ordered such that their means are non-decreasing, i.e.,
\begin{equation}
    \label{ranking no regression}
    \mu_i \leq \mu_j \, \textrm{ whenever } \, i \leq j,
\end{equation}
we construct random variables $L^\alpha_{\bY ,i}$ and $U^\alpha_{\bY ,i}$ such that
\begin{equation}
\label{bound mean}
\p \Bigl(L^\alpha_{\bY ,i} \leq \mu_i \leq U^\alpha_{\bY ,i} \, \textrm{ for all } i \in \{1, \dots , n \}\Bigl) \, \geq \, 1- \alpha,
\end{equation}
for any given confidence level $1-\alpha \in (0,1)$. The resulting calibration band $$\left(L^\alpha_{\bY ,i}, U^\alpha_{\bY ,i}\right)_{i=1}^n$$ is {\it data-dependent} as it depends on the realizations of the random vector $\bY$, and its construction relies on the ordering assumed in \eqref{ranking no regression}, as well as on stochastic ordering properties and convolution formulas of the EDF that are discussed in the next section.

\section{Stochastic orders within the exponential dispersion family}

\label{section properties}

The construction of calibration bands on the mean within the EDF is motivated by the same idea used by Dimitriadis et al.~\cite{Dimitriadis} in order to construct calibration bands for independent binary random variables. Binary random variables have the nice property that we can aggregate them to binomial random variables, which satisfy certain stochastic ordering properties. The calibration bands constructed by Dimitriadis et al.~\cite{Dimitriadis} are based on these aggregations and stochastic orderings. In this section, we outline the assumptions and properties needed to generalize these ideas to the entire EDF by extracting similar stochastic orders. For this, we start with an assumption on the effective domain $\Theta$ and the $\sigma$-finite measures $\nu_i$ that define the supports of the responses $Y_i$.

\begin{ass}
    \label{assumption EDF}
    We assume that the effective domain has a non-empty interior $\interior{\Theta}$ and that the $\sigma$-finite measures $\nu_i$ are not a single point mass.
\end{ass}

This assumption excludes any trivial case of the EDF, and it implies that the effective domain $\Theta$ is a (possibly infinite) interval in $\R$ with a non-empty interior; we refer to J{\o}rgensen \cite{Jorgenssen1, Jorgenssen2}. Moreover, under Assumption \ref{assumption EDF}, the cumulant function $\kappa$ is smooth and strictly convex on the interior of the effective domain~$\mathring{\Theta}$, which implies that the derivative of the cumulant function $\kappa'$ can be inverted on this range. There is thus a one-to-one correspondence between the canonical parameter $\theta \in \interior{\Theta}$ and the mean of the random variable $Y \sim \textnormal{EDF}(\theta, v, \varphi, \kappa(\cdot))$ that is given by
\begin{equation*}
    \E[Y] = \kappa'(\theta).
\end{equation*}
This correspondence can be expressed as
\begin{equation} \label{canonical link}
    h(\E[Y]) = \theta,
\end{equation}
where $h = (\kappa')^{-1} : \kappa'(\interior{\Theta}) \to \interior{\Theta}$ is a strictly increasing function called the \textit{canonical link} of the chosen distribution within the EDF. In order to exploit this bijective map between means and canonical parameters, we make the following assumption and call the set $\kappa'(\interior{\Theta})$ the {\it mean parameter space}.
\begin{ass}
    \label{ass effective domain}
    We assume that the canonical parameters of all the considered random variables lie in the interior of the effective domain, i.e., $\theta_i \in \interior{\Theta}$ for all $1 \leq i \leq n$.
\end{ass}

Under Assumptions \ref{assumption EDF} and \ref{ass effective domain}, one can now derive various stochastic ordering relations that will be used to construct the calibration bands on the means in Section \ref{section no regression}. To do so, we introduce the usual stochastic order and the likelihood ratio order as in Shaked--Shanthikumar \cite{Stochastic Orders}.

\begin{defi}
    The usual stochastic order and the likelihood ratio order are defined as follows:
    \begin{itemize}
        \item A random variable $X$ is said to be smaller than a random variable $Y$ in the usual stochastic order, write
    $X \leq_{st} Y$, if
    \begin{equation*}
        \p(X \leq x) \geq \p(Y \leq x), \qquad \textrm{for all } x \in \R.
    \end{equation*}
        \item A random variable $X$ (with density $f$) is said to be smaller than a random variable $Y$ (with density $g$) in the likelihood ratio order, write $X \leq_{lr} Y$, if
    \begin{equation}
        \label{likelihood ratio order}
        t \mapsto \frac{g(t)}{f(t)}
    \end{equation}
    is a non-decreasing function in $t$, for $t$ being in the union of the supports of $f$ and $g$, and where $a/0$ is taken to be equal to $\infty$, whenever $a > 0$.
    \end{itemize}
\end{defi}

Note that the densities in \eqref{likelihood ratio order} are understood with respect to $\sigma$-finite measures on $\R$. Theorem 1.C.1 of Shaked--Shanthikumar \cite{Stochastic Orders} states that the likelihood ratio order is weaker than the usual stochastic order. That is, for any two random variables $X$ and $Y$ satisfying $X \leq_{lr} Y$, we have $X \leq_{st} Y$.
This implication leads to a first stochastic ordering result within the EDF; all proofs are provided in the appendix.
\begin{prop} \label{stochastic order EDF}
    Suppose that Assumptions \ref{assumption EDF} and \ref{ass effective domain} hold and let $\mu_1 \leq \mu_2$ be in the mean parameter space $\kappa'(\interior{\Theta})$. Then, for any volume $v > 0$, dispersion parameter $\varphi > 0$ and cumulant function $\kappa$, the random variables $Y_1 \sim \textnormal{EDF}(h(\mu_1), v, \varphi, \kappa(\cdot))$ and $Y_2~\sim~\textnormal{EDF}(h(\mu_2), v, \varphi, \kappa(\cdot))$ satisfy $Y_1 \leq_{st} Y_2$.
\end{prop}

Denote the distribution of a random variable $Y \sim \textnormal{EDF}(h(\mu), v, \varphi, \kappa(\cdot))$ for $\mu \in \kappa'(\interior{\Theta})$ by 
\begin{equation*}
    F(y; h(\mu), v, \varphi, \kappa(\cdot)) = \p(Y \leq y),
\end{equation*}
and the left-continuous, right-limit distribution of this random variable by
\begin{equation*}
    F^{*}(y; h(\mu), v, \varphi, \kappa(\cdot)) = \p(Y < y).
\end{equation*}
For fixed $y \in \R$, $v > 0$, $\varphi > 0$ and cumulant function $\kappa$, the stochastic ordering result in Proposition \ref{stochastic order EDF} implies that the functions
\begin{equation*}
    \mu \in \kappa'(\interior{\Theta}) \mapsto F(y; h(\mu), v, \varphi, \kappa(\cdot)),
\end{equation*}
and
\begin{equation*}
    \mu \in \kappa'(\interior{\Theta}) \mapsto F^{*}(y; h(\mu), v, \varphi, \kappa(\cdot)),
\end{equation*}
are non-increasing in $\mu$. This observation leads to the construction of the bounds on the mean in the next proposition. 

\begin{prop} \label{bound one EDF}
    Suppose that Assumptions \ref{assumption EDF} and \ref{ass effective domain} hold. Let $Y \sim \textnormal{EDF}(\theta, v, \varphi, \kappa(\cdot))$ for $\theta \in \interior{\Theta}$, $\delta \in (0,1)$, and define the random variables
    \begin{equation}\label{lower bound delta}
        l^{\delta}(Y, v, \varphi, \kappa(\cdot)) = \inf \left\{ \left. \mu \in \kappa'(\interior{\Theta}) \, \right| \, F^{*}\left(Y; h(\mu), v, \varphi, \kappa(\cdot)\right) \leq 1-\delta \right\},
    \end{equation}
    and
    \begin{equation}\label{upper bound delta}
        u^{\delta}(Y, v, \varphi, \kappa(\cdot)) = \sup \left\{ \left. \mu \in \kappa'(\interior{\Theta}) \, \right| \, F\left(Y; h(\mu), v, \varphi, \kappa(\cdot)\right) \geq \delta \right\}.
    \end{equation}
    Then,
    \begin{equation*}
        \p \left( \E[Y] \geq l^{\delta}(Y, v, \varphi, \kappa(\cdot)) \right) \geq 1- \delta \quad \textrm{and} \quad \p \left( \E[Y] \leq u^{\delta}(Y, v, \varphi, \kappa(\cdot)) \right) \geq 1- \delta.
    \end{equation*}
\end{prop}

Proposition \ref{bound one EDF} provides lower and upper bounds holding for the case of a single response $Y \sim \textnormal{EDF}(\theta, v, \varphi, \kappa(\cdot))$. These bounds directly depend on the value of $\delta \in (0,1)$ as for any $y \in \R$, the interval $$\left[l^{\delta}(y, v, \varphi, \kappa(\cdot)),u^{\delta}(y, v, \varphi, \kappa(\cdot))\right]$$ is wide for small values of $\delta \in (0,1/2]$ and narrow for large values of $\delta \in (0,1/2]$. Additionally, note that for $\delta \geq 1/2$, we might even have that the lower bound exceeds the upper bound, i.e.,
$$l^{\delta}(y, v, \varphi, \kappa(\cdot)) \geq u^{\delta}(y, v, \varphi, \kappa(\cdot)).$$

To lift this result to the case of $n$ independent responses being ordered such that their canonical parameters (or means) are increasing, see \eqref{ranking no regression}, we aim at using Proposition \ref{bound one EDF} to derive lower and upper bounds on the weighted partial sums
\begin{equation}
    \label{sum j k}
    Z_{j:k} = \frac{1}{v_{j:k}}\sum_{i=j}^k v_i Y_i,
\end{equation}
 for any $1 \leq j \leq k \leq n$, and where we define the aggregated volumes as
 \begin{equation}
    \label{aggregated volume}
     v_{j:k} = \sum_{i=j}^k v_i.
 \end{equation}

 For this, we use the following stochastic bounds on the random variables $Z_{j:k}$.
\begin{lemma} \label{stochastic bound j k}
    Let $Y_j, \dots, Y_k$ be independent $\textnormal{EDF}(\theta_i, v_i, \varphi, \kappa(\cdot))$ distributed random variables for given volumes $v_i > 0$ and indices $j\leq i \leq k$ such that $\theta_j \leq \dots \leq \theta_k$. Under Assumptions \ref{assumption EDF} and \ref{ass effective domain}, the weighted sum $Z_{j:k}$ in \eqref{sum j k}
    satisfies $$Z_{j:k}^{-} \leq_{st} Z_{j:k} \leq_{st} Z_{j:k}^{+},$$ for
    \begin{equation}  \label{random variables jk}
        Z_{j:k}^{-} \sim \textnormal{EDF}(\theta_j, v_{j:k}, \varphi, \kappa(\cdot)) \quad \textrm{ and } \quad Z_{j:k}^{+} \sim \textnormal{EDF}(\theta_k, v_{j:k}, \varphi, \kappa(\cdot)).
    \end{equation}
\end{lemma}

\medskip

Using these stochastic bounds and Proposition \ref{bound one EDF}, we can now derive bounds on the means of the weighted partial sums $Z_{j:k}$.

\begin{prop} \label{last corollary}
    Suppose that Assumptions \ref{assumption EDF} and \ref{ass effective domain} hold. Let $Y_j, \dots, Y_k$ be independent $\textnormal{EDF}(\theta_i, v_i, \varphi, \kappa(\cdot))$ distributed random variables for given volumes $v_i > 0$ and indices $j\leq i \leq k$ such that $\theta_j \leq \dots \leq \theta_k$. Moreover, let $\delta \in (0,1)$, define $Z_{j:k}$ and $v_{j:k}$ as in \eqref{sum j k}-\eqref{aggregated volume} and denote by $\mu_j$ and $\mu_k$ the means of $Y_j$ and $Y_k$, respectively. Then, we have
    \begin{equation*}
        \p\left(  \mu_j \leq  u^{\delta}(Z_{j:k}, v_{j:k}, \varphi, \kappa(\cdot))   \right) \geq 1- \delta \quad \textrm{and} \quad \p\left(  \mu_k \geq  l^{\delta}(Z_{j:k}, v_{j:k}, \varphi, \kappa(\cdot))   \right) \geq 1- \delta,
    \end{equation*}
    for the random variables $l^{\delta}(Z_{j:k}, v_{j:k}, \varphi, \kappa(\cdot))$ and $u^{\delta}(Z_{j:k}, v_{j:k}, \varphi, \kappa(\cdot))$ defined in \eqref{lower bound delta}-\eqref{upper bound delta}.
\end{prop}

As we will outline in the next section, Proposition \ref{last corollary} is in fact at the core of the construction of the calibration bands on the mean holding for $n$ independent responses because for any indices $j \le k$, it provides upper and lower bounds on the true means $\mu_j$ and $\mu_k$ that only depend on the realizations of the responses $Y_j,\dots, Y_k$, the given volumes $v_j, \dots, v_k$ as well as the dispersion parameter $\varphi$ and cumulant function $\kappa$.

\section{Construction of the calibration bands}

\label{section no regression}

\subsection{Main result}
The aim of this section is to construct calibration bands on the mean, as defined in \eqref{bound mean}, for independent responses $(Y_i)_{i=1}^n \sim \textrm{EDF}(\theta_i, v_i, \varphi, \kappa(\cdot))$ that are ordered such that their canonical parameters fulfill $\theta_1 \leq \dots \leq \theta_n$. This construction makes use of sets of ordered pairs $\mathcal{J}~\subseteq~\{1, \dots, n\}^2$ that we define as sets satisfying
\begin{equation*}
    (j,k) \in \mathcal{J} \implies j \leq k.
\end{equation*}
By using a union bound argument, a corollary of Proposition \ref{last corollary} is that
\begin{equation*}
    \p \Bigl(\mu_j \leq  u^{\delta}(Z_{j:k}, v_{j:k}, \varphi, \kappa(\cdot) ) \textrm{ and } \mu_k \geq  l^{\delta}(Z_{j:k}, v_{j:k}, \varphi, \kappa(\cdot)) \Bigl) \geq 1-2\delta,
\end{equation*}
for any pair $(j,k) \in \mathcal{J}$, because the complement of the above event is nested in the event where at least one of the true means fails to lie above or below the constructed lower and upper bounds, respectively. Similarly, we have
\begin{equation}
    \label{prob bound J}
    \p \Bigl(\mu_j \leq  u^{\delta}(Z_{j:k}, v_{j:k}, \varphi, \kappa(\cdot) ) \textrm{ and } \mu_k \geq  l^{\delta}(Z_{j:k}, v_{j:k}, \varphi, \kappa(\cdot)) \textrm{ for all } (j,k) \in \mathcal{J} \Bigl) \geq 1-2|\mathcal{J}|\delta,
\end{equation}
where the bounds on the means $\mu_j$ and $\mu_k$ now hold simultaneously for all pairs $(j,k) \in \mathcal{J}$. This last inequality allows us to construct calibration bands on the mean of $n$ independent responses, as stated in the next theorem.
\begin{theo}
    \label{theo calibration bands}
    Suppose that Assumptions \ref{assumption EDF} and \ref{ass effective domain} hold. Let $Y_1, \dots, Y_n$ be independent $\textnormal{EDF}(\theta_i, v_i, \varphi, \kappa(\cdot))$ distributed random variables for given volumes $v_i > 0$ and indices $1 \leq i \leq n$ such that $\theta_1 \leq \dots \leq \theta_n$. 
    Moreover, let $\mathcal{J} \subseteq \{1, \dots, n\}^2$ be any set of ordered pairs. By writing $\mu_i = \E[Y_i]$ for $1 \leq i \leq n$, we have for any given confidence level $1-\alpha \in (0,1)$ that 
    \begin{equation}
    \label{inequality discrete case}
    \p\Bigl(  L^\alpha_{\bY, i} \leq \mu_i \leq  U^\alpha_{\bY, i} \, \textnormal{ for all } \, i \in \{1, \dots, n\} \Bigl) \geq 1- \alpha, 
\end{equation}
with
\begin{equation}\label{L_i general}
    L^\alpha_{\bY, i} = \sup_{(j,k) \in \mathcal{J} \, : \, \theta_i \geq \theta_k} l^{\delta}(Z_{j:k}, v_{j:k}, \varphi, \kappa(\cdot)),
\end{equation}
and
\begin{equation}\label{U_i general}
    U^\alpha_{\bY, i} = \inf_{(j,k) \in \mathcal{J} \, : \, \theta_i \leq \theta_j}u^{\delta}(Z_{j:k}, v_{j:k}, \varphi, \kappa(\cdot)),
\end{equation}
for $ \delta = \alpha/(2|\mathcal{J}|)$.
\end{theo}
We emphasize that the construction of the calibration band in Theorem \ref{theo calibration bands} only depends on the realizations of the random variables $(Y_i)_{i=1}^n$, the volumes $(v_i)_{i=1}^n$, the dispersion parameter $\varphi$ and the cumulant function $\kappa$. In particular, it does not depend on the means $(\mu_i)_{i=1}^n$, but only on their rankings. Moreover, it holds for any sample size, i.e., it does not rely on asymptotic sample size considerations.


\subsection{Choice of the set of ordered pairs and binning of the observations}

\label{section computational time}

The statement in Theorem \ref{theo calibration bands} holds for any set of ordered pairs $\mathcal{J} \subseteq \{1, \dots, n\}^2$ and, in fact, the resulting calibration band directly depends on the choice of this set. Moreover, note that in principle, $\mathcal{J}$ might only contain a few pairs, which could lead to take the supremum and the infimum of empty sets in \eqref{L_i general} and \eqref{U_i general}. In such cases, we adopt the convention
\begin{equation*}
    \inf \emptyset = \sup_{\theta \in \interior{\Theta}} \kappa'(\theta) \qquad \textrm{and} \qquad  \sup \emptyset = \inf_{\theta \in \interior{\Theta}} \kappa'(\theta).
\end{equation*}
That is, the underlying lower and upper bounds are equal to the infimum and the supremum of the mean parameter space, respectively. An intuitive choice for $\mathcal{J}$ is the set of all possible combinations of ordered pairs that we denote by
\begin{equation*}
        \mathcal{J}^{full} = \{(j,k) \in \{1, \dots, n\}^2 \, | \, j \leq k\}.
\end{equation*}
In this case, we call the constructed band the \textit{full calibration band}. Many other choices are possible, and we want to discuss two contrasting factors that create a trade-off situation. On the one hand, for a fixed $\delta \in (0,1)$, the lower and upper bounds in \eqref{L_i general} and \eqref{U_i general} lead to a wider band than the full calibration band for any smaller set of ordered pairs $\mathcal{J} \subseteq \mathcal{J}^{full}$. This suggests that the set $\mathcal{J}$ should be large. On the other hand, the map $\delta \mapsto l^{\delta}(Z_{j:k}, v_{j:k}, \varphi, \kappa(\cdot))$ is non-decreasing, whereas the map $\delta \mapsto u^{\delta}(Z_{j:k}, v_{j:k}, \varphi, \kappa(\cdot))$ is non-increasing. Consequently, for a fixed set of ordered pairs $\mathcal{J}$, the resulting band becomes narrower as $\delta$ increases. However, since the value of $\delta$ is directly determined by the number of elements in the set $\mathcal{J}$, via the relation $\delta = \alpha/(2|\mathcal{J}|)$ in Theorem \ref{theo calibration bands}, a large set of ordered pairs leads to a low value for $\delta$ and vice versa. This creates a trade-off and there is thus no optimal choice for the set of ordered pairs in general.

In their construction of calibration bands for the binary case, Dimitriadis et al.~\cite{Dimitriadis} suggest to use a slightly modified version of $\mathcal{J}^{full}$ that we call $$\mathcal{J}^{distinct} = \{(j,k) \in \mathcal{D}^2 \, | \, j \leq k\},$$
where $\mathcal{D}$ is any largest subset of $\{1, \dots, n \}$ such that there are no ties in the canonical parameters, i.e., $\theta_i \neq \theta_j$ for all $i \neq j \in \mathcal{D}$. Note that using the convolution property of the EDF, one can always merge observations associated to the same canonical parameter and appropriately adapt the volumes $v_i$ before constructing the calibration band using $\mathcal{J}^{distinct}$, we refer to Corollary 2.15 of Wüthrich--Merz \cite{WM2023}.

Another consideration is the computational time required for constructing the band, which corresponds to $\mathcal{O}(|\mathcal{J}|)$. Indeed, using $\mathcal{J}^{full}$ as the set of ordered pairs leads to a computational time of $\mathcal{O}(n^2)$, which might be problematic for large datasets.
One way to improve the run time is to reduce the amount of pairs in the set $\mathcal{J}$. Another way is to reduce the number of observations by binning them even if there are no ties in the canonical parameters. We come back to those methods through the numerical examples in Section \ref{numerical examples}. 

\subsection{Crossings inside the calibration bands}

Although we call the simultaneous lower and upper bounds on the means $(L^\alpha_{\bY ,i}, U^\alpha_{\bY ,i})_{i=1}^n$ derived in Theorem \ref{theo calibration bands} a calibration band, we point out that, in general, we might have $U^\alpha_{\bY ,i} < L^\alpha_{\bY ,i}$ for some indices $1 \leq i \leq n$. This phenomenon was already observed by Dimitriadis et al.~\cite{Dimitriadis} in the binary case, and these authors argue that this typically happens when the ranking of the responses violates \eqref{ranking no regression}, i.e., when a ranking obtained from empirical data is not fully accurate. 
In order to construct calibration bands that do not exhibit any crossings, Dimitriadis et al.~\cite{Dimitriadis} propose to take the pointwise minimum (maximum) of the lower (upper) band with the isotonic regressor of the responses that is defined by
\begin{equation}
    \label{isotonic estimator}
    \widehat{\bmu}^{Iso}(\bY, \bv) = \mathop{\arg \min}\limits_{\bmu \in \R^n} \Big\{ \sum_{i=1}^n v_i(Y_i-\mu_i)^2: \mu_1 \leq \dots \leq \mu_n \Big\}, \, \textrm{ for } \bv =(v_1, \dots, v_n)^\top \in \R^n.
\end{equation}
The vector $\widehat{\bmu}^{Iso}(\bY, \bv)$ provides a non-parametric estimator of the means of $\bY$ that satisfies the ordering in \eqref{ranking no regression} by construction. In the case where the calibration band shows some crossings, we follow the proposition of Dimitriadis et al.~\cite{Dimitriadis}, and we use the following modified band
\begin{equation}
    \label{cal not crossing}
    \tilde{L}^\alpha_{\bY ,i} = \min \left(L^\alpha_{\bY ,i}, \widehat{\mu}^{Iso}_i(\bY, \bv) \right) \quad \textrm{and} \quad \tilde{U}^\alpha_{\bY ,i} = \max \left(U^\alpha_{\bY ,i}, \widehat{\mu}^{Iso}_i(\bY,\bv) \right),
\end{equation}
for $1 \leq i \leq n$. Of course, the main result of this section still holds because this makes the interval wider. That is, for any given confidence level $1-\alpha \in (0,1)$, we still have 
    \begin{equation*}
    \p\Bigl(  \tilde{L}^\alpha_{\bY ,i} \leq \mu_i \leq  \tilde{U}^\alpha_{\bY ,i} \, \textnormal{ for all } \, i \in \{1, \dots, n\} \Bigl) \geq 1- \alpha. 
\end{equation*}

\section{Explicit calibration bands for selected distributions}

\label{sec explicit bands}
The calibration band derived in Theorem \ref{theo calibration bands} can be constructed for any member of the EDF under Assumptions \ref{assumption EDF} and \ref{ass effective domain}. To do so, the evaluation of the lower and upper bounds $l^{\delta}(Z_{j:k}, v_{j:k}, \varphi, \kappa(\cdot))$ and $u^{\delta}(Z_{j:k}, v_{j:k}, \varphi, \kappa(\cdot))$ requires the use of a root-finding algorithm. For some members of the EDF, these bounds can be calculated in closed form. We give the resulting expressions for the binomial, Poisson, negative binomial, gamma and normal cases in this section. These expressions are derived using closed form characterizations for the quantiles of the above distributions.

We point out that the explicit calibration bands presented in this section could also be derived using fiducial distributions.
Fiducial distributions were introduced by Fisher \cite{Fisher1, Fisher2} in the 1930s, who aimed at providing a framework for constructing probability distributions of unknown parameters based on available observations. The use of fiducial distributions has been shown to lead to some contradictory results, see for example Chapter 5.4 in Sprott \cite{Sprott}. Therefore, such distributions may only be used under specific assumptions that are discussed by Pedersen \cite{Pedersen}. Veronese--Mellili \cite{Veronese--Mellili} compute the fiducial distributions of selected members of the EDF, including the above mentioned examples, and they show that these fiducial distributions satisfy those assumptions. As a consequence, all the explicit calibration bands derived in this section could actually also be derived using fiducial distributions.

\subsection{Discrete distributions}

\label{sec discrete distributions}

The binomial, Poisson and negative binomial distributions are members of the EDF since any random variable $N$ belonging to one of these distributions can be written as 
\begin{equation}
\label{duality transform}
 N = \frac{vY}{\varphi},   
\end{equation}
 for $Y \sim\textrm{EDF}(\theta, v, \varphi, \kappa(\cdot))$ and for a carefully chosen effective domain, volume, dispersion parameter and cumulant function, we refer to Section 3.3 of J{\o}rgenssen \cite{Jorgenssen2} and Section 2.2 of Wüthrich--Merz \cite{WM2023}. Note that the transformation \eqref{duality transform} is called the duality transformation as it provides a duality between the random variables $$Y \sim\textrm{EDF}(\theta, v, \varphi, \kappa(\cdot)) \quad \textrm{and} \quad N = \frac{vY}{\varphi}.$$ The former random variable $Y$, whose density is given in \eqref{density EDF}, is said to be in the {\it reproductive} form of the EDF, whereas the latter random variable $N$ is said to be in the {\it additive} form of the EDF, see Section 3.1 of J{\o}rgenssen \cite{Jorgenssen2}. Using the duality transformation, Theorem \ref{theo calibration bands} can be used to derive calibration bands for members of the additive form of the EDF under Assumptions \ref{assumption EDF} and \ref{ass effective domain}.
      We show in the next result that the resulting bands for the binomial, Poisson and negative binomial cases can be given in closed form using explicit expressions for the quantiles of those distributions. Note that the calibration bands derived by Dimitriadis et al.~\cite{Dimitriadis} for the binary case are contained in the binomial case, below.

\begin{prop}
    \label{prop discrete distributions}
    Suppose that Assumptions \ref{assumption EDF} and \ref{ass effective domain} hold. Let $N_1, \dots, N_n$ be independent members of the EDF in the additive form, i.e., there exist canonical parameters $(\theta_i)_{i=1}^n$, volumes $(v_i)_{i=1}^n$, a dispersion parameter $\varphi > 0$ and a cumulant function $\kappa$ such that $Y_i = \varphi N_i/v_i$, for $Y_i \sim \textrm{EDF}(\theta_i, v_i, \varphi, \kappa(\cdot))$ and $1 \leq i \leq n$. By writing $\mu_i = \E[Y_i]$ and assuming that $\mu_1 \leq \dots \leq \mu_n$, we have for any set of ordered pairs $\mathcal{J} \subseteq \{1, \dots, n\}^2$ and any confidence level $1-\alpha \in (0,1)$ that
    \begin{equation*}
    \p\Bigl(  L^\alpha_{\bY, i} \leq \mu_i \leq  U^\alpha_{\bY, i} \, \textnormal{ for all } \, i \in \{1, \dots, n\} \Bigl) \geq 1- \alpha, 
\end{equation*}
    where the lower and upper bounds $L^\alpha_{\bY, i}$ and $U^\alpha_{\bY, i}$ are defined in \eqref{L_i general}-\eqref{U_i general}. These bounds can be explicitly expressed in the following three cases using the weighted partial sums $Z_{j:k}$ and aggregated volumes $v_{j:k}$ in \eqref{sum j k}-\eqref{aggregated volume}:
    \begin{itemize}
         \item Binomial case. The lower and upper bounds are given by
         \begin{equation}
    \label{L_i bin}
    L^\alpha_{\bY ,i} = \sup_{(j,k) \in \mathcal{J} \, : \, \mu_i \geq \mu_k} q_{B}(\delta; v_{j:k}Z_{j:k}/\varphi, 1+v_{j:k}/\varphi-v_{j:k}Z_{j:k}/\varphi) \, \mathds{1}_{\{Z_{j:k} > 0\}},
\end{equation}
and
\begin{equation}
    \label{U_i bin}
    U^\alpha_{\bY ,i} = \inf_{(j,k) \in \mathcal{J} \, : \, \mu_i \leq \mu_j}  q_{B}(1-\delta; 1+v_{j:k}Z_{j:k}/\varphi, v_{j:k}/\varphi-v_{j:k}Z_{j:k}/\varphi) \, \mathds{1}_{\{Z_{j:k} < 1\}} + \mathds{1}_{\{Z_{j:k} = 1\}},
\end{equation}
for $ \delta = \alpha/(2|\mathcal{J}|)$, and where $q_B(\delta; \alpha, \beta)$ denotes the $\delta$-quantile of a beta distribution with parameters $\alpha, \beta > 0$.
         \item Poisson case. The lower and upper bounds are given by
         \begin{equation}
    \label{L_i pois}
    L^\alpha_{\bY ,i} = \sup_{(j,k) \in \mathcal{J} \, : \, \mu_i \geq \mu_k} \frac{ \varphi q_{\Gamma}(\delta; v_{j:k}Z_{j:k}/\varphi, 1)}{v_{j:k}} \, \mathds{1}_{\{ Z_{j:k} > 0 \}},
\end{equation}
and
\begin{equation}
    \label{U_i pois}
    U^\alpha_{\bY ,i} = \inf_{(j,k) \in \mathcal{J} \, : \, \mu_i \leq \mu_j}  \frac{\varphi q_{\Gamma}(1-\delta; 1+v_{j:k}Z_{j:k}/\varphi, 1)}{v_{j:k}},
\end{equation}
for $ \delta = \alpha/(2|\mathcal{J}|)$, and where $q_\Gamma(\delta; \gamma, c)$ denotes the $\delta$-quantile of a gamma distribution with shape parameter $\gamma > 0$ and scale parameter $c > 0$.
\item Negative binomial case. The lower and upper bounds are given by
\begin{equation}
    \label{L_i neg bin}
    L^\alpha_{\bY ,i} = \sup_{(j,k) \in \mathcal{J} \, : \, \mu_i \geq \mu_k} \frac{q_{B}(\delta; v_{j:k}Z_{j:k}/\varphi, v_{j:k}/\varphi)}{1-q_{B}(\delta; v_{j:k}Z_{j:k}/\varphi, v_{j:k}/\varphi)} \, \mathds{1}_{\{Z_{j:k} > 0\}},
\end{equation}
and
\begin{equation}
    \label{U_i neg bin}
    U^\alpha_{\bY ,i} = \inf_{(j,k) \in \mathcal{J} \, : \, \mu_i \leq \mu_j}  \frac{q_{B}(1-\delta; 1+v_{j:k}Z_{j:k}/\varphi, v_{j:k}/\varphi)}{1-q_{B}(1-\delta; 1+v_{j:k}Z_{j:k}/\varphi, v_{j:k}/\varphi)},
\end{equation}
for $ \delta = \alpha/(2|\mathcal{J}|)$.
\end{itemize}
 \end{prop}

\subsection{Continuous distributions}

The gamma and normal distributions are members of the EDF. This time, note that any random variable $Y$ belonging to one of these distributions can directly be written as an $\textrm{EDF}(\theta, v, \varphi, \kappa(\cdot))$ random variable for a carefully chosen effective domain, volume, dispersion parameter and cumulant function, we refer to Section 3.3 of J{\o}rgenssen \cite{Jorgenssen2} or Section 2.2 of Wüthrich--Merz \cite{WM2023}. That is, the normal and the gamma distributions can directly be expressed in the reproductive form of the EDF. Moreover, these distributions satisfy Assumptions \ref{assumption EDF} and \ref{ass effective domain}, which allows us to derive calibration bands on the mean of gamma or normal responses. As above, a closed form expression for the calibration bands can be obtained using the quantiles of these distributions.

\begin{prop}
    \label{prop continuous distributions}
    Let $\mathcal{J} \subseteq \{1, \dots, n\}^2$ be any set of ordered pairs. For any confidence level $1-\alpha \in (0,1)$, the calibration band in Theorem \ref{theo calibration bands} can be explicitly expressed in the following two cases using the weighted partial sums $Z_{j:k}$ and aggregated volumes $v_{j:k}$ in \eqref{sum j k}-\eqref{aggregated volume}:

    \begin{itemize}
         \item Gamma case. The lower and upper bounds in \eqref{L_i general}-\eqref{U_i general} are given by
\begin{equation}
    \label{L_i gamma}
    L^\alpha_{\bY ,i} = \sup_{(j,k) \in \mathcal{J} \, : \, \mu_i \geq \mu_k} \frac{v_{j:k}/\varphi}{ q_\Gamma(1-\delta;v_{j:k}/\varphi,Z_{j:k})} \, \mathds{1}_{\{Z_{j:k} > 0\}},
\end{equation}
and
\begin{equation}
    \label{U_i gamma}
    U^\alpha_{\bY ,i} = \inf_{(j,k) \in \mathcal{J} \, : \, \mu_i \leq \mu_j}  \frac{v_{j:k}/\varphi}{ q_\Gamma(\delta;v_{j:k}/\varphi,Z_{j:k})} \, \mathds{1}_{\{Z_{j:k} > 0\}},
\end{equation}
for $ \delta = \alpha/(2|\mathcal{J}|)$.
         \item Normal case. The lower and upper bounds in \eqref{L_i general}-\eqref{U_i general} are given by
 \begin{equation}
    \label{L_i normal}
    L^\alpha_{\bY ,i} = \sup_{(j,k) \in \mathcal{J} \, : \, \mu_i \geq \mu_k} Z_{j:k} - \frac{\Phi^{-1}(1-\delta)}{\sqrt{v_{j:k}/\varphi}},
\end{equation}
and
\begin{equation}
    \label{U_i normal}
    U^\alpha_{\bY ,i} = \inf_{(j,k) \in \mathcal{J} \, : \, \mu_i \leq \mu_j}  Z_{j:k} - \frac{\Phi^{-1}(\delta)}{\sqrt{v_{j:k}/\varphi}},
\end{equation}
for $ \delta = \alpha/(2|\mathcal{J}|)$, and where $\Phi^{-1}(\delta)$ denotes the $\delta$-quantile of the standard normal distribution.
\end{itemize}
 \end{prop}

\begin{rem}
    \label{remark}
    \textnormal{Let $Y_1, \dots, Y_n$ be independent $\mathcal{N}(\mu_i, \sigma_i^2)$ random variables for known standard deviations $\sigma_i > 0$ and indices $1 \leq i \leq n$ such that $\mu_1 \leq \dots \leq \mu_n$. The underlying aggregated volumes and weighted partial sums are given by 
    \begin{equation*}
        v_{j:k} = \sum_{i=j}^k \frac{1}{\sigma_i^2} \quad \textrm{and} \quad Z_{j:k} = \frac{1}{v_{j:k}}\sum_{i=j}^k v_i Y_i,
    \end{equation*}
    where we set $\varphi = 1$, because the distribution of an $\textrm{EDF}(\theta, v, \varphi, \kappa(\cdot))$ distributed random variable only depends on the ratio $v/\varphi$. These sums correspond to weighted sums of normal random variables and these weights are determined by Theorem \ref{theo calibration bands} for the general EDF case. Note that large weights are given to responses $Y_i$ that are associated to a small variance and vice versa. 
    The resulting weighted partial sums $Z_{j:k}$ are thus scaled sums of $\mathcal{N}(\mu_i/\sigma_i^2, 1/\sigma_i^2)$ random variables. 
    Since the normal distribution has the nice property that any weighted sum of independent normal responses is again normal, other weights could in principle be chosen as, for example,
    \begin{equation*}
        \tilde{v}_{j:k} = \sum_{i=j}^k \frac{1}{\sigma_i} \quad \textrm{and} \quad \tilde{Z}_{j:k} = \frac{1}{\tilde{v}_{j:k}}\sum_{i=j}^k \frac{Y_i}{\sigma_i},
    \end{equation*}
    which results in weighted partial sums $\tilde{Z}_{j:k}$ being scaled sums of normal random variables that all have variance 1. The resulting calibration band can be expressed as
    \begin{equation*}
    \tilde{L}^\alpha_{\bY ,i} = \sup_{(j,k) \in \mathcal{J} \, : \, \mu_i \geq \mu_k} \tilde{Z}_{j:k} - \frac{\Phi^{-1}(1-\delta) \sqrt{k-j+1}}{\tilde{v}_{j:k}},
\end{equation*}
and
\begin{equation*}
    \tilde{U}^\alpha_{\bY ,i} = \inf_{(j,k) \in \mathcal{J} \, : \, \mu_i \leq \mu_j}  \tilde{Z}_{j:k} - \frac{\Phi^{-1}(\delta) \sqrt{k-j+1}}{\tilde{v}_{j:k}}.
\end{equation*}
This new calibration band is in general different from the one derived in Proposition \ref{prop continuous distributions}. However, both calibration bands coincide in the case of a constant variance for the responses, i.e., when $\sigma_i = \sigma$, for $1 \leq i \leq n$.
    }
\end{rem}

\subsection{Comparison with Yang--Barber's calibration bands}

In the literature, calibration bands on the mean have first been constructed by Yang--Barber \cite{Yang_Barber}, under the assumption that the responses $(Y_i)_{i=1}^n$ satisfy the additive relation
\begin{equation*}
    Y_i = \mu_i + \sigma \varepsilon_i,
\end{equation*}
for means $\mu_1 \leq \dots \leq \mu_n$, for some fixed and known $\sigma > 0$ and for independent and zero-mean random variables $\varepsilon_i$ that are sub-Gaussian, i.e.,
\begin{equation*}
    \p(|\varepsilon_i| > t) \leq 2 e^{-t^2/2}, \quad \textrm{ for all } t >0 \textrm{ and for all } 1 \leq i \leq n.
\end{equation*}
The construction of their calibration bands makes use of the isotonic regressor of $\bY$ defined in \eqref{isotonic estimator} in order to introduce the empirical partial sums
\begin{equation*}
    Z_{j:k}^{Iso} = \frac{1}{k-j+1} \sum_{i = j}^k \widehat{\mu}^{Iso}_i(\bY, \bone),
\end{equation*}
for $\bone = (1, \dots,1)^\top \in \R^n$. Yang--Barber \cite{Yang_Barber} show that for the set of ordered pairs 
$$\mathcal{J}^{full} = \{(j,k) \in \{1, \dots, n\}^2 \, | \, j \leq k\},$$
and for any confidence level $1-\alpha \in (0,1)$, we have
\begin{equation*}
    \p\Bigl(  L^{\alpha,YB}_{\bY ,i} \leq \mu_i \leq  U^{\alpha,YB}_{\bY ,i} \, \textnormal{ for all } \, i \in \{1, \dots, n\} \Bigl) \geq 1- \alpha,
\end{equation*}
with
\begin{equation}
\label{L_i YB}
    L^{\alpha,YB}_{\bY ,i} = \sup_{(j,k) \in \mathcal{J}^{full} \, : \, \mu_i \geq \mu_k} Z_{j:k}^{Iso} - \frac{\sqrt{2 \sigma^2 \log(1/\delta)}}{\sqrt{k-j+1}},
\end{equation}
and
\begin{equation}
\label{U_i YB}
    U^{\alpha,YB}_{\bY ,i} = \inf_{(j,k) \in \mathcal{J}^{full} \, : \, \mu_i \leq \mu_j}  Z_{j:k}^{Iso} + \frac{\sqrt{2 \sigma^2 \log(1/\delta)}}{\sqrt{k-j+1}},
\end{equation}
for $ \delta = \alpha/(2|\mathcal{J}^{full}|) = \alpha/(n^2+n)$. 
\medskip

A particular case of the framework used by Yang--Barber \cite{Yang_Barber} arises by taking i.i.d.~Gaussian random variables $\varepsilon_i \sim \mathcal{N}(0,1)$. In this case, we obtain independent responses $Y_i \sim \mathcal{N}(\mu_i, \sigma^2)$, for $1 \leq i \leq n$. 
Similar to Dimitriadis et al.~\cite{Dimitriadis} for the binary case, we show that our calibration band is narrower than Yang--Barber's band in this setting.

\begin{theo}
    \label{our bands better}
    Let $Y_1, \dots, Y_n$ be independent $\mathcal{N}(\mu_i, \sigma^2)$ random variables for a known standard deviation $\sigma > 0$ and indices $1 \leq i \leq n$ such that $\mu_1 \leq \dots \leq \mu_n$. The bands derived in Proposition \ref{prop continuous distributions} using the set of ordered pairs $\mathcal{J}^{full}$ satisfy
    \begin{equation*}
       L^{\alpha}_{\bY ,i}  \geq L^{\alpha,YB}_{\bY ,i} \quad \textrm{and} \quad U^{\alpha}_{\bY ,i} \leq U^{\alpha,YB}_{\bY ,i},
    \end{equation*}
    for all $1 \leq i \leq n$ and for any confidence level $1-\alpha \in (0,1)$.
\end{theo}

The proof of this theorem relies on Proposition B1 of Dimitriadis et al.~\cite{Dimitriadis} that characterizes the pairs $(i,j) \in \mathcal{J}^{full}$ for which the maximum in \eqref{L_i YB} and the minimum in \eqref{U_i YB} are attained. It is provided in the appendix.

\section{Extension to regression modeling}

\label{sec regression}

\subsection{Regression modeling within the exponential dispersion family}

The calibration bands on the mean derived in Theorem \ref{theo calibration bands} can be extended to regression modeling. To this end, let $(\Omega, \mathcal{F}, \p)$ be the underlying probability space and consider an independent sample $(Y_i, \bX_i)_{i=1}^n$ with responses $Y_i$ and i.i.d.~features $\bX_i$ satisfying 
\begin{equation*}
    Y_i \, | \, \bX_i = \bx_i ~\sim~\textrm{EDF}(\theta (\bx_i), v_i, \varphi, \kappa(\cdot)),
\end{equation*}
for given volumes $v_i > 0$, as well as a dispersion parameter $\varphi > 0$ and a cumulant function $\kappa$ that do not depend on $i$. We denote the support of the features $\bX_i$ by $\mathcal{X}$ and call it the {\it feature space}.
The goal of a regression on the mean is to estimate the true mean function
\begin{equation}
    \label{mean function}
    \mu^{*} : \mathcal{X} \to \kappa'(\interior{\Theta}), \quad \bx \mapsto \kappa'(\theta(\bx)),
\end{equation}
where the map $\theta: \mathcal{X} \to \interior{\Theta}$ is unknown. Note that this true mean function is a strictly increasing map of the canonical parameter $\bx \mapsto \theta(\bx)$ due to the strict convexity of the cumulant function $\kappa$ under Assumptions \ref{assumption EDF} and \ref{ass effective domain}, see \eqref{mean function}. In particular, the true mean function does not depend on the volume and the dispersion parameter. Therefore, one can write
\begin{equation}
    \label{true mean function}
    \mu^{*}(\bX) = \E[Y \, | \, \bX], \quad \p\textrm{-a.s.},
\end{equation}
for any pair $(Y, \bX)$ satisfying
\begin{equation*}
    Y \, | \, \bX = \bx ~\sim~\textrm{EDF}(\theta (\bx), v, \varphi, \kappa(\cdot)),
\end{equation*}
regardless of the volume $v > 0$ and the dispersion parameter $\varphi > 0$, we refer to Section \ref{section properties}. That is, the true mean function $\mu^{*}:\X \to \kappa'(\interior{\Theta})$ maps each feature $\bx \in \mathcal{X}$ to the conditional expectation of the response $Y$, given this feature is observed. 

\subsection{Construction of the calibration bands}

\label{sec construction of calibration band}
In regression modeling, a calibration band on the mean denotes a set of lower and upper bounds such that the probability that the true mean function in \eqref{mean function} lies simultaneously inside these bounds for almost every (a.e.)~feature $\bx \in \mathcal{X}$ exceeds a given confidence level. As in Section \ref{section no regression}, where we required the responses to be ordered such that their canonical parameters are increasing, the construction of this band is based on the assumption of knowing a ranking function that indicates the ordering of the true mean function for a.e.~$\bx \in \mathcal{X}$.

\label{known ranking}
\begin{ass} 
    \label{ranking pi}
    There exist a measurable ranking function $\pi : \mathcal{X} \to \R$ and a version of the true mean function $\mu^{*}_\pi : \X \to \R$, i.e., a regression function satisfying $$\mu^{*}_\pi(\bX) = \mu^*(\bX), \quad \p\textrm{-a.s.},$$
    such that for any two features $\bx, \bx' \in \mathcal{X}$,
    \begin{equation}
        \label{eq ranking pi}
        \pi(\bx) \leq \pi(\bx') \implies \mu^{*}_\pi(\bx) \leq \mu^{*}_\pi(\bx').
    \end{equation}
\end{ass}

Because the conditional mean in \eqref{true mean function} is only given $\p\textrm{-a.s.}$, we only require $\mu^{*}_\pi$ to align with $\mu^{*}$ $\p\textrm{-a.s.}$ In other words, the ranking function in \eqref{eq ranking pi} can be chosen such that it complies with the ranking of the true mean function $\mu^{*}$ for a.e.~feature $\bx \in \X$.

The existence of such a ranking function is clear as $\pi (\bx) = \mu^{*}(\bx)$ fulfills \eqref{eq ranking pi}. In fact, there are infinitely many ranking functions since, for example, any positive affine transformation of a ranking function is again a ranking function. The crucial point is that we assume to know at least one of these functions. 
Moreover, note that the above assumption is actually equivalent to saying that there exists a non-decreasing map $G : \R \to \R$ such that
        \begin{equation*}
            G(\pi(\bx)) = \mu_\pi^{*}(\bx),
        \end{equation*} 
for every $\bx \in \mathcal{X}$. Under Assumption \ref{ranking pi}, we aim at constructing a data-dependent calibration band
\begin{equation*}
    \left( L^\alpha_{\pi, (Y_i,\bX_i)_{i=1}^n} (\bx), U^\alpha_{\pi, (Y_i,\bX_i)_{i=1}^n} (\bx) \right)_{\bx \in \mathcal{X}},
\end{equation*}
i.e., a band depending on the realizations of the responses $(Y_i)_{i=1}^n$ and the features $(\bX_i)_{i=1}^n$ such that 
\begin{equation*}
    \p \Bigl(  L^\alpha_{\pi, (Y_i,\bX_i)_{i=1}^n} (\bx) \leq \mu_\pi ^{*}(\bx) \leq U^\alpha_{\pi, (Y_i,\bX_i)_{i=1}^n}(\bx)  \, \textrm{ for all } \bx \in \mathcal{X}  \Bigl) \, \geq \, 1- \alpha,
\end{equation*}
for any given confidence level $1-\alpha \in (0,1)$.  To do so, we further make the following assumption.
\begin{ass}
\label{ass reg cond prob}
The map
\begin{equation*}
    \Q : \Omega \times \mathcal{F} \to [0,1], \quad \Q (\omega, A) = \E[ \mathds{1}_A \, | \, \bX_1, \dots, \bX_n] (\omega),
\end{equation*}
is a regular conditional probability of $\p$, given the features $\bX_1, \dots, \bX_n$.
\end{ass}
This assumption fails to hold in general, and we refer to Section 3.2 of Rao--Swift \cite{Rao--Swift} for necessary conditions ensuring the existence of this regular conditional probability. Moreover, we emphasize that
Assumption \ref{ass reg cond prob} means that the map
\begin{equation*}
    A \in \mathcal{F} \mapsto \Q_{(\bx_i)_{i=1}^n} (A) = \p(A \, | \, \bX_1 = \bx_1, \dots, \bX_n = \bx_n),
\end{equation*}
is a probability measure on $(\Omega, \mathcal{F})$ for any realization $\bx_1, \dots, \bx_n$ of the features $\bX_1, \dots, \bX_n$. In particular, we have for any $A \in \mathcal{F}$,
\begin{equation}
    \label{tower prop integral}
    \begin{split}
    \p(A) &= \int \p(A \, | \, \bX_1 = \bx_1, \dots, \bX_n = \bx_n) \, d\p(\bx_1, \dots, \bx_n) \\
    &= \int \Q_{(\bx_i)_{i=1}^n} (A) \, d\p(\bx_1, \dots, \bx_n).
    \end{split}
\end{equation}

Denote by $\bx_1, \dots, \bx_n$ the observed features. The calibration band constructed in Theorem \ref{theo calibration bands} can now be extended to regression modeling under the probability measure $\Q_{(\bx_i)_{i=1}^n}$. For this, let $\tau_{\bx_1, \dots, \bx_n} : \{1, \dots, n\} \to \{1, \dots, n\}$ be any permutation on the indices such that for the given ranking function $\pi$, we have
\begin{equation}
    \label{permutation p}
    \pi(\bx_{\tau_{
\bx_1, \dots, \bx_n}(1)}) \leq \dots \leq \pi(\bx_{\tau_{\bx_1, \dots, \bx_n}(n)}).
\end{equation}
We point out that such a permutation always exists and depends on the given ranking function $\pi$. However, as for the ordering assumed in \eqref{ranking no regression}, the map $\tau_{
\bx_1, \dots, \bx_n}$ is in general not unique. We use it in order to rank the responses according to their conditional means and define the following weighted partial sums that depend on the responses $(Y_i)_{i=1}^n$, the features $(\bX_i)_{i=1}^n$ and the ranking function $\pi$ through the permutation $\tau_{\bx_1, \dots, \bx_n}$. We drop the subscript of the permutation function for convenience and write
\begin{equation*}
    Z_{j:k} = \frac{1}{v_{j:k}}\sum_{i=j}^k v_{\tau(i)} Y_{\tau(i)} ,
\end{equation*}
for $1 \leq j \leq k \leq n$, with aggregated volumes
 $$v_{j:k} = \sum_{i=j}^k v_{\tau(i)}.$$
 
\begin{theo}
    \label{cor calibration bands}
    Suppose that Assumptions \ref{assumption EDF}, \ref{ass effective domain}, \ref{ranking pi} and \ref{ass reg cond prob} hold. Let $(Y_i, \bX_i)_{i=1}^n$ be independent random variables such that $Y_i \, | \, \bX_i = \bx_i \sim \textnormal{EDF}(\theta(\bx_i), v_i, \varphi, \kappa(\cdot))$ for i.i.d.~features $(\bX_i)_{i=1}^n$ and given volumes $v_i > 0$.
    Moreover, let $\mathcal{J} \subseteq \{1, \dots, n\}^2$ be any set of ordered pairs. Then, for a.e.~realization of the features $\bx_1, \dots, \bx_n$ and any given confidence level $1-\alpha \in (0,1)$, we have
    \begin{equation}
    \label{prob bound Q_x}
    \Q_{(\bx_i)_{i=1}^n} \Bigl(  L^\alpha_{\pi, (Y_i, \bX_i)_{i=1}^n}(\bx) \leq \mu_\pi^{*}(\bx) \leq  U^\alpha_{\pi, (Y_i, \bX_i)_{i=1}^n}(\bx) \, \textnormal{ for all } \, \bx \in \mathcal{X} \Bigl) \geq 1- \alpha, 
    \end{equation}
where
\begin{equation*}
    L^\alpha_{\pi, (Y_i, \bX_i)_{i=1}^n}(\bx) = \sup_{(j,k) \in \mathcal{J}} \left[ l^{\delta}(Z_{j:k}, v_{j:k}, \varphi, \kappa(\cdot)) \, \mathds{1}_{\{\pi(\bx) \geq \pi(\bX_{\tau(k)})\}} + \inf_{\theta \in \interior{\Theta}} \kappa'(\theta) \, \mathds{1}_{\{\pi(\bx) < \pi(\bX_{\tau(k)})\}} \right],
\end{equation*}
and
\begin{equation*}
    U^\alpha_{\pi, (Y_i, \bX_i)_{i=1}^n}(\bx) = \inf_{(j,k) \in \mathcal{J}} \left[ u^{\delta}(Z_{j:k}, v_{j:k}, \varphi, \kappa(\cdot)) \, \mathds{1}_{\{\pi(\bx) \leq \pi(\bX_{\tau(j)})\}} + \sup_{\theta \in \interior{\Theta}} \kappa'(\theta) \, \mathds{1}_{\{\pi(\bx) > \pi(\bX_{\tau(j)})\}} \right],
\end{equation*}
for $\bx \in \mathcal{X}$ and $\delta = \alpha/(2|\mathcal{J}|)$.
\end{theo}

We emphasize that the construction of this calibration band on the mean is similar to Theorem \ref{theo calibration bands} and relies on the ranking function $\pi : \mathcal{X} \to \R$. Moreover, the statement in \eqref{prob bound Q_x} can be rewritten as
\begin{equation*}    
    \p\Bigl( \left. L^\alpha_{\pi, (Y_i, \bX_i)_{i=1}^n}(\bx) \leq \mu_\pi^{*}(\bx) \leq  U^\alpha_{\pi, (Y_i, \bX_i)_{i=1}^n}(\bx) \, \textnormal{ for all } \, \bx \in \mathcal{X} \,  \right| \, \bX_1 = \bx_1, \dots, \bX_n = \bx_n \Bigl) \geq 1- \alpha.
    \end{equation*}
 That is, given a.e.~realization of the features $(\bx_i)_{i=1}^n$, the probability that the realizations of the underlying responses lead to a calibration band being able to fully bound the mean function $\mu_\pi^{*} : \mathcal{X} \to \R$ exceeds $1-\alpha$. Due to \eqref{tower prop integral}, a corollary of Theorem \ref{cor calibration bands} is that
\begin{equation}
    \label{NEW}
    \p \Bigl(  L^\alpha_{\pi, (Y_i, \bX_i)_{i=1}^n}(\bx) \leq \mu_\pi^{*}(\bx) \leq  U^\alpha_{\pi, (Y_i, \bX_i)_{i=1}^n}(\bx) \, \textnormal{ for all } \, \bx \in \mathcal{X} \Bigl) \geq 1- \alpha, 
\end{equation}
for any given confidence level $1-\alpha \in (0,1)$. We emphasize, however, that the conditional probability bound in Theorem \ref{cor calibration bands} is stronger than this inequality as it holds for a.e.~fixed and known realization of the features $(\bX_i)_{i=1}^n$, i.e., when only the responses $(Y_i)_{i=1}^n$ are random. As the mean function $\mu_\pi^{*} : \mathcal{X} \to \R$ was assumed to be a version of the true mean function $\mu^{*} : \mathcal{X} \to \kappa'(\interior{\Theta})$, note that the calibration band constructed in this section provides a bound on the true mean function for a.e.~feature $\bx \in \X$ with a probability exceeding the confidence level $1-\alpha$. This leads to the construction of the statistical tests being introduced in the next section.

\section{Statistical testing of calibration and auto-calibration}

\label{sec tests}

\subsection{The auto-calibration property}

In regression modeling, the true mean function in \eqref{mean function} is approximated by a regression function that we denote by $\widehat{\mu}:\X \to \R$. This regression function is said to be calibrated if it matches the true mean function for a.e.~realization of the features $\bx \in \mathcal{X}$, i.e.,
\begin{equation}
    \label{cond calibration}
    \widehat{\mu}(\bx) = \mu^{*}(\bx), \quad \textrm{for a.e. } \bx \in \mathcal{X}.
\end{equation}
As the true mean function often exhibits a complex behaviour w.r.t.~the features $\bx \in \mathcal{X}$ and the mean estimation task is performed over a finite sample of (noisy) observations, it is in general impossible to aim for a calibrated regression function in practice. A related notion was introduced in the literature under the name of \textit{auto-calibration}. It is defined as follows, see, e.g., Krüger--Ziegel \cite{Kruger}.

\begin{defi}
    \label{def auto calibration}
    A regression function $\widehat{\mu}:\X \to \R$ is auto-calibrated for $(Y, \bX)$ if
    \begin{equation*}
        \widehat{\mu}(\bX) = \E[Y \, | \, \widehat{\mu}(\bX)], \quad \p\textrm{-a.s.}
    \end{equation*}
\end{defi}

Note that for any calibrated regression function $\widehat{\mu}:\X \to \R$, we have
\begin{equation}
    \label{cal implies auto-cal}
    \E[Y \, | \, \widehat{\mu}(\bX)] = \E[\E[Y \, | \, \bX] \, | \, \widehat{\mu}(\bX)] = \E[\mu^{*}(\bX) \, | \, \widehat{\mu}(\bX)] = \widehat{\mu}(\bX), \quad \p\textrm{-a.s.,}
\end{equation}
where in the first equality, we used that $\sigma(\widehat{\mu}(\bX)) \subseteq \sigma(\bX)$ and in the last equality, that the regression function $\widehat{\mu}:\X \to \R$ is calibrated. In other words, \eqref{cal implies auto-cal} shows that any calibrated regression function is auto-calibrated. The auto-calibration property means that by conditioning on a given mean estimate, the conditional expectation of the response matches this estimate. An auto-calibrated regression function thus guarantees that if we partition the feature space $\X$ into subsets according to the estimated regression values $\left(\widehat{\mu}(\bx)\right)_{\bx \in \X}$, the mean of all the responses within such a subset matches the estimated mean for this subset, leading to locally unbiased estimates. While this notion is weaker than calibration, it is of particular interest in several applications, where mean estimates within given specific groups have to be unbiased. For example, this is the case in insurance pricing, where the auto-calibration of a regression function is a minimal requirement, as it ensures that each cohort of individuals paying a certain price is self-financing; we refer to Wüthrich--Ziegel \cite{Wuethrich_Ziegel}.

\subsection{Statistical tests for calibration}

\label{stat tests for calibration}

The calibration band derived in Theorem \ref{cor calibration bands} can be used to construct a statistical test for calibration with confidence level $1-\alpha \in (0,1)$ as for any ranking function $\pi : \X \to \R$, \eqref{cond calibration} is equivalent to 
\begin{equation*}
    \widehat{\mu} (\bx) = \mu_\pi^{*}(\bx), \quad \textrm{for a.e. } \bx \in \mathcal{X}.
\end{equation*}
Moreover, the set
\begin{equation*}
    \Bigl\{  L^\alpha_{\pi, (Y_i, \bX_i)_{i=1}^n}(\bx) \leq \mu_\pi^{*}(\bx) \leq  U^\alpha_{\pi, (Y_i, \bX_i)_{i=1}^n}(\bx) \, \textnormal{ for all } \, \bx \in \mathcal{X}  \Bigl\}
\end{equation*}
is included in the set
\begin{equation*}
    \Bigl\{  L^\alpha_{\pi, (Y_i, \bX_i)_{i=1}^n}(\bx) \leq \widehat {\mu}(\bx) \leq  U^\alpha_{\pi, (Y_i, \bX_i)_{i=1}^n}(\bx) \, \textnormal{ for a.e.~} \, \bx \in \mathcal{X}  \Bigl\},
\end{equation*}
when $\widehat{\mu} : \X \to \R$ is calibrated and the former holds with probability greater or equal than $1-\alpha$, see \eqref{NEW}. Therefore, by writing $(y_i, \bx_i)_{i=1}^n$ for the observed realizations of the responses and features and $$\left(L^\alpha_{\pi, (y_i, \bx_i)_{i=1}^n}(\bx), U^\alpha_{\pi, (y_i, \bx_i)_{i=1}^n}(\bx)\right)_{\bx \in \X},$$ for the resulting calibration band, we reject the null-hypothesis
\begin{equation}
    \label{classical test}
    \mathbb{H}_0 : \widehat{\mu}(\bx) = \mu^{*}(\bx)  \textrm{ for a.e.~} \bx \in \mathcal{X},
\end{equation}
with confidence level $1-\alpha$ whenever the set
\begin{equation}
    \label{condition test 1}
    \X^{out} = \left\{ \bx \in \X \, \left| \, \widehat{\mu}(\bx) \notin \left[ L^\alpha_{\pi, (y_i, \bx_i)_{i=1}^n}(\bx),  U^\alpha_{\pi, (y_i, \bx_i)_{i=1}^n}(\bx) \right] \right. \right\},
\end{equation}
satisfies $\p(\bX \in \X^{out}) > 0$. That is, we reject the null-hypothesis of calibration of a model whenever the set of features for which the mean estimates fall outside the band has non-zero probability under the distribution of the features $\bX \in \X$.
\medskip

We propose a procedure in order to graphically determine the decision induced by this statistical test. First, the calibration band can be plotted against the ranking function, which results in two non-decreasing step functions delimiting the band. Then, the pairs $(\pi(\bx), \widehat{\mu}(\bx))_{\bx \in \X}$ can be drawn in the same plot for all features $\bx \in \X$, and we reject the null-hypothesis $\mathbb{H}_0$ whenever the set of pairs falling outside the calibration band is large enough, i.e., whenever it corresponds to a set of features being larger than a nullset. We call such a plot a {\it calibration plot}, and we emphasize that as the distribution of $\bX$ is unknown in practice, the decision of rejecting or not the above null-hypothesis requires some assumptions on the support of the random variable $\widehat{\mu}(\bX)$, we come back to this in Section \ref{numerical examples}, below.

We also point put that the decision of the statistical test depends on the ranking function through the constructed calibration band. In practice, there are a few cases where a ranking function is known, and in such cases, statistical techniques under order restrictions could be used for mean estimation; we refer to Barlow et al.~\cite{Barlow} and Robertson et al.~\cite{Robertson}. Note that one of these techniques is isotonic regression, which has the nice property to lead to empirically auto-calibrated regression functions. Most of the time, however, we do not have access to any ranking function giving the ordering of the true mean function over the feature space $\mathcal{X}$. In such cases, the ranking function needs to be approximated. We start our discussion from a related work, where Wüthrich--Ziegel \cite{Wuethrich_Ziegel} propose a method to restore the auto-calibration of a given regression function $\widehat{\mu}:\X \to \R$. For this, they perform an isotonic regression by using the regression function itself as a ranking function and they call their method the \textit{isotonic recalibration} step because in a first step, a regression function is estimated and under the assumption that it provides the correct ranking, this ranking is lifted to be (empirically) auto-calibrated in a second (isotonic recalibration) step. We make the same choice here and take $\widehat{\mu} : \X \to \R$ as a ranking function.
\begin{ass}
    \label{ranking mu hat}
    The regression function $\widehat{\mu}:\X \to \R$ satisfies that there exists a version of the true mean function $\mu^{*}_{\widehat{\mu}} : \X \to \R$, i.e., a regression function satisfying $$\mu^{*}_{\widehat{\mu}}(\bX) = \mu^*(\bX), \quad \p\textrm{-a.s.},$$
    such that for any two features $\bx, \bx' \in \mathcal{X}$,
    \begin{equation*}
        \widehat{\mu}(\bx) \leq \widehat{\mu}(\bx') \implies \mu^{*}_{\widehat{\mu}}(\bx) \leq \mu^{*}_{\widehat{\mu}}(\bx').
    \end{equation*}
\end{ass} 
This choice was also implicitly made by Dimitriadis et al.~\cite{Dimitriadis} for the binary case and, actually, we emphasize that Assumption \ref{ranking mu hat} is nested in the null-hypothesis \eqref{classical test}, i.e., it holds under $\mathbb{H}_0$.
This implies, in particular, that the above statistical test can always be applied by taking the regression function itself as a ranking function.

Finally, as pointed out by Dimitriadis et al.~\cite{Dimitriadis}, an opposite statistical test can be constructed where the calibration property now lies in the alternative. This test allows one to quantify deviations from calibration as its null-hypothesis reads as
\begin{equation}
    \label{opposite test}
    \mathbb{H}_0 : \left| \widehat{\mu}(\bx)- \mu^{*}(\bx) \right| > \varepsilon\textrm{ for a.e.~} \bx \in \mathcal{X},
\end{equation}
for some given $\varepsilon > 0$. This hypothesis is rejected with confidence level $1-\alpha$, whenever the set
\begin{equation}
    \label{condition test 2}
    \mathcal{X}_{\varepsilon}^{in}  = \left \{\bx \in \mathcal{X} : \Bigl[ L^\alpha_{\pi, (y_i, \bx_i)_{i=1}^n}(\bx),  U^\alpha_{\pi, (y_i, \bx_i)_{i=1}^n}(\bx)\Bigl] \, \subseteq \, \Bigl[ \widehat{\mu}(\bx)-\varepsilon, \widehat{\mu}(\bx)+\varepsilon  \Bigl]  \right\}
\end{equation}
satisfies $\p \left(\bX \in \mathcal{X}_{\varepsilon}^{in} \right) > 0$. As above, note that the evaluation of this test can be done graphically by plotting the calibration band and the estimated means against the ranking function. Moreover, the previous discussion about the choice of the ranking function equivalently applies here. 

\subsection{Statistical tests for auto-calibration}

\label{Statistical tests for auto-calibration}

In order to construct statistical tests for the auto-calibration property of a given regression function $\widehat{\mu}:\X \to \R$, we assume that Assumption \ref{ranking mu hat} holds. That is, the regression function manages to correctly provide the ordering of the true mean function. 
Interestingly, we can show that under this assumption, calibration is equivalent to auto-calibration, see Proposition 5.1 in Denuit--Trufin \cite{Denuit--Trufin}.
This means that any auto-calibrated regression function managing to correctly express the ranking of the means over the feature space is actually equal to the true mean function for a.e.~feature $\bx \in \mathcal{X}$. Thus, the tests \eqref{classical test} and \eqref{opposite test} do not only provide a test for calibration but also for auto-calibration under Assumption \ref{ranking mu hat}.
The first test consists in rejecting the auto-calibration of a regression function $\widehat{\mu}:\X \to \R$ with null-hypothesis
\begin{equation*}
    \mathbb{H}_0 : \E[Y \, | \, \widehat{\mu}(\bX) = \widehat{\mu}(\bx)] = \widehat{\mu}(\bx) \textrm{ for a.e.~} \bx \in \mathcal{X},
\end{equation*}
whenever the set $\X^{out}$ in \eqref{condition test 1} satisfies $\p(\bX \in \X^{out}) > 0$. The second test consists in rejecting the null-hypothesis
\begin{equation*}
    \mathbb{H}_0 : \Bigl| \E[Y \, | \, \widehat{\mu}(\bX) = \widehat{\mu}(\bx)] - \widehat{\mu}(\bx) \Bigl| > \varepsilon \textrm{ for a.e.~} \bx \in \mathcal{X},
\end{equation*}
whenever the set $\X_{\varepsilon}^{in}$ in \eqref{condition test 2} satisfies $\p(\bX \in \X_\varepsilon^{in}) > 0$. Finally, we conclude this section by highlighting that the derived statistical tests for calibration and auto-calibration also apply to the framework of Section \ref{section no regression}. Indeed, by choosing as feature space $\mathcal{X} = \{1, \dots, n\}$ and an appropriate ranking function $\pi : \mathcal{X} \to \{1, \dots, n\}$ such that $\theta(\bx_i) =\theta_{\pi(\bx_i)}$ holds for all $1 \leq i \leq n$, the mean estimation task in Section \ref{section no regression} can be expressed as a regression modeling problem.

\subsection{Impact of the dispersion parameter}

\label{sec dispersion}
The calibration band on the mean derived in Theorem \ref{cor calibration bands} holds for EDF responses, for which the dispersion parameter $\varphi$ and the cumulant function $\kappa$ are known and fixed. The cumulant function uniquely determines the distribution of the responses, while the dispersion parameter characterizes their variances. Indeed, the variance of an EDF random variable $Y \sim \textrm{EDF}(\theta, v, \varphi, \kappa(\cdot))$ for $\theta \in \interior{\Theta}$ is given by
\begin{equation*}
    \textrm{Var} (Y) = \frac{\varphi}{v} \kappa''(\theta),
\end{equation*}
and it can equivalently be expressed in terms of the mean $\mu = \E[Y]$ using the canonical link in \eqref{canonical link} through
\begin{equation*}
    \textrm{Var} (Y) = \frac{\varphi}{v} \kappa''(h(\mu)) = \frac{\varphi}{v} V(\mu),
\end{equation*}
where $V = \kappa'' \circ h$ is the {\it variance function} of the chosen distribution within the EDF. The larger the dispersion parameter $\varphi$ is, the larger the variance of the underlying response will be. Thus, the constructed calibration bands depend on the value of $\varphi$, and this dependence actually lies in the lower and upper bounds in \eqref{lower bound delta} and \eqref{upper bound delta}.

Using a suitable dispersion estimate is thus crucial in order to construct the above statistical tests as, in practice, the dispersion parameter is often unknown. There are various methods for estimating this parameter, see Section 5.3.1 of Wüthrich--Merz \cite{WM2023}. We introduce one of these methods here, which consists in computing the {\it Pearson's estimate} that is given by
\begin{equation*}
    \widehat{\varphi}^P = \frac{1}{n-q} \sum_{i=1}^n \frac{(Y_i-\widehat{\mu}(\bX_i))^2}{V(\widehat{\mu}(\bX_i))/v_i},
\end{equation*}
where $q$ denotes the number of unknown parameters that are estimated in order to derive the regression function $\widehat{\mu}:\X \to \R$. Note that when this regression function is obtained from a well-specified generalized linear model (GLM), Pearson's estimate has the advantage of providing a consistent estimator for $\varphi$, see Section 8.3.6 of McCullagh--Nelder \cite{Nelder--Wedderburn}.

\section{Numerical examples}

This section provides numerical examples where we construct calibration bands on the mean of given responses. Our goal is first to highlight the impact of different factors on the resulting calibration bands, as the choice of the confidence level and the set of ordered pairs, or the influence of binning observations. Then, we study a lime trees real dataset and construct a calibration band on the mean of the foliage biomass of the trees. We show that for this small dataset, the statistical test for auto-calibration introduced in Section \ref{Statistical tests for auto-calibration} manages to detect violations of auto-calibration in contrast to the test of Denuit et al.~\cite{Denuit2}. After that, we discuss the impact of estimating the dispersion parameter on the resulting calibration band for simulated inverse Gaussian responses. We then consider a popular French motor third party liability real dataset and construct a calibration band on the claim frequency of the insured drivers. We show that for this large dataset of more than half a million insurance policies, the resulting band allows us to detect violations of calibration and that the isotonic recalibration step proposed by Wüthrich--Ziegel \cite{Wuethrich_Ziegel} addresses this issue.  Finally, we discuss the power of the statistical test for calibration presented in Section \ref{stat tests for calibration} by considering the same simulated example as Wüthrich \cite{Wüthrich}.

\label{numerical examples}

\subsection{Example 1 : calibration bands for simulated normal responses}

In this first example, we consider simulated normal responses and aim at assessing the calibration of mean estimates that are obtained using the isotonic estimator introduced in \eqref{isotonic estimator}. For this, we sample $n = 2000$ independent normal random variables $Y_i \sim \mathcal{N}(\mu_i, \sigma_i)$, where the means are equally spaced over the range $[1500, 2500]$, i.e., 
\begin{equation*}
    \mu_i = \E[Y_i] = 1500+ \frac{i-1}{n-1} \cdot 1000, \quad \textrm{for }1 \leq i \leq n,
\end{equation*}
and where the standard deviations are chosen to satisfy $\sigma_i = 0.5 \mu_i$. This choice of the parameters $(\mu_i, \sigma_i)_{i=1}^{n}$ leads to a coefficient of variation being constant for all responses as we have
\begin{equation*}
    \textrm{Vco}(Y_i) = \frac{\sqrt{\textrm{Var}(Y_i)}}{\E[Y_i]} = \frac{\sigma_i}{\mu_i} = \frac{1}{2} , \quad \textrm{for }1 \leq i \leq n.
\end{equation*}

\begin{figure}[htb!]
    \centering
    \includegraphics[width = 0.7\linewidth]{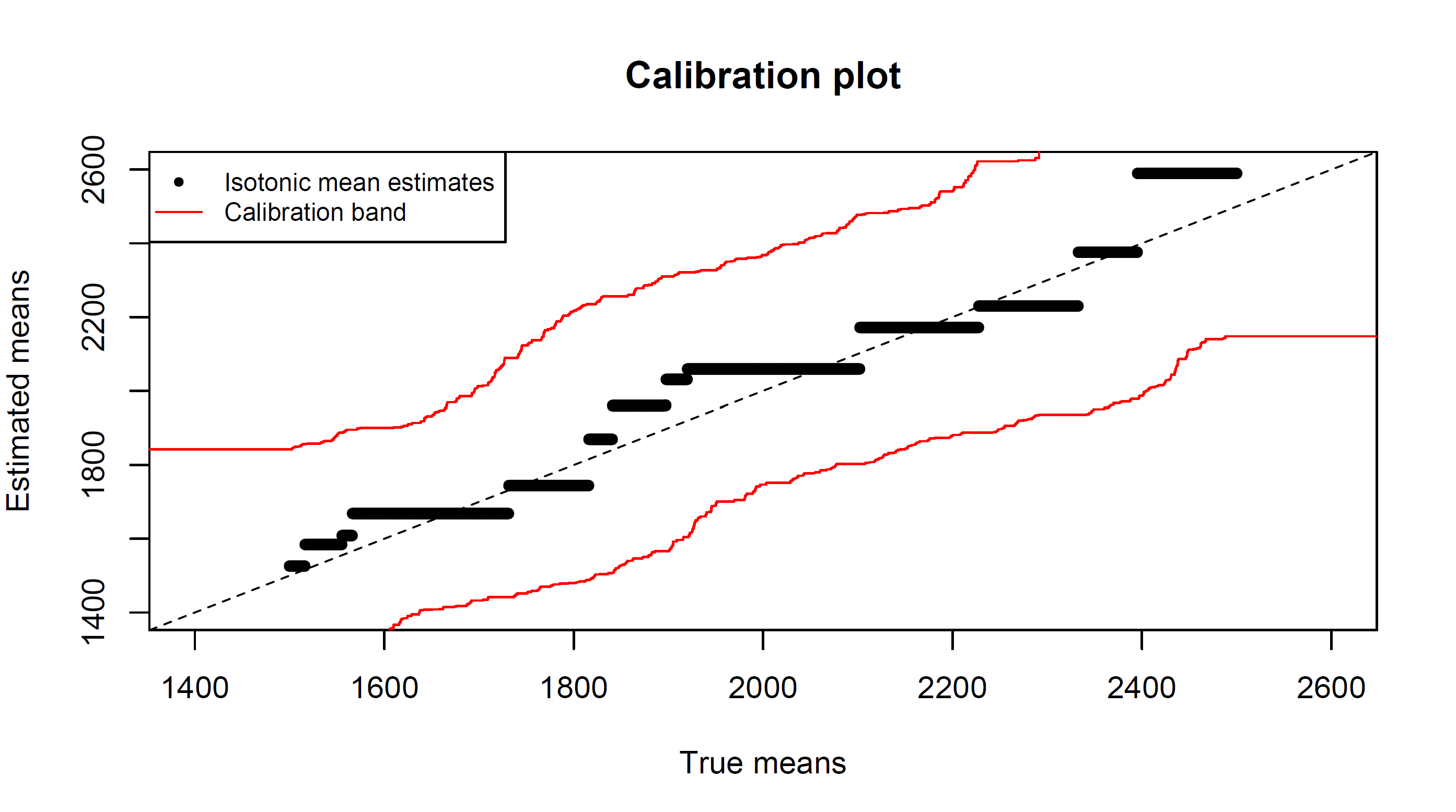}
    \caption{Calibration plot of the isotonic mean estimates of independent normal responses $(Y_i)_{i=1}^{n}$. The calibration band is plotted in red, whereas the mean estimates are drawn in black.}
    \label{Fig:Plot1}
\end{figure}
Then, we estimate the means of the responses using the simulated responses $\bY = (Y_1, \dots, Y_n)^\top$ based on the isotonic mean estimator $\widehat{\bmu}^{Iso}(\bY, \bv)$, i.e., we set
\begin{equation*}
    \widehat{\mu}_i =   \widehat{\mu}^{Iso}(\bY, \bv)_i , \quad \textrm{for }1 \leq i \leq n,
\end{equation*}
where $\bv = (1/\sigma_1^2, \dots, 1/\sigma_n^2)^\top$, we refer to Remark \ref{remark} for the choice of the volumes $\bv$. In order to assess the calibration of the mean estimates $(\widehat{\mu}_i)_{i=1}^n$, we construct a full calibration band on the mean of the above responses using the ordering of their true means and a confidence level of $1-\alpha = 0.95$. The resulting calibration plot is provided in Figure \ref{Fig:Plot1}. As all the mean estimates (black dots) lie within the band (red lines), the conclusion of the test in \eqref{classical test} is to not reject the calibration assumption of the isotonic mean estimator $\widehat{\bmu}^{Iso}(\bY, \bv)$ in this example. Note that the decision of the performed test depends on the confidence level and the set of ordered pairs used to construct the calibration band. We show the impact of these factors on the band below and we additionally look at the effect of binning observations.

\subsubsection{Sensitivity with respect to the confidence level}

\begin{figure}[htb!]
    \centering
    \includegraphics[width = 0.7\linewidth]{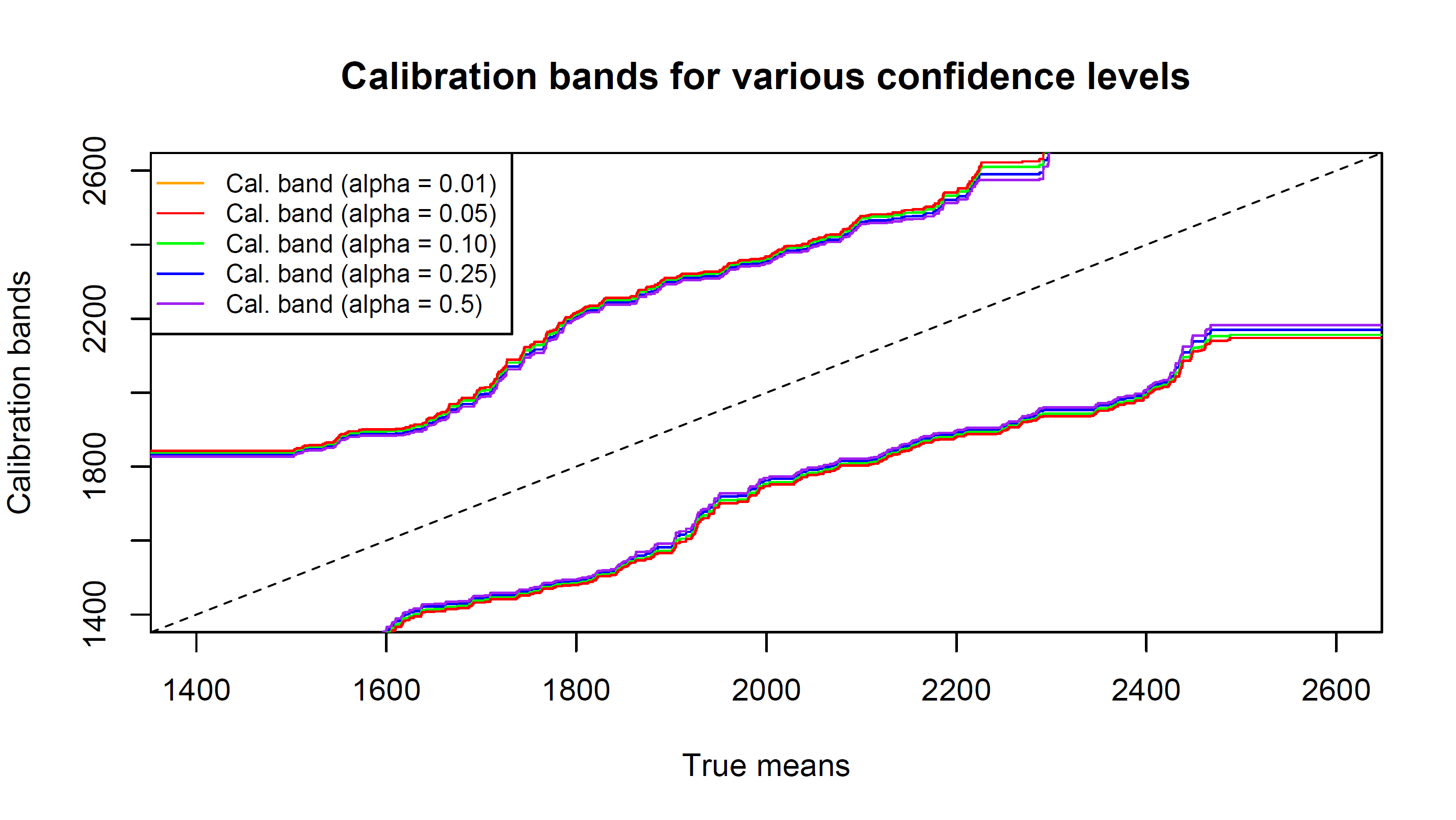}
    \caption{Calibration bands on the mean of independent normal responses $(Y_i)_{i=1}^{n}$ for various confidence levels.}
    \label{Fig:Plot_1_alpha_}
\end{figure}

The width of the calibration band depends on the chosen confidence level. We evaluate this impact in this example by providing full calibration bands on the mean of the above simulated responses for various confidence levels $1-\alpha \in \{0.99, 0.95, 0.9, 0.75, 0.5\}$ in Figure \ref{Fig:Plot_1_alpha_}.

As expected, the calibration bands get narrower as the value of $\alpha$ increases. However, we point out that the value of the confidence level seems not to lead to significant impacts on the width of the calibration band. That is, the constructed bands are not very sensitive to the confidence level, which indicates that the union bound inequality in \eqref{inequality discrete case} is not very sharp in this example.

\subsubsection{Impact of the chosen set of ordered pairs}

\label{sec set ordered pairs}

\begin{figure}[htb!]
    \centering
    \includegraphics[width=0.7\linewidth]{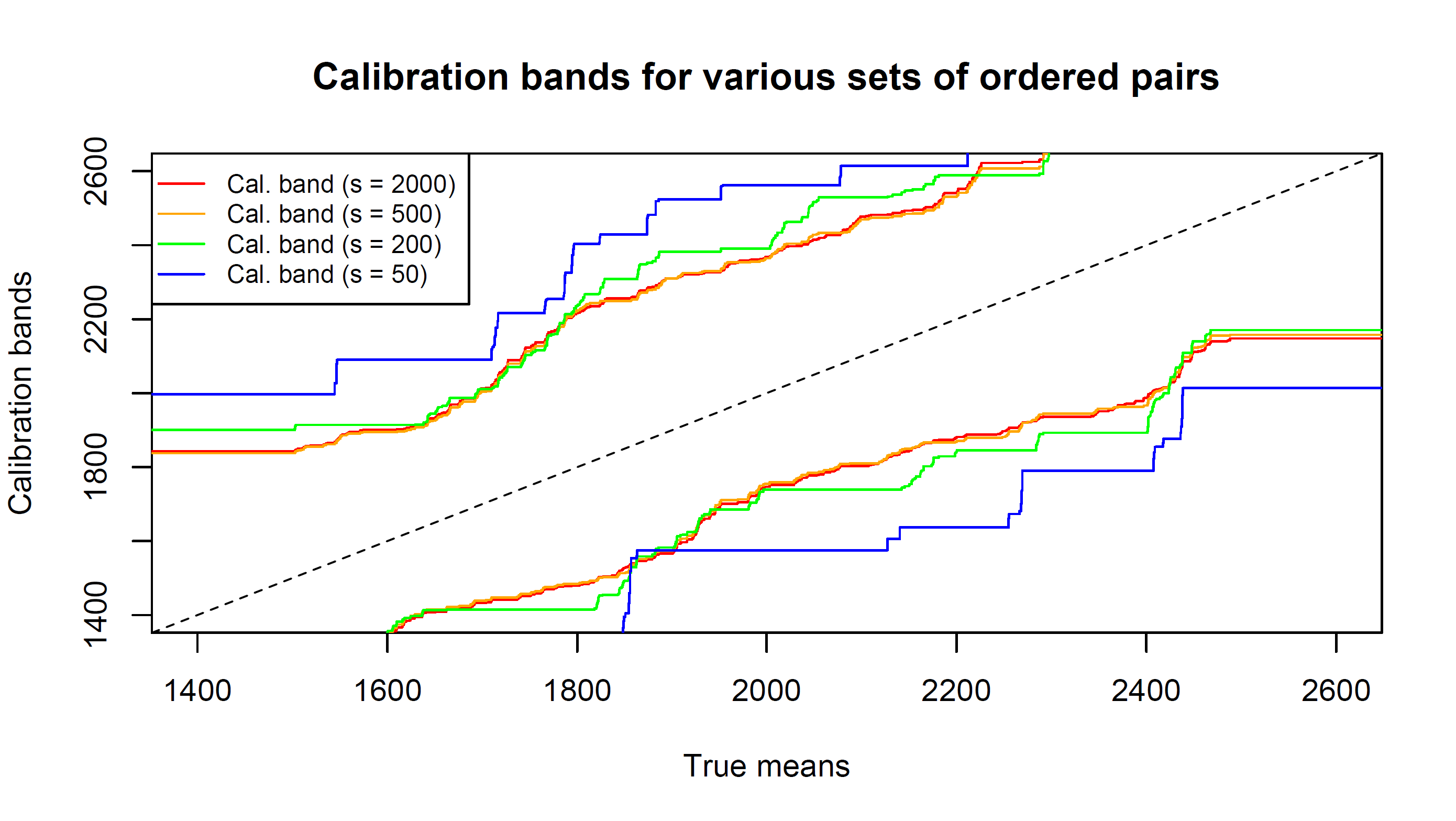}
    \caption{Calibration bands on the mean of independent normal responses $(Y_i)_{i=1}^{n}$ that are constructed using the set $\mathcal{J}^{nbh}_{s}$ for various sizes $s$.}
    \label{Fig:Normal_set_ordered_pairs_nbh}
\end{figure}

We now study the impact of constructing calibration bands with sets of ordered pairs that are different from $\mathcal{J}^{full}$. For this, we use again the above simulated observations and fix the confidence level at $1-\alpha = 0.05$. As discussed in Section \ref{section computational time}, note that the choice of a smaller set of ordered pairs enables us to reduce the computational time required to construct the band. We first consider a set of nearest neighbours (nbh) given by
\begin{equation*}
    \mathcal{J}^{nbh}_{s} = \left\{(j,k) \in  \left. \mathcal{J}^{full} \, \right| \,  k-j \leq s \right\},
\end{equation*}
for $s \in \N$. The use of such a set of ordered pairs can be justified by the intuition that the means of the weighted partial sums $Z_{j:k}$ are too far from $\mu_j$ and $\mu_k$ when the difference $k-j$ is large, implying the bounds in Proposition \ref{last corollary} not to be very sharp. We construct calibration bands using the set $\mathcal{J}^{nbh}_{s}$ for different sizes $s \in \{50, 200, 500, 2000\}$ and plot these bands in Figure \ref{Fig:Normal_set_ordered_pairs_nbh}. Moreover, the computational time required to construct the bands is shown in Table \ref{Tab:Normal_set_ordered_pairs_nbh}. 

\begin{table}[htb!]
    \centering
    \begin{tabular}{l c c c c}
      \toprule
      $s$ & $2000$ & $500$ & $200$ & $50$\\
      \midrule
      Time (s) & $119.81$ & $25.16$ & $9.86$ & $2.72$\\
      \bottomrule
    \end{tabular}
    \caption{Time (seconds) required to construct the calibration bands on the mean of independent normal responses $(Y_i)_{i=1}^{n}$ using the set $\mathcal{J}^{nbh}_{s}$ for various sizes $s$.}
    \label{Tab:Normal_set_ordered_pairs_nbh}
\end{table}
Interestingly, it seems that although small sizes $s$ lead to small sets of ordered pairs $ \mathcal{J}^{nbh}_{s}$, the constructed calibration bands seem close to each other for $s \in \{200, 500, 2000\}$, whereas the band is significantly wider for the case $s = 50$. This might be due to the underlying aggregated volumes $v_{j:k}$ that fail to be large enough in order to obtain a suitable band in the latter case. 

\begin{figure}[htb!]
    \centering
    \includegraphics[width=0.7\linewidth]{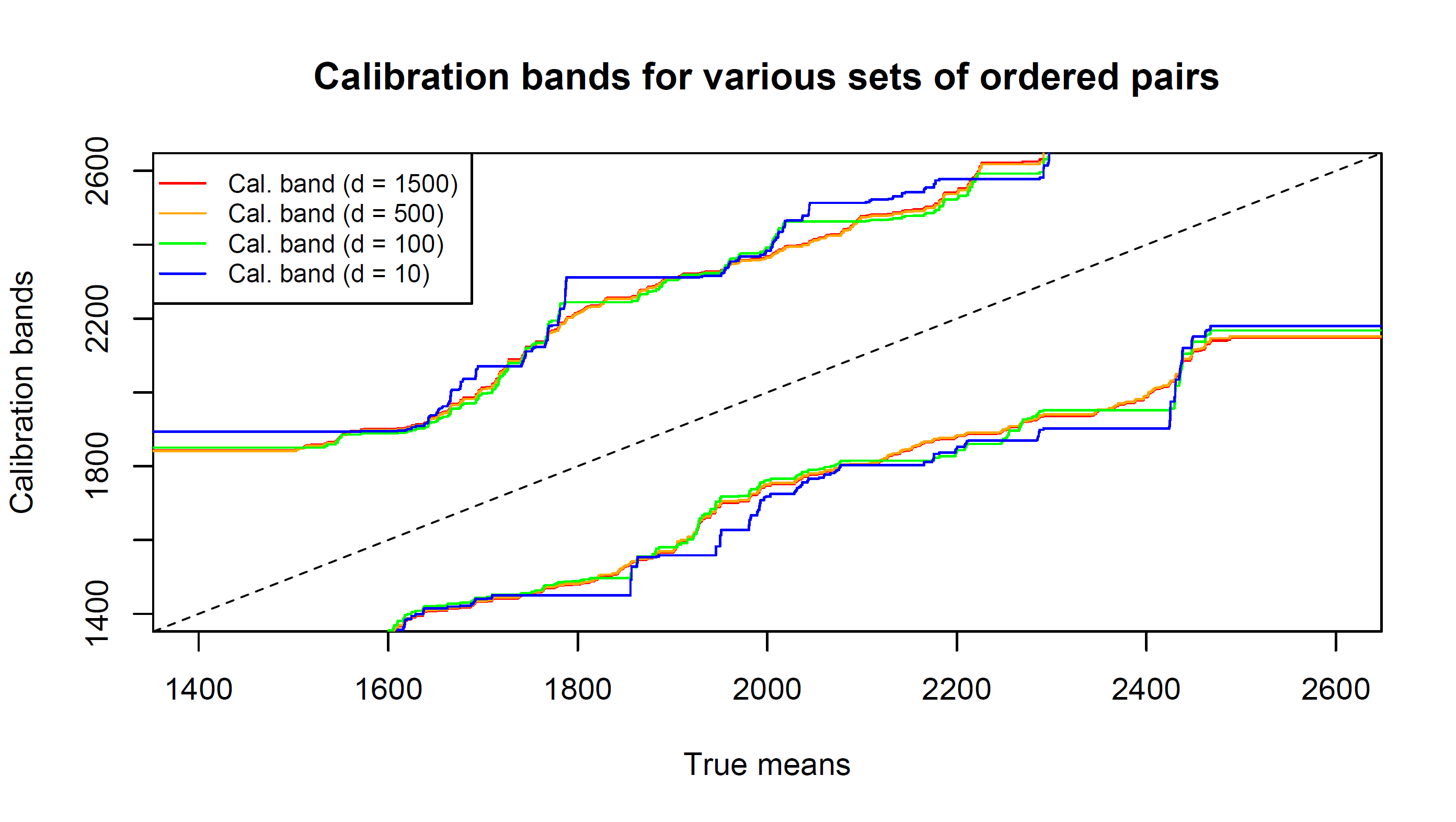}
    \caption{Calibration bands on the mean of independent normal responses $(Y_i)_{i=1}^{n}$ that are constructed using the set $\mathcal{J}^{dist}_{d}$ for various distances $d$.}
    \label{Fig:Normal_set_ordered_pairs_estimates_approx}
\end{figure}

Next, we consider another set of ordered pairs based on the distance (dist) between available mean estimates $(\widehat{\mu}_i)_{i=1}^n$
\begin{equation}
    \label{J dist}
    \mathcal{J}^{dist}_{d} = \left\{(j,k) \in  \left. \mathcal{J}^{full} \, \right| \, \left| \widehat{\mu}_i - \widehat{\mu}_j\right| \leq d \right\},
\end{equation}
for $d \in \R$. The idea behind such a choice is to only take pairs into account, for which we believe that the means of the responses are close to each other. By using the above isotonic mean estimates in order to define the set $\mathcal{J}^{dist}_{d}$, i.e., by setting $\widehat{\mu}_i = \widehat{\mu}^{Iso}(\bY, \bv)_i$ for $1 \leq i \leq n$ in \eqref{J dist}, we construct calibration bands for various distances $d \in \{10, 100, 500, 1500\}$ and plot them in Figure \ref{Fig:Normal_set_ordered_pairs_estimates_approx}. The computational time required to construct these bands is provided in Table \ref{Tab:Normal_set_ordered_pairs_estimates_approx}. This time, we notice that all the calibration bands have a similar width for almost all distances $d$, except for the case $d = 10$ where the band seems to be wider at some specific locations. That is, for both restricted sets of ordered pairs we consider in this example, our results show that constructing bands using smaller sets of ordered pairs than $\mathcal{J}^{full}$ is computationally more efficient, and these smaller sets lead to suitable calibration bands as long as their size is not too small.

\begin{table}[htb!]
    \centering
    \begin{tabular}{l c c c c}
      \toprule
      $d$ & $1500$ & $500$ & $100$ & $10$\\
      \midrule
      Time (s) & $124.62$ & $58.98$ & $12.62$ & $6.54$ \\
      \bottomrule
    \end{tabular}
    \caption{Time (seconds) required to construct the calibration bands on the mean of independent normal responses $(Y_i)_{i=1}^{n}$ using the set $\mathcal{J}^{dist}_{d}$ for various distances $d$.}
    \label{Tab:Normal_set_ordered_pairs_estimates_approx}
\end{table}

\subsubsection{Impact of binning observations}

\label{sec binning}
\begin{figure}[htb!]
    \centering
    \includegraphics[width = 0.7\linewidth]{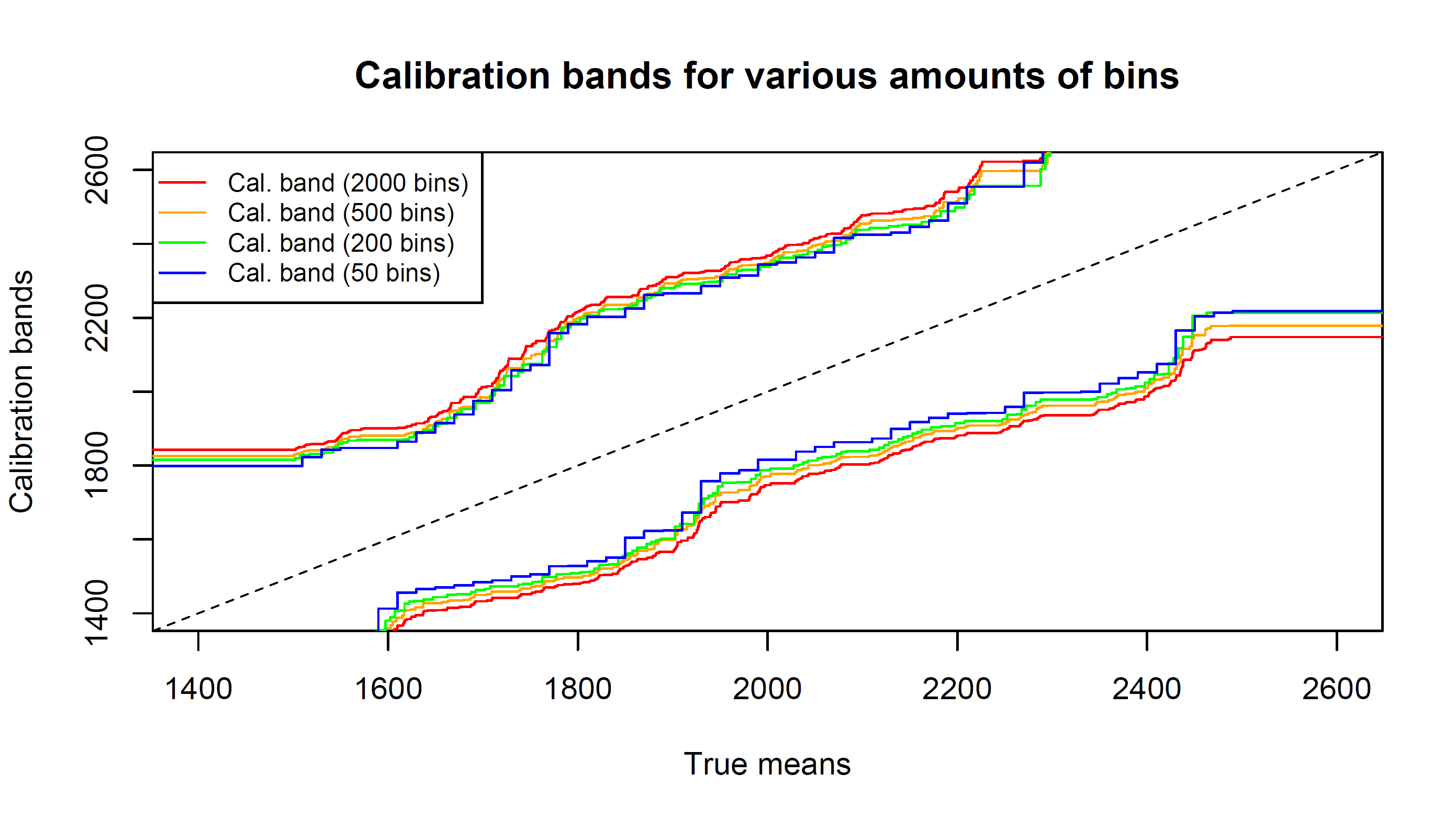}
    \caption{Calibration bands on the mean constructed by binning independent normal responses $(Y_i)_{i=1}^{n}$ for different amounts of bins.}
    \label{Fig:Plot_2_Bin}
\end{figure}

Another method for reducing the computational time needed to construct the calibration bands consists in binning observations, we refer to Section \ref{section computational time}. To understand the impact of this method on the resulting bands in this example, we set the confidence level to $1-\alpha = 0.95$, choose $\mathcal{J}^{full}$ as the set of ordered pairs, and use the same independent normal observations $(y_i)_{i=1}^{n}$ as above. Moreover, we create $L$ equally sized bins in order to define new observations $(\tilde{y}_l)_{l=1}^L$ that satisfy
\begin{equation}
    \label{binned normal}
    \tilde{y}_l = \sum_{k=n(l-1)/L+1}^{nl/L} \frac{v_k y_k}{v_{n(l-1)/L+1:nl/L}},
\end{equation}
with
\begin{equation*}
    v_i = \frac{1}{\sigma_i^2}, \quad \textrm{for } 1 \leq i \leq n.
\end{equation*}
That is, the new observations are weighted sums of the original ones, with weights being equal to the volumes of the original observations, see Remark \ref{remark}. Under the assumption that the means of the responses within a given bin are equal, note that this weighting implies that the binned responses belong to the EDF due to the convolution formula in Corollary 2.15 of Wüthrich--Merz \cite{WM2023}. In general, we emphasize that the new observations are not realizations of EDF random variables as this assumption might be violated. In this example, however, this assumption is not needed as binned normal responses are always normally distributed, regardless of the chosen weights. Nonetheless, we still choose the volumes of the original observations as weights in \eqref{binned normal} and construct calibration bands using as rankings the true means of the binned responses for $L \in \{50, 200, 500, 2000\}$. The resulting bands are provided in Figure \ref{Fig:Plot_2_Bin} and the computational times required to construct them are given in Table \ref{Tab:Plot_2_Bin}. 

\begin{table}[htb!]
    \centering
    \begin{tabular}{l c c c c}
      \toprule
      $l$ & $2000$ & $500$ & $200$ &  $50$\\
      \midrule
      Time (s) & $122.66$ & $4.31$ & $1.15$ &  $0.63$\\
      \bottomrule
    \end{tabular}
    \caption{Time (seconds) required to construct the calibration bands on the mean of  binned independent normal responses $(Y_i)_{i=1}^{n}$ for different amounts of bins.}
    \label{Tab:Plot_2_Bin}
\end{table}
We notice that all bin sizes lead to pretty similar calibration bands in this example and actually, we can even observe in Figure \ref{Fig:Plot_2_Bin} that the bands get narrower the smaller the number of bins is. As discussed in Section \ref{section computational time}, this might be a consequence of having a small number of observations, implying the set of ordered pairs $\mathcal{J}^{full}$ to be small. In fact, the ratio
\begin{equation*}
    \frac{\Phi^{-1}\left(\frac{0.05}{2000^2+2
    000}\right)}{\Phi^{-1}\left(\frac{0.05}{50^2+50}\right)}  = 1.355
\end{equation*}
hints that the band constructed using $2000$ bins should be approximately 1.355 wider than the band constructed using only 50 bins, see \eqref{L_i normal} and \eqref{U_i normal}. However, this is is not the case in Figure \ref{Fig:Plot_2_Bin} due to the role played by the weighted partial sums $Z_{j:k}$ and the aggregated volumes $v_{j:k}$ that are used to construct the calibration bands. 


Together with Section \ref{sec set ordered pairs}, this section shows that binning observations or choosing suitable sets of ordered pairs can be an interesting technique to reduce the computational time required to construct the calibration bands. In this example, it seems that the chosen sets of ordered pairs enable to reduce the running time of the construction at the cost of having slightly wider calibration bands, whereas binning observations leads to even narrower bands. More generally, we point out that the choice of the set of ordered pairs or the size of the bins has to be carefully made in practice as the true means of the responses are unknown and might not be evenly distributed over the range of interest. Moreover, note that while the inequality \eqref{inequality discrete case} holds for any chosen set of ordered pairs, it does not hold anymore when the observations are binned since the distribution of the underlying binned responses is unknown. The latter method has the advantage of leading to low computational times, while using all the observations and constructing the band with large aggregated volumes.
On the contrary, by reducing the number of elements in the set of ordered pairs, the resulting band has to be constructed using small aggregated volumes, and this is actually the reason why the bands become too wide for small sets of ordered pairs in Figures \ref{Fig:Normal_set_ordered_pairs_nbh} and \ref{Fig:Normal_set_ordered_pairs_estimates_approx}. Therefore, we recommend to use the binning method for large datasets. We will follow this choice in Section \ref{MTPL sec} where we consider a portfolio of more than half a million insurance policies.

\subsection{Example 2 : calibration bands for a small sample-sized real dataset}

\label{sec ex 2}

We consider in this example a lime trees real dataset\footnote{The dataset can be downloaded by running       {\tt library(GLMsData); data(lime)} in {\sf R}.} available from the
{\sf R} \cite{R Core Team} package {\tt GLMsData} hosted by Dunn--Smyth \cite{Dunn--Smyth}. This dataset contains measurements from $n= 385$ small-leaved lime trees growing in Russia such as the foliage biomass (in kg), the tree diameter (in cm) and the origin of the tree ({\tt Coppice, Natural} or {\tt Planted}). Our goal is to estimate the foliage biomass using the tree diameter and the origin of the tree. As in Examples 11.4 and 11.7 in Dunn--Smyth \cite{Dunn--Smyth2}, we fit a gamma and an inverse Gaussian GLM with a logarithmic link using a continuous feature providing the logarithm of the tree diameter and a categorical feature providing the origin of the tree. In both GLMs, we further include an interaction term between the two features and we denote the resulting regression functions by $\hat{\mu}^{gamma} : \X \to (0,\infty)$ and $\hat{\mu}^{IG} : \X \to (0,\infty)$, where 
\begin{equation}
\label{feature space ex 2}
\X = (0, \infty) \times \{\textrm{Coppice, Natural, Planted} \},
\end{equation}
stands for the feature space. By construction, both models have the same number of parameters and the volumes are $\bv = (1, \dots, 1)^\top \in \R^n$, because we did not introduce any weights in the fitting procedure. Based on AIC, we clearly give preference to the gamma GLM, see Table \ref{Tab:AIC}.
\begin{table}[htb!]
    \centering
    \begin{tabular}{l c c }
      \toprule
      & gamma GLM & inverse Gaussian GLM \\ 
      \midrule
      AIC & $750.33$ & $1089.50$ \\ 
      \bottomrule
    \end{tabular}
    \caption{AIC of the gamma and inverse Gaussian GLMs.}
    \label{Tab:AIC}
\end{table}

\begin{figure}[htb!]
    \centering
    \includegraphics[width = 0.7\linewidth]{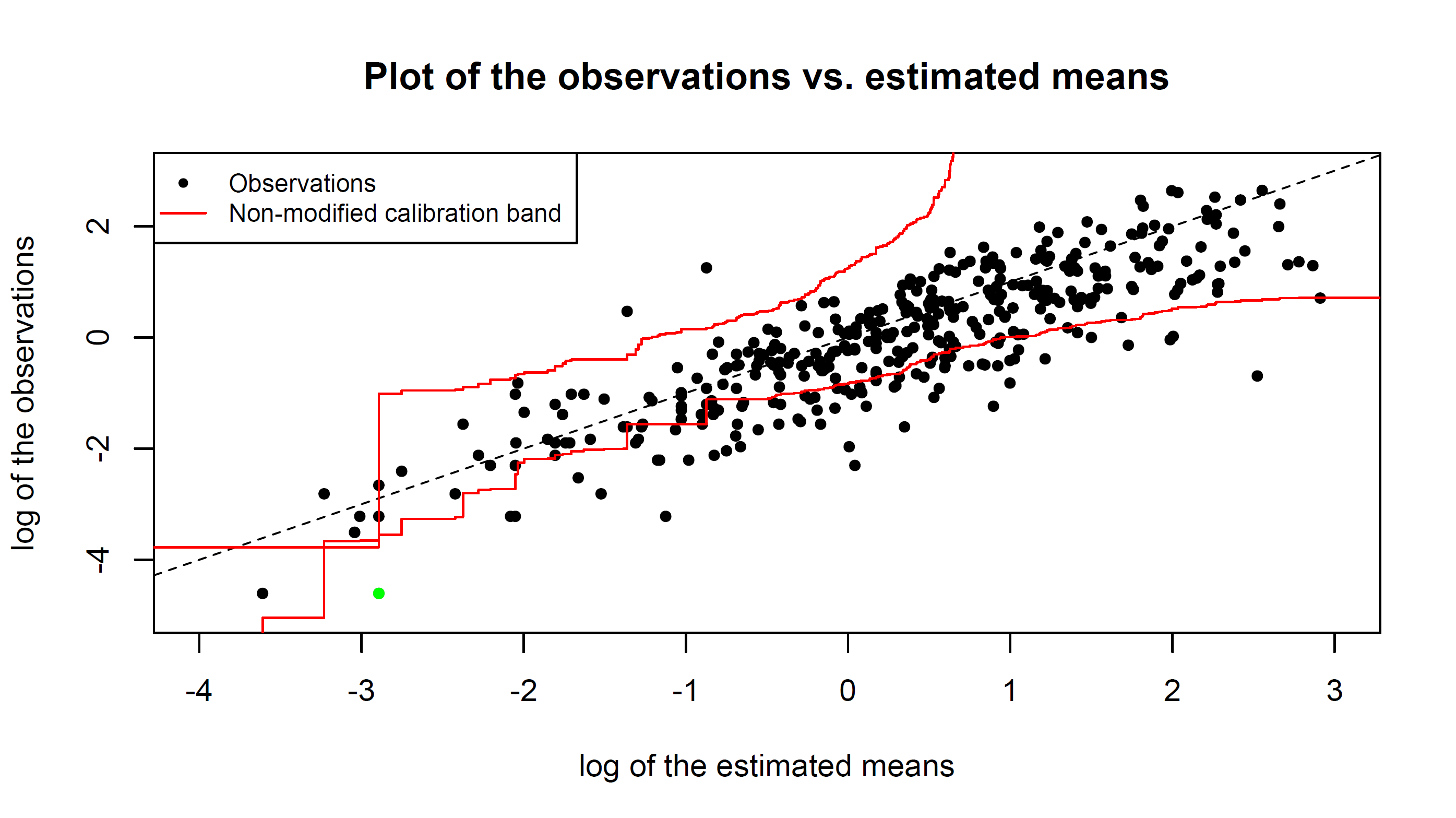}
    \caption{Plot of the observations versus the inverse Gaussian GLM mean estimates on the log scale. The non-modified calibration band is plotted in red, whereas the observations are drawn in black. The seventh observation is highlighted in green.}
    \label{Fig:observations}
\end{figure}

\begin{figure}[htb!]
    \centering
    \includegraphics[width = 0.7\linewidth]{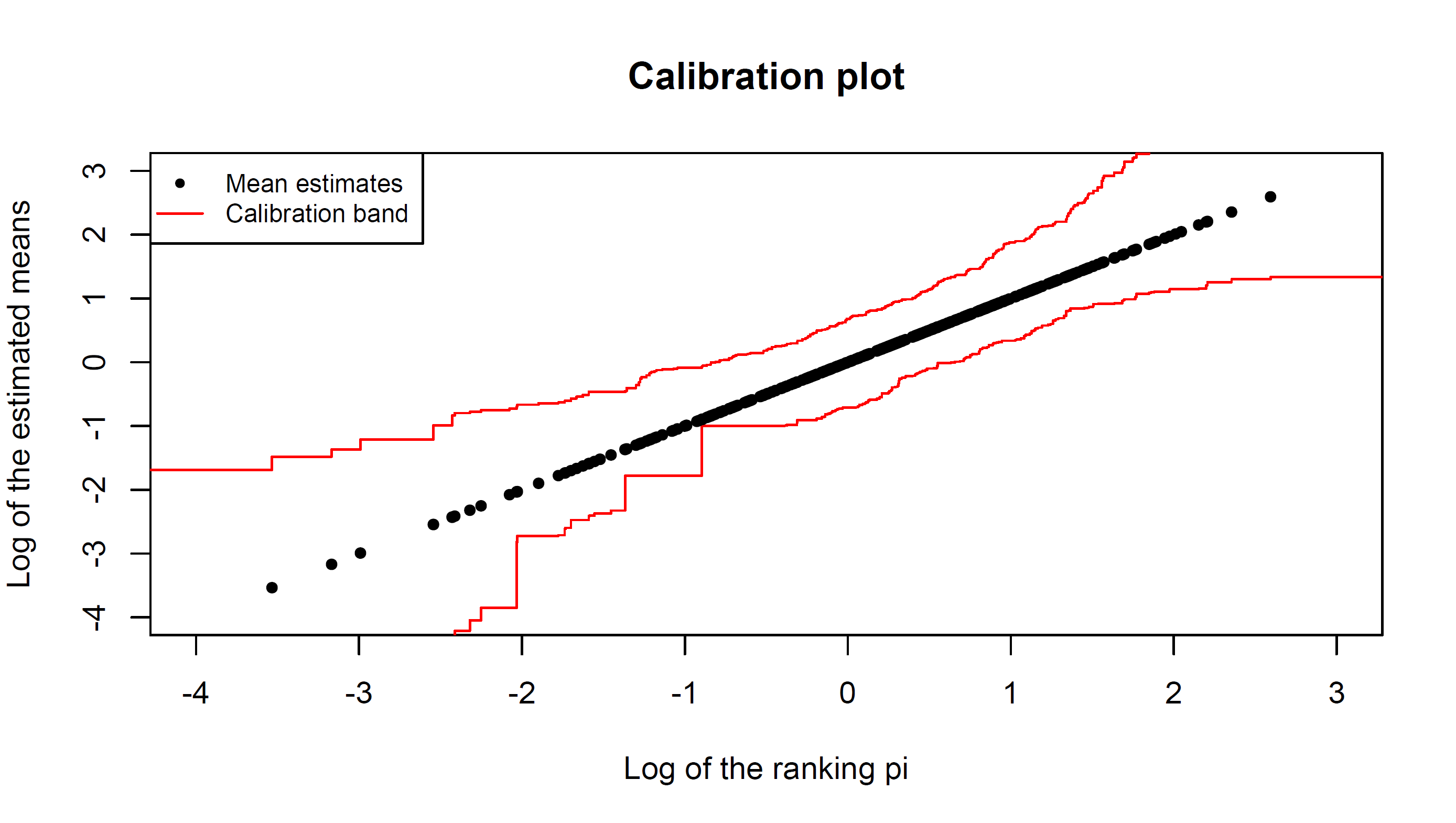}
    \caption{Calibration plot of the gamma GLM on the log scale. The calibration band is plotted in red, whereas the mean estimates are drawn in black.}
    \label{Fig:Gamma GLM}
\end{figure}

\begin{figure}[htb!]
    \centering
    \includegraphics[width = 0.7\linewidth]{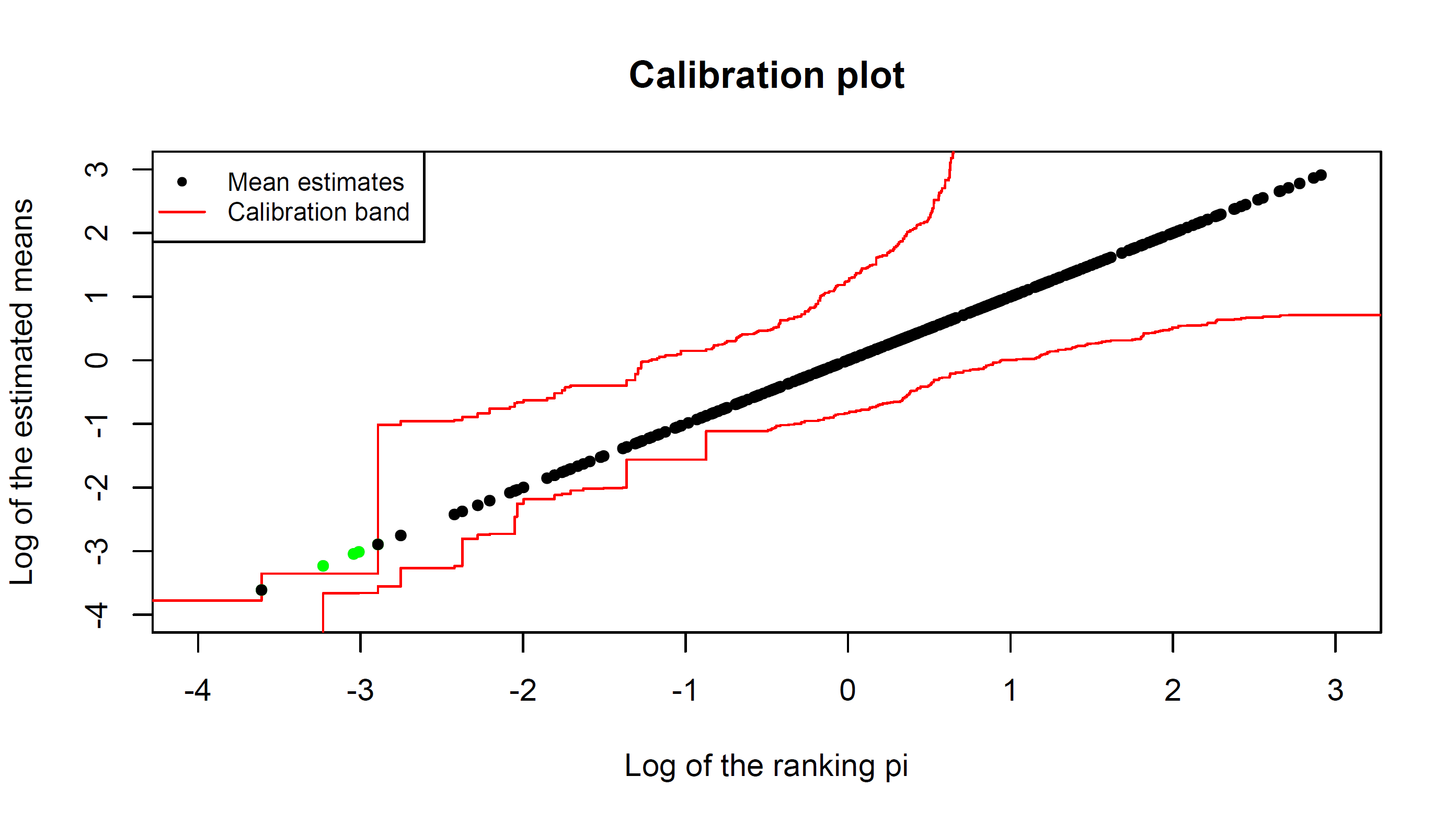}
    \caption{Calibration plot of the inverse Gaussian GLM on the log scale. The calibration band is plotted in red. The mean estimates $(\widehat{\mu}^{IG}(\bx_i))_{i=1}^n$ falling within the band are drawn in black, whereas those falling outside the bands are drawn in green.}
    \label{Fig:IG GLM}
\end{figure}

In order to assess the auto-calibration of the above models, we construct full calibration bands on the mean of the foliage biomass by using the GLMs themselves as ranking functions and choosing the confidence level $1-\alpha = 0.95$. The gamma and inverse Gaussian bands are shown in Figures \ref{Fig:Gamma GLM} and \ref{Fig:IG GLM}. For the latter band, note that the constructed upper and lower bounds for small mean estimates are crossing, see Figure \ref{Fig:observations}. 
We thus use the method in \eqref{cal not crossing} in order to obtain a non-crossing band in Figure \ref{Fig:IG GLM}. Those crossings may happen because the ranking of the estimated means may not match the ranking of the smallest observations. The assumed variance function when using the inverse Gaussian GLM fulfills
$$
V(\mu) \propto \mu^3,
$$
which implies that the lower and upper bounds in \eqref{L_i general} and \eqref{U_i general} are extremely close to the weighted partial sums $Z_{j:k}$ in  \eqref{sum j k} for small values of the sums. Outliers are thus rare and we see in Figure \ref{Fig:observations} that in this example, the seventh observation, highlighted in green, is responsible for a low upper band on the left end of the interval.

The conclusion of the calibration plot in Figure \ref{Fig:Gamma GLM} is not to reject the auto-calibration of the gamma GLM at a confidence level $1-\alpha = 0.95$ as all the mean estimates fall within the constructed calibration band. However, the plot in Figure \ref{Fig:IG GLM} leads to the rejection of the auto-calibration of the inverse Gaussian GLM because at the bottom left of the plot, some of the mean estimates fall outside the constructed band, those are plotted in green. In order to perform these tests, we assumed that $\widehat{\mu}^{gamma}(\bX)$ and $\widehat{\mu}^{IG}(\bX)$ are absolutely continuous random variables with strictly positive density over their supports.

Interestingly, the conclusion of the test derived by Denuit et al.~\cite{Denuit2} is in this case different as the auto-calibration of both models is not rejected at the level $1-\alpha = 0.95$. To perform the latter test, non-parametric Monte Carlo methods have to be used to compute the p-value $\hat{p}_{auto}$, and the null-hypothesis of auto-calibration is rejected at confidence level $1-\alpha$ whenever $\hat{p}_{auto} < \alpha$. The p-values obtained for $\hat{p}_{auto}$ for both GLMs are summarized in Table \ref{Tab:p_auto} by performing $B = 500$ Monte-Carlo simulations, we refer the reader to Section 4 of Denuit et al.~\cite{Denuit2} for more technical details. While these values are both larger than $\alpha = 0.05$, note that the value for the gamma GLM is close to $1$, whereas the value for the inverse Gaussian GLM is close to the chosen significance level. This indicates that the gamma GLM is indeed closer to auto-calibration than the inverse Gaussian GLM.
\begin{table}[htb!]
    \centering
    \begin{tabular}{l c c }
      \toprule
      & gamma GLM & inverse Gaussian GLM \\ 
      \midrule
      $\hat{p}_{auto}$ & $0.89$ & $0.08$ \\ 
      \bottomrule
    \end{tabular}
    \caption{The p-values $\hat{p}_{auto}$ for the gamma and inverse Gaussian GLMs.}
    \label{Tab:p_auto}
\end{table}
For this small-sample sized real dataset, our test manages to detect violations of auto-calibration of the inverse Gaussian GLM in contrast to the test of Denuit et al.~\cite{Denuit2}. The reason might be that the latter test relies on asymptotic results and should thus only be used on large datasets, whereas our test adapts to the distribution of the responses and can be used for any sample size as the construction of the calibration band does not rely on any asymptotic result. This constitutes an advantage over the other methods when it comes to assess the auto-calibration of a small dataset.

\subsection{Example 3 : calibration bands for simulated inverse Gaussian responses}

In Section \ref{sec dispersion}, we emphasized the importance of using suitable dispersion estimates in order to construct calibration bands as all of our results are based on the assumption of a known and fixed dispersion estimate. Throughout this example, we construct calibration bands on the mean of simulated inverse Gaussian responses that rely on different dispersion estimates. Our goal is to show that although the dispersion parameter $\varphi$ is unknown, and has to be estimated, its influence on the resulting bands is comparably small, allowing us to use the statistical tests of Section \ref{sec tests} in practice.

\begin{figure}[htb!]
\begin{center}
\begin{minipage}[t]{0.35\textwidth}
\begin{center}
\includegraphics[width=\textwidth]{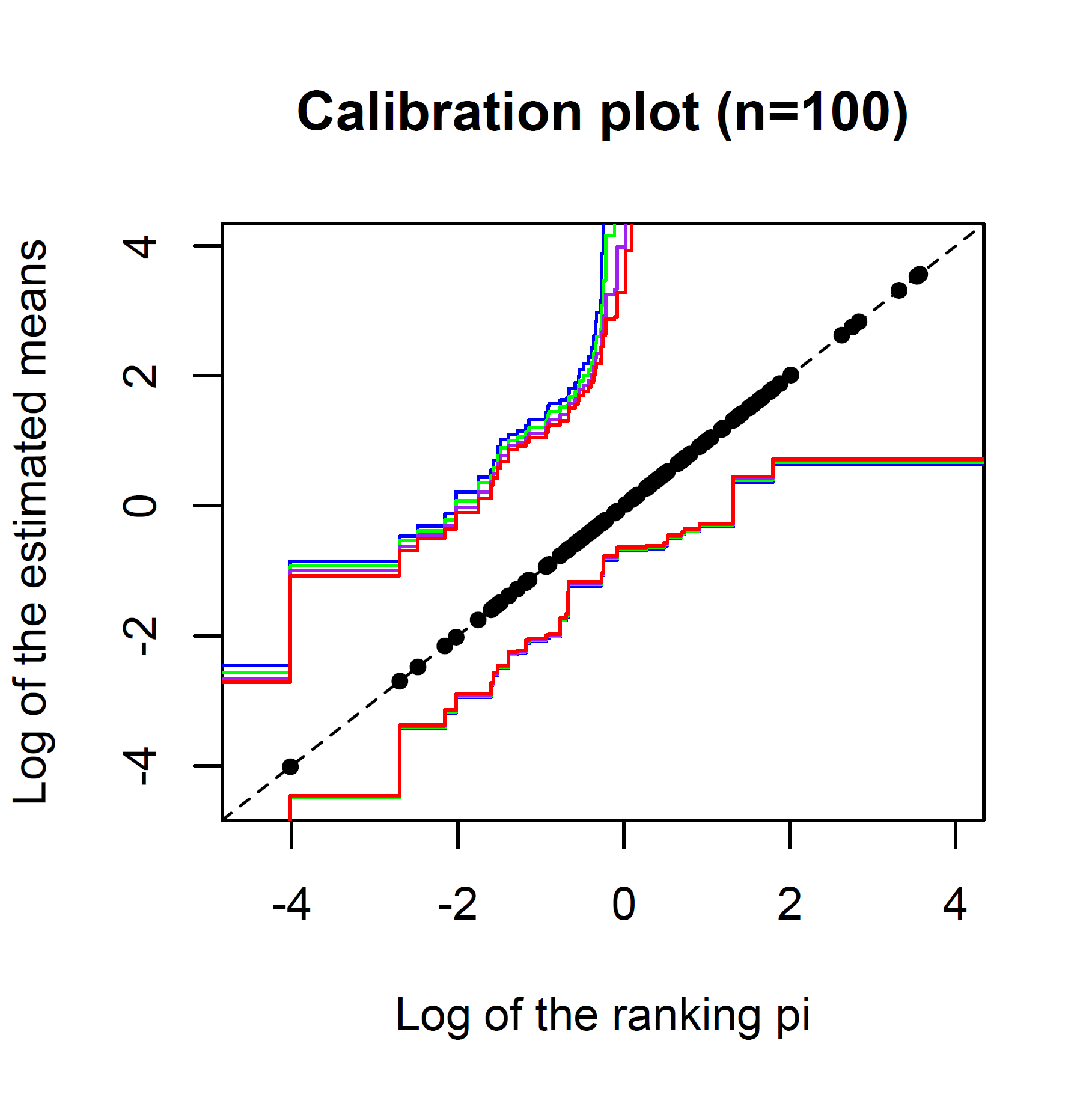}
\end{center}
\end{minipage}
\begin{minipage}[t]{0.35\textwidth}
\begin{center}
\includegraphics[width=\textwidth]{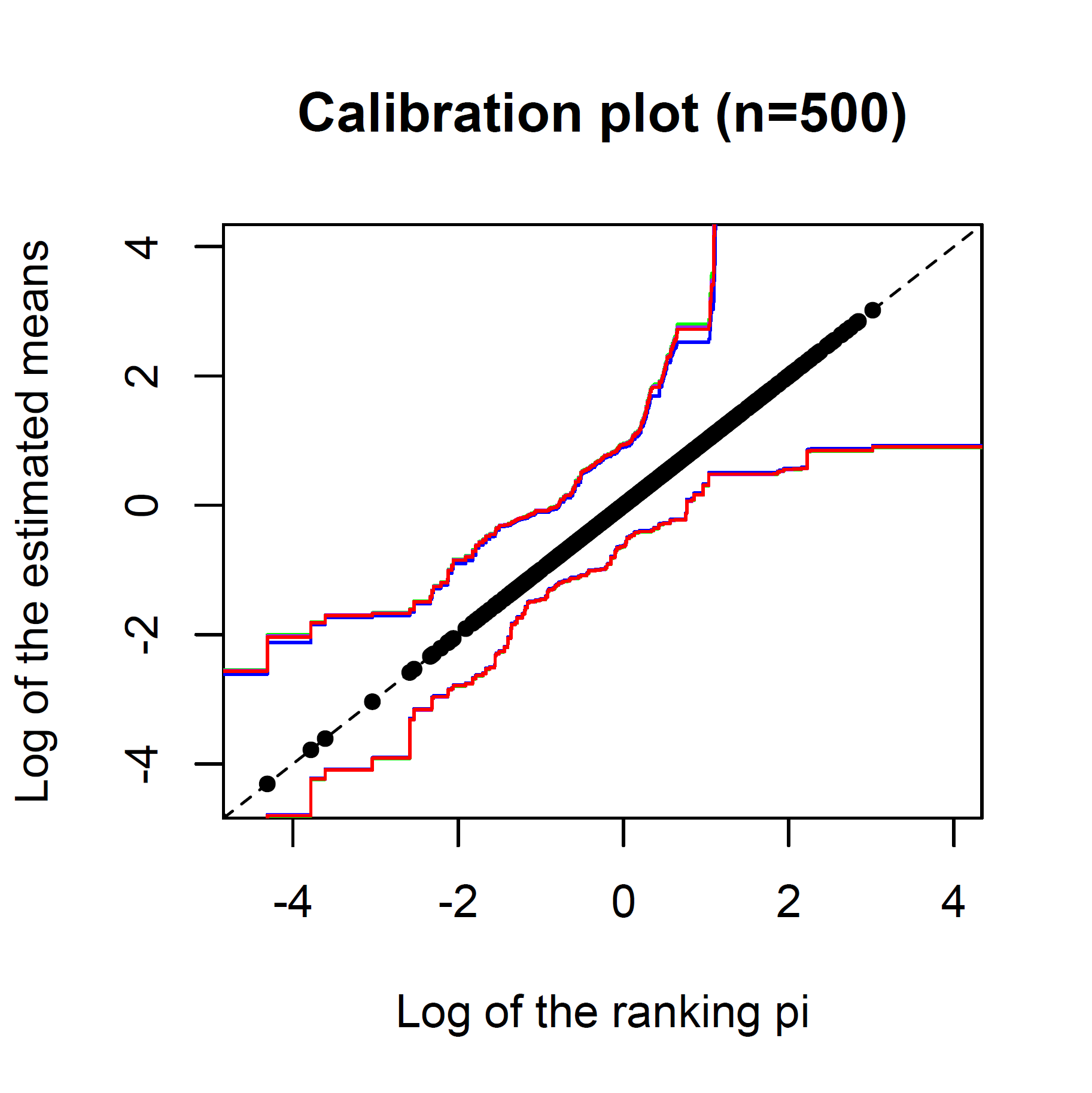}
\end{center}
\end{minipage}
\begin{minipage}[t]{0.35\textwidth}
\begin{center}
\includegraphics[width=\textwidth]{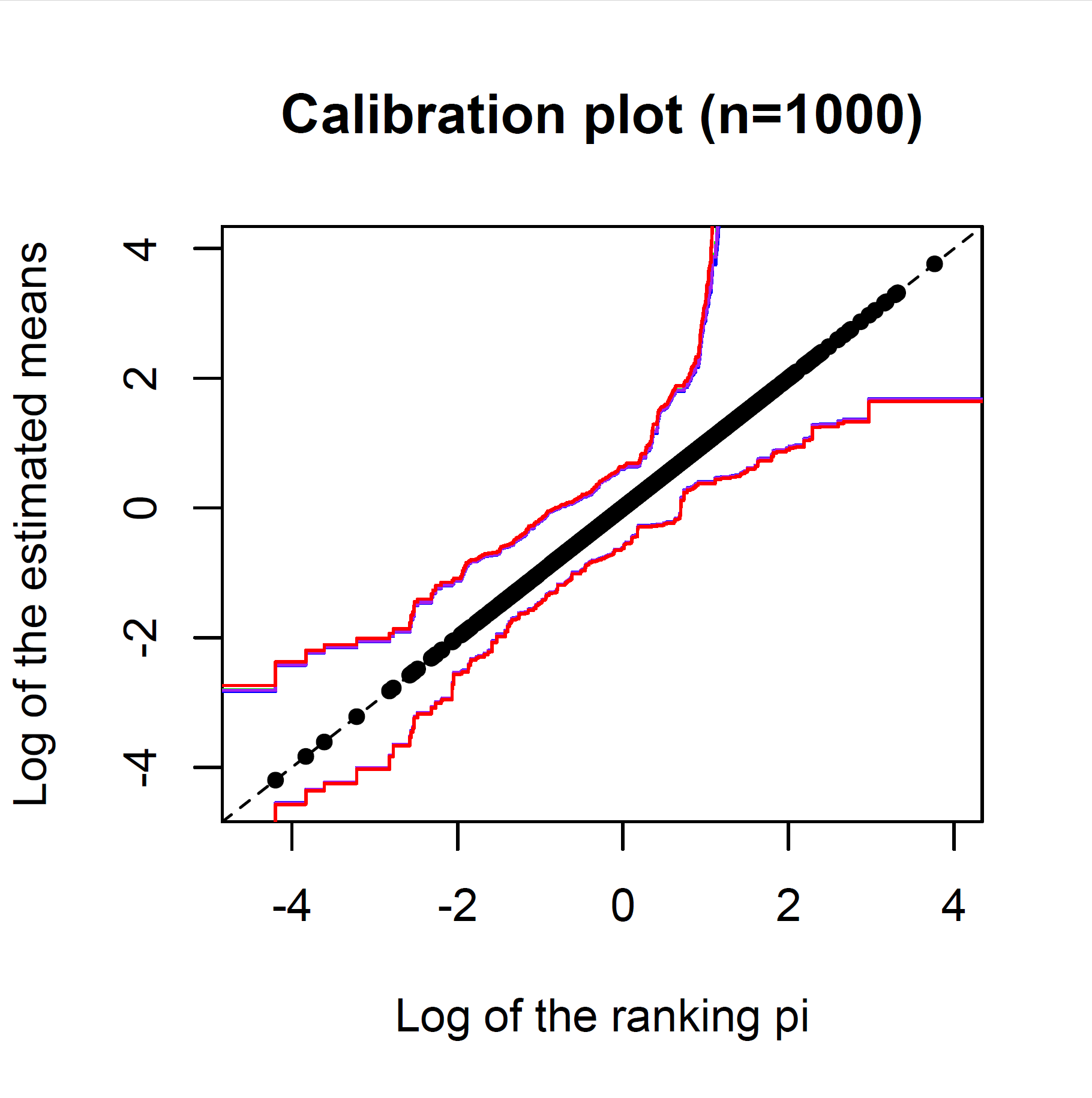}
\end{center}
\end{minipage}
\begin{minipage}[t]{0.35\textwidth}
\begin{center}
\includegraphics[width=\textwidth]{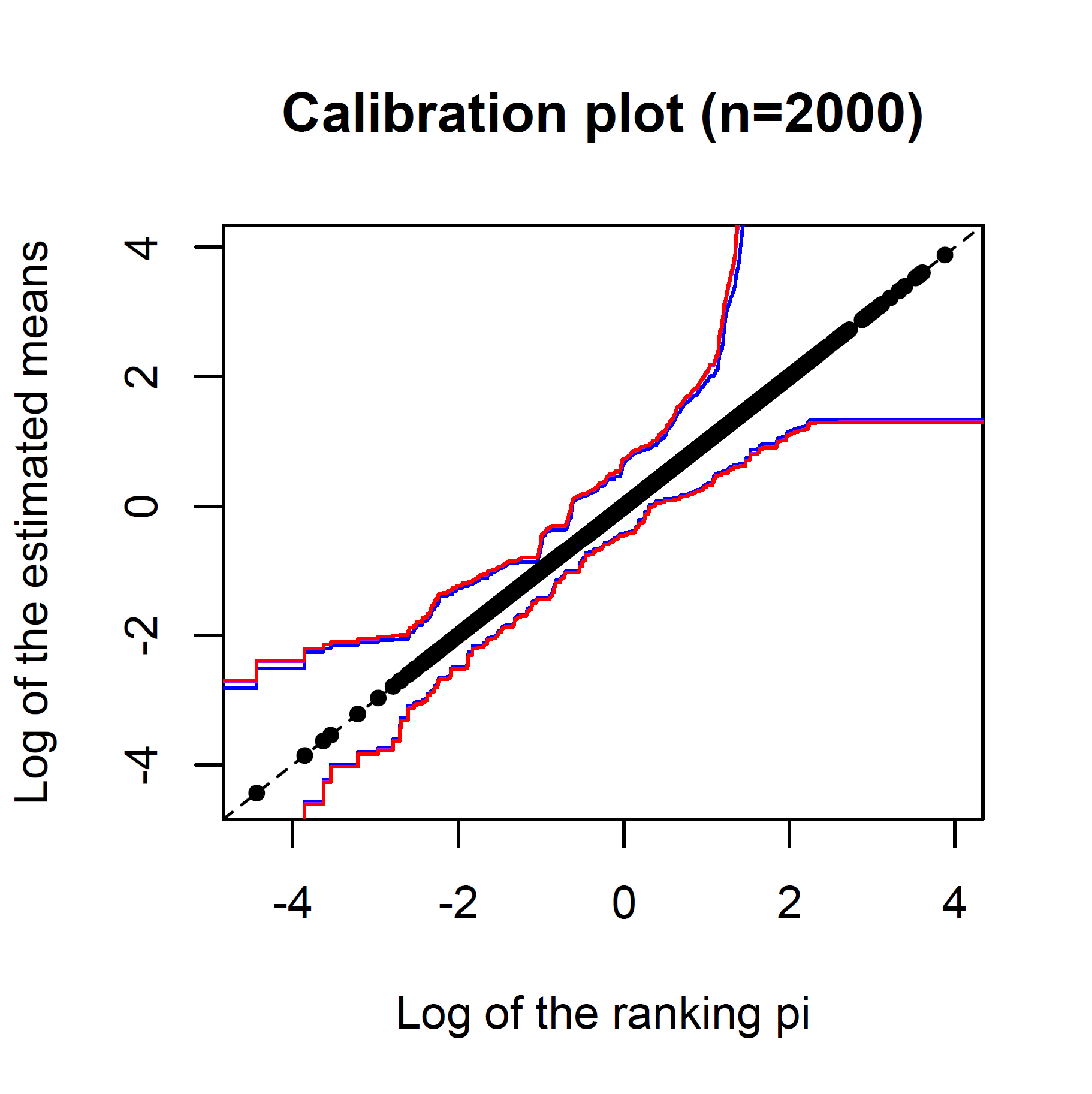}
\end{center}
\end{minipage}
\end{center}
\vspace{-.7cm}
\caption{Calibration plots of the regression function $\tilde{\mu}^{gamma} : \X \to (0, \infty)$ on the log scale for four different sample sizes. For each plot, four different calibration bands are constructed using $\pi(\cdot) = \tilde{\mu}^{gamma}(\cdot)$ as a ranking function and the true dispersion parameter (red) as well as the Pearson estimate (blue), the deviance estimate (green) and the MLE estimate (purple) for the dispersion. The mean estimates are drawn in black.}
\label{Fig:4 plots}
\end{figure}

In order to obtain realizations of $n$ independent inverse Gaussian responses, we first replicate the feature space of the above real dataset by simulating $n$ feature components $\tilde{\bx}_i$ lying in the set $\X$ given in \eqref{feature space ex 2}. To this end, we assume that the tree diameters follow a gamma distribution, with scale and shape parameters chosen based on the empirical mean and standard deviation of the observed diameters in the above real dataset. In addition, the origin component of the simulated features is generated using the empirical distribution of origins from the same dataset. This leads to $n$ new features $(\tilde{\bx}_i)_{i=1}^n$. We then choose a dispersion parameter of $\varphi = 1.26$, corresponding to the Pearson estimate of the above dataset. Finally, we assume that the responses satisfy
$$
\tilde{Y}_i \sim \textnormal{IG}(\hat{\mu}^{IG}(\tilde{\bx}_i), \varphi), \quad 1\leq i \leq n.
$$
This provides us with an example where the true model is given by an inverse Gaussian GLM with an interaction term between the two feature components, see Section \ref{sec ex 2}. As the inverse Gaussian distribution has the advantage that the MLE dispersion estimate can be expressed in closed form, we aim at constructing four different full calibration bands with confidence level $1-\alpha = 0.95$ by using successively the true dispersion parameter, the Pearson dispersion estimate, the deviance dispersion estimate and the MLE dispersion estimate, we refer to Chapter 11.6 in Dunn--Smyth \cite{Dunn--Smyth} for more details about these estimates. Using the constructed bands, we then assess the calibration of a newly fitted gamma GLM $\tilde{\mu}^{gamma} : \X \to (0, \infty)$ on the new dataset $(\tilde{y}_i, \tilde{\bx})_{i=1}^n$ for $n \in \{100, 500, 1000, 2000\}$, and the results are provided in Figure \ref{Fig:4 plots}. Although all the treated datasets are small, we notice that the bands nearly coincide for all sample sizes. This example shows that when the model is well-specified, i.e.~when the data generating process matches the model's assumptions, using an estimate for the dispersion instead of the true dispersion parameter to construct calibration bands does not have a major impact on the resulting bands, even for small datasets. The values of the dispersion estimates are shown in Table \ref{Tab:dispersion estimates}, where we see that even a 20\% difference with respect to the true parameter $\varphi = 1.26$ does not lead to any significant impact on the bands. Finally, note that while the statistical tests we derive can be applied for all sample sizes, we see in Figure \ref{Fig:4 plots} that the calibration band gets narrower the larger the sample size is, leading to more powerful statistical tests.

\begin{table}[htb!]
    \centering
    \begin{tabular}{l c c c c}
      \toprule
      Sample sizes & $n=100$ & $n=500$ & $n=1000$ & $n=2000$\\
      \midrule
      Pearson estimate & $1.50$ & $1.18$ & $1.13$ & $1.04$ \\
      Deviance estimate & $1.40$ & $1.28$ & $1.16$ & $1.24$ \\
      MLE estimate & $1.32$ & $1.27$ & $1.15$ & $1.24$\\
      \bottomrule
    \end{tabular}
    \caption{Dispersion estimates for various sample sizes; the true value is $\varphi = 1.26$.}
    \label{Tab:dispersion estimates}
\end{table}
\subsection{Example 4 : calibration bands for a large sample-sized real dataset}
\label{MTPL sec}

After considering a small dataset in Section \ref{sec ex 2}, we study in this example a French motor third party liability (MTPL) real dataset available from the
{\sf R} \cite{R Core Team} package {\tt CASdatasets} hosted by Dutang--Charpentier \cite{Dutang--Charpentier}. This dataset contains information on insurance policies and claim frequency of more than half a million French car drivers. 
We follow Listing 13.1 in Wüthrich--Merz \cite{WM2023} in order to clean the data, leading to a portfolio of $n = 678,007$ insurance policies and $26,383$ claims\footnote{The cleaned dataset can be downloaded under \url{https://people.math.ethz.ch/~wueth/Lecture/freMTPL2freq.rda}}. For each policy $1 \leq i \leq n$, the resulting dataset provides the number of claims $N_i \in \N$ that occurred during an exposure period $v_i \in (0,1]$ (years-at-risk) and features containing information of the policyholder as, for example, the age of the driver, the brand and power of their cars, or the region of residence. The total exposure at risk is equal to $358,359$ years, indicating that some policyholders were covered for a period of less than one year, and as one might expect in motor liability insurance, most policies do not lead to any claim, see Table \ref{Tab:description data}. We refer to Section 13.1 in Wüthrich--Merz \cite{WM2023} for an extended description of the dataset. 

\begin{table}[htb!]
    \centering
    \begin{tabular}{l c c c c c c}
      \toprule
      Number of claims occurred for each policy& 0 & 1 & 2 & 3 & 4 & 5\\
      \midrule
      Number of policies & $653,069$ & $23,571$ & $1298$ & $62$ & $5$ & $2$ \\
      Total exposure & $341,090$ & $16,315$ & $909$ & $42$ & $2$ & $1$ \\
      \bottomrule
    \end{tabular}
    \caption{Number of policies and total exposure within the portfolio that is split with respect to the number of claims occurred for each policy.}
    \label{Tab:description data}
\end{table}
In this section, we aim at modeling the claim frequency of each policyholder. For this, we follow Listings 5.1-5.2 in Wüthrich--Merz \cite{WM2023} in order to pre-process the available features. Moreover, we consider a subset of the features by only keeping the information about the policyholder and not their cars. That is, we use $3$ continuous feature components and $2$ categorical ones\footnote{The used features are : {\sf BonusMalusGLM, DensityGLM, AreaGLM, DrivAgeGLM, Region}, see Section 5.2.4 in Wüthrich--Merz \cite{WM2023}.}, leading to a feature space 
\begin{equation*}
    \X \subset \R^3 \times \{0,1\}^6 \times \{0,1\}^{21}.
\end{equation*}
After pre-processing the categorical variables using dummy coding, we fit a Poisson GLM with the canonical link on the whole dataset in order to estimate the claim frequency of each policyholder with feature $\bx \in \X$, given by the true mean function $\mu^{*}: \X \to (0, \infty)$. That is, we assume that
\begin{equation*}
   N_i \sim \textrm{Poi} (\mu^{*}(\bx_i)v_i),\quad \textrm{for }1 \leq i \leq n,
\end{equation*}
where $\bx_i \in \X$ is the considered feature of the policy $i$. 
We call the resulting estimated regression function  $\widehat{\mu}^{Poi} : \X \to (0, \infty)$, and this function satisfies
\begin{equation*}
    \min_{1 \leq i \leq n} \widehat{\mu}^{Poi}(\bx_i) = 0.024 \quad \textrm{and} \quad  \max_{1 \leq i \leq n} \widehat{\mu}^{Poi}(\bx_i) = 1.292.
\end{equation*}
This means that the model predicts that an accident occurs on average once every 40 years for some drivers, while for others, it predicts that more than one accident occurs each year on average. In order to assess whether the obtained regression function is calibrated, we construct a calibration band on the claim frequency. As the dataset is large, we bin the responses according to their estimated means for computational reasons. To do so, we first define $Y_i = N_i/v_i$ and use the convolution formula for the reproductive form of the EDF in order to derive $L = 5000$ new responses
\begin{equation*}
    \tilde{Y}_l = \frac{\sum_{i=1}^n v_iY_i \, \mathds{1}_{\left\{\widehat{\mu}^{Poi}(\bx_i) \in \mathcal{I}_l\right\}}}{\sum_{i=1}^n v_i \, \mathds{1}_{\left\{\widehat{\mu}^{Poi}(\bx_i) \in \mathcal{I}_l \right\}}}, \quad \textrm{for }1 \leq l \leq L, 
\end{equation*}
where the intervals $\mathcal{I}_l =[a_{l-1}, a_{l})$ are delimited by some partition $(a_l)_{l=0}^L$ of $[0.024,1.292]$ such that the volumes of all the binned responses are approximately equal. This can be achieved by ranking the responses $(Y_i)_{i=1}^n$ according to their mean estimates and using a weighted quantile binning. This procedure leads to the new volumes
\begin{equation*}
    \tilde{v}_l = \sum_{i=1}^n v_i \, \mathds{1}_{\left\{\widehat{\mu}^{Poi}(\bx_i) \in \mathcal{I}_l \right\}} \in [70.7, 72.64] , \quad \textrm{for }1 \leq l \leq L.
\end{equation*}

\begin{figure}[htb!]
\begin{center}
\includegraphics[width=0.7\textwidth]{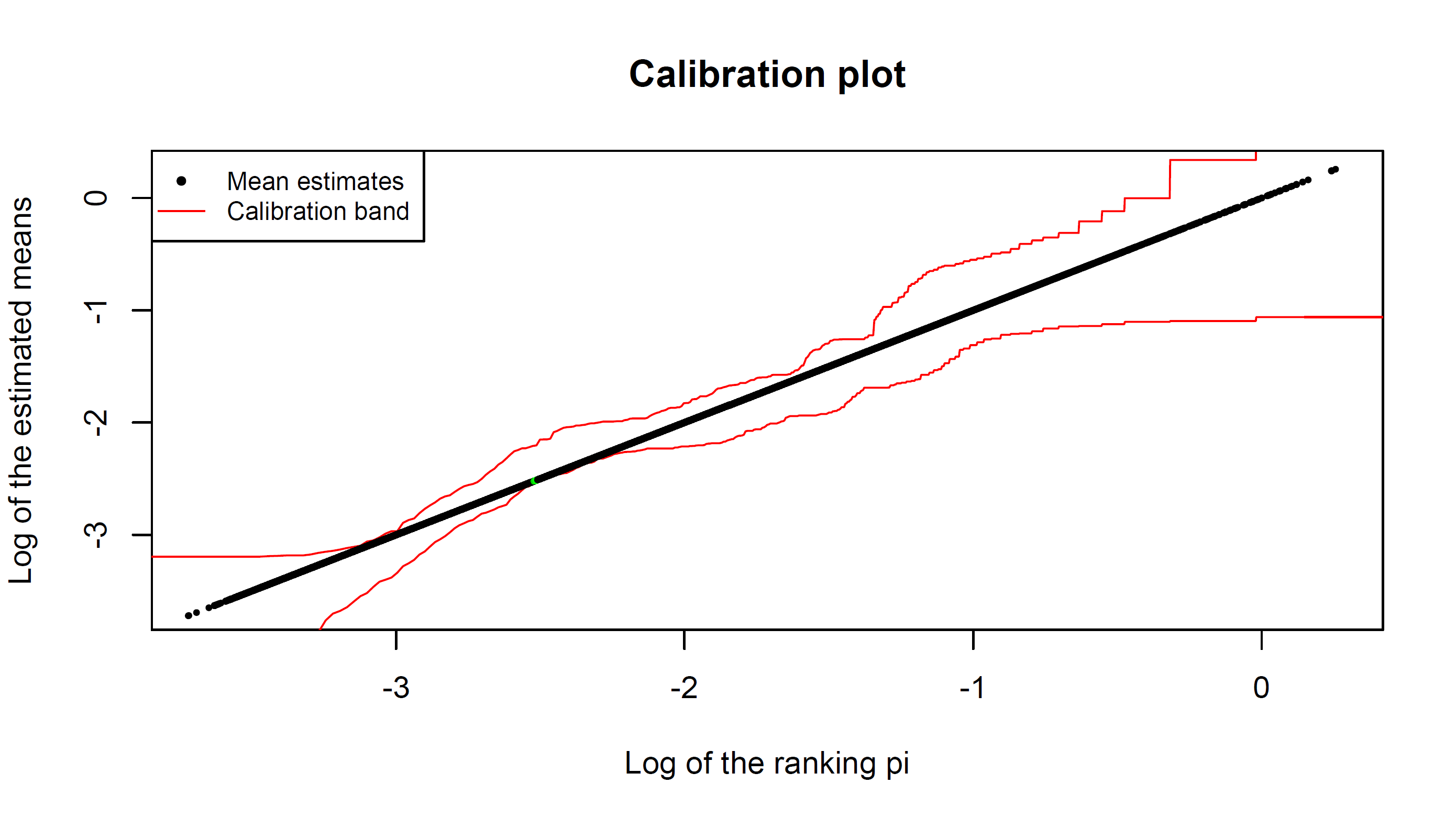}
\end{center}
\begin{center}
\begin{minipage}[t]{0.35\textwidth}
\begin{center}
\includegraphics[width=\textwidth]{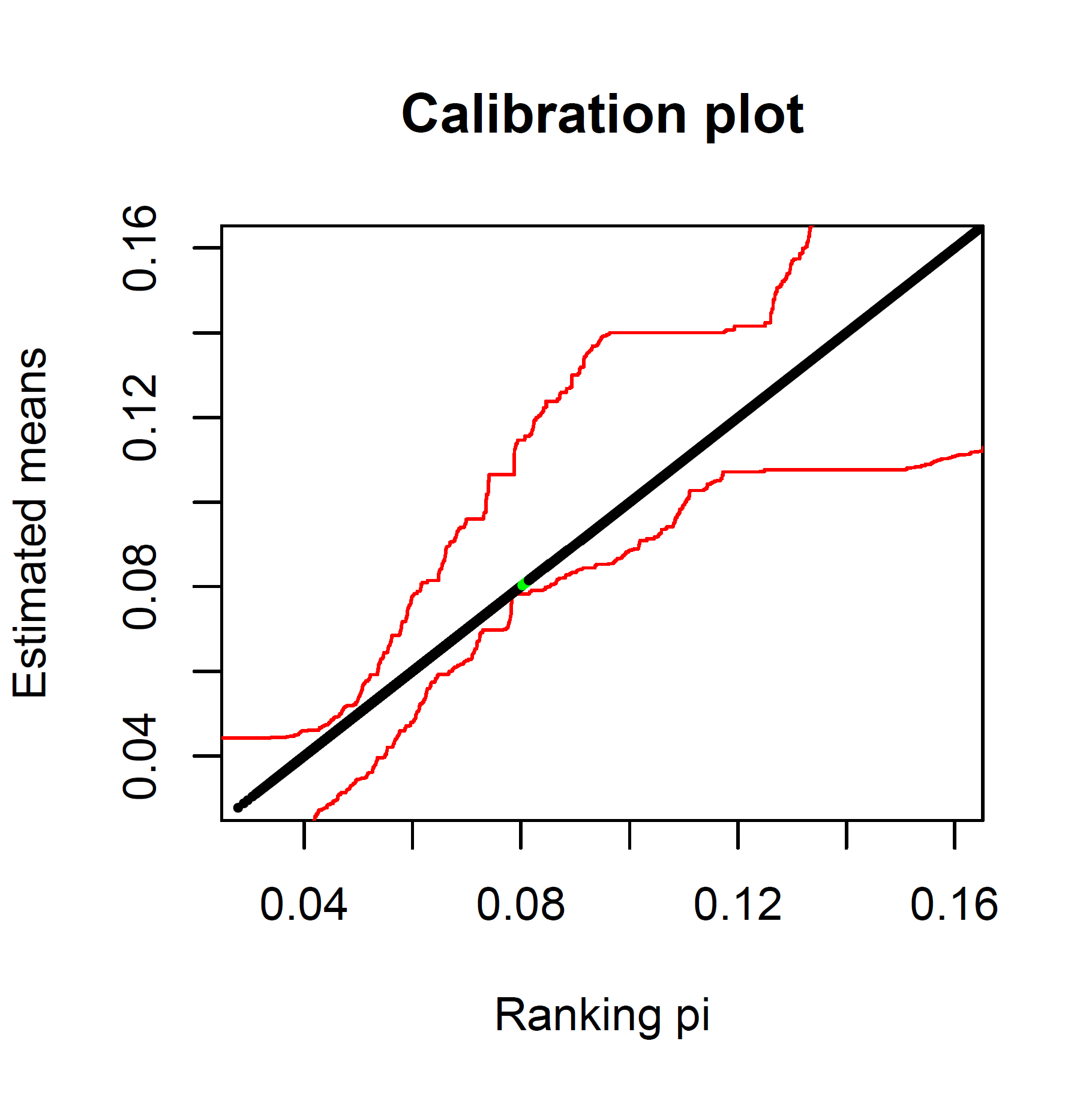}
\end{center}
\end{minipage}
\begin{minipage}[t]{0.35\textwidth}
\begin{center}
\includegraphics[width=\textwidth]{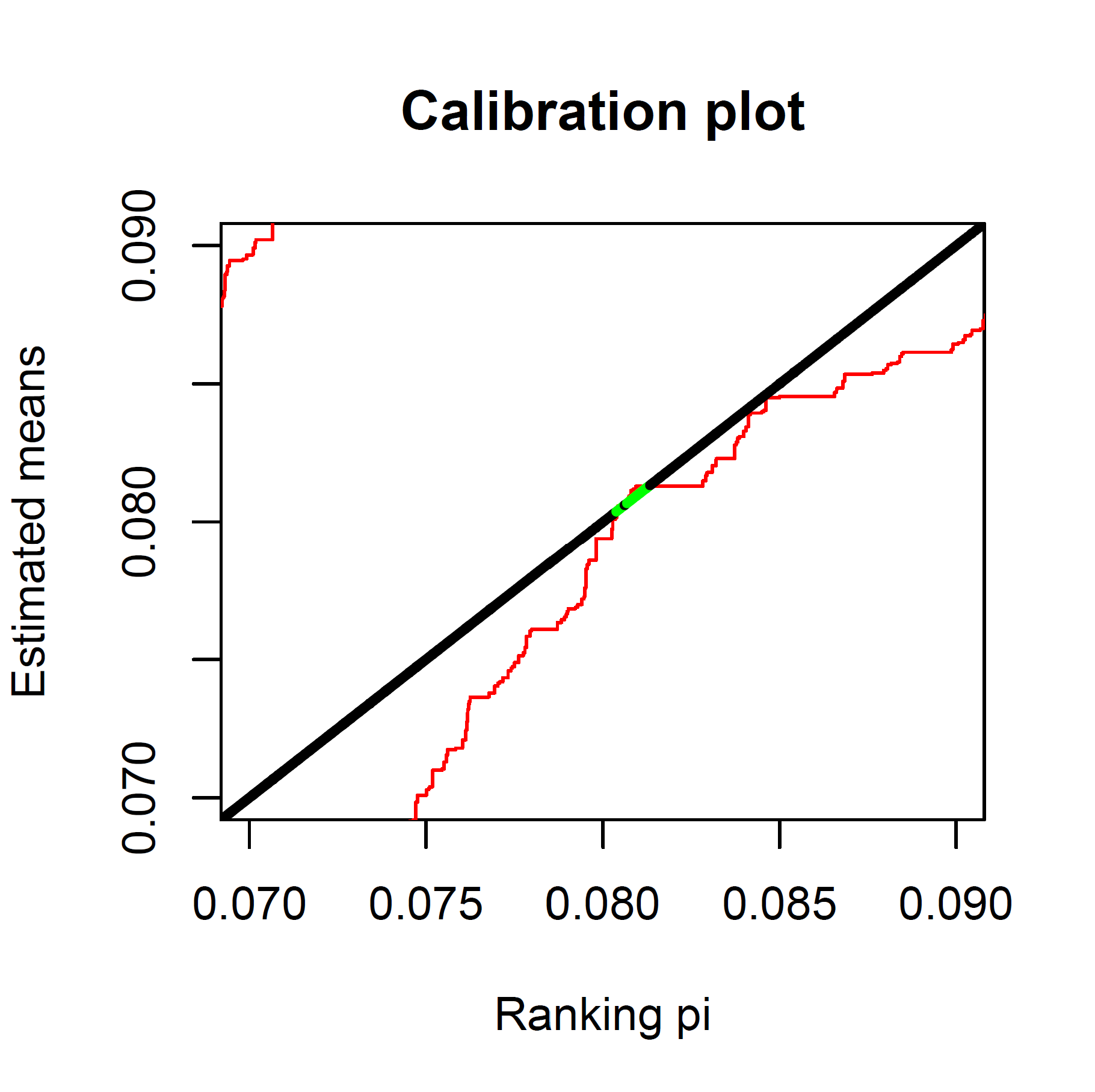}
\end{center}
\end{minipage}
\end{center}
\vspace{-.7cm}
\caption{Calibration plot of the regression function $\widehat{\mu}^{Poi} : \X \to (0, \infty)$ on the log scale (above). Zoomed versions of this plot are provided on the linear scale (below). The calibration band is constructed using $\pi(\cdot) = \widehat{\mu}^{Poi}(\cdot)$ as a ranking function and is plotted in red. The mean estimates $(\widehat{\mu}^{Poi}(\bx_i))_{i=1}^n$ falling within the band are drawn in black, whereas those falling outside the bands are drawn in green.}
\label{Fig:Poisson_GLM}
\end{figure}

The same binning applied to the mean estimates $(\widehat{\mu}^{Poi}(\bx_i))_{i=1}^n$ allows us to derive the mean estimates of the new responses $\left(\widehat{\mu}^{Poi}_l\right)_{l=1}^L$.
Using the realizations of the binned responses as observations and $\pi(\cdot) = \widehat{\mu}^{Poi}(\cdot)$ as a ranking function, we construct a full calibration band for the confidence level $1-\alpha = 0.95$. This band is drawn in Figure \ref{Fig:Poisson_GLM}, where we additionally plot the mean estimates $(\widehat{\mu}^{Poi}(\bx_i))_{i=1}^n$ against the corresponding rankings. Note that the resulting points all lie on the diagonal due to our choice of the ranking function. 

As we use the log scale for the upper plot in Figure \ref{Fig:Poisson_GLM}, we notice that the band is narrow for small means and wide for large means. The reason for this is that the aggregated exposure of the policies for which the estimated mean is below 0.2 corresponds to $97.5\%$ of the total exposure of the portfolio. In other words, the aggregated volumes used to compute the bounds fail to be large enough for mean estimates exceeding $0.2$. The lower plots in Figure \ref{Fig:Poisson_GLM} show that some mean estimates fall outside the calibration band, those are plotted in green. By assuming that $\widehat{\mu}^{Poi}(\bX)$ is an absolutely continuous random variable with strictly positive density over its support, the conclusion of the statistical test derived in Section \ref{Statistical tests for auto-calibration} is thus to reject calibration at a confidence level of $1-\alpha = 0.95$. 
This decision indicates that the regression function $\widehat{\mu}^{Poi}: \X \to (0, \infty)$ is too far from the true mean function although the violation only happens for a small part of the support of $\widehat{\mu}^{Poi}(\bX)$ in Figure \ref{Fig:Poisson_GLM}. We emphasize, however, that as the inequality in \eqref{prob bound J} might not be very sharp in the construction of the band, the lower left plot in Figure \ref{Fig:Poisson_GLM} hints that the regression function might not be sufficiently calibrated for other mean estimates too, which are close to the boundary of the band. 

\begin{figure}[htb!]
\begin{center}
\includegraphics[width=0.7\textwidth]{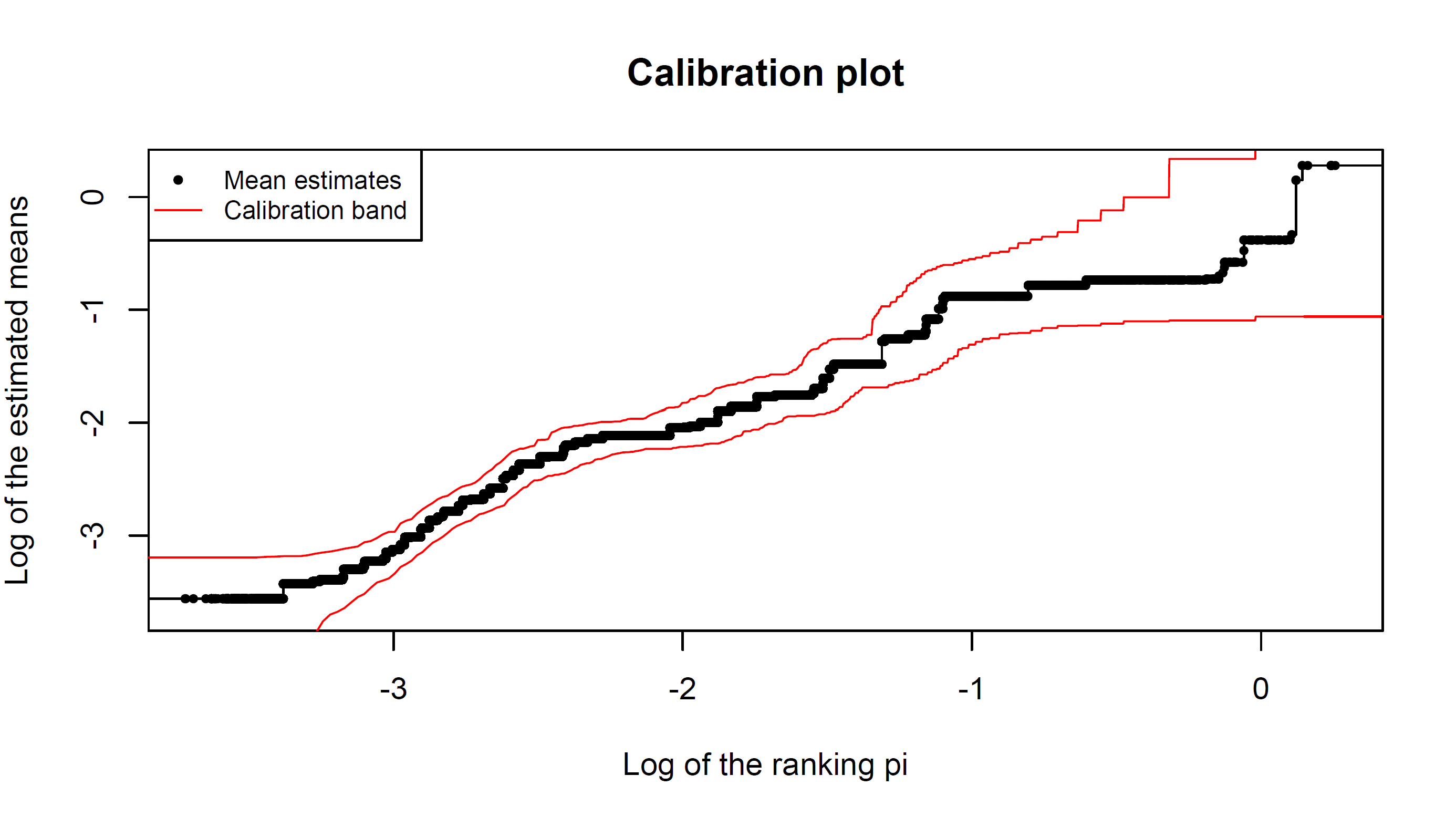}
\end{center}
\vspace{-.7cm}
\caption{Calibration plot of the regression function $\widehat{\mu}^{Poi}_{rec} : \X \to (0, \infty)$ on the log scale. The calibration band is constructed using $\pi(\cdot) = \widehat{\mu}^{Poi}(\cdot)$ as a ranking function and is plotted in red, whereas the mean estimates $(\widehat{\mu}^{Poi}_{rec}(\bx_i))_{i=1}^n$ are drawn in black.}
\label{Fig:Isotonic_GLM}
\end{figure}

Our next goal is to improve the obtained regression function. For this, we assume that it provides the correct ordering of the true mean function, but the decision of the above statistical test indicates a violation of the auto-calibration of  $\widehat{\mu}^{Poi}: \X \to (0, \infty)$, see Section \ref{Statistical tests for auto-calibration}. Therefore, we construct another regression function by applying the isotonic recalibration step proposed by Wüthrich--Ziegel \cite{Wuethrich_Ziegel} to the Poisson GLM. We call the new resulting regression function $\widehat{\mu}^{Poi}_{rec} : \X \to (0, \infty)$ and point out that the isotonic recalibration step is performed using the ranking provided by $\left(\widehat{\mu}^{Poi}(\bx_i)\right)_{i=1}^n$ and the exposures $(v_i)_{i=1}^n$ as weights. The obtained regression function $\widehat{\mu}^{Poi}_{rec} : \X \to \R$ is provided in Figure \ref{Fig:Isotonic_GLM}, where we draw the same calibration band as above, i.e., we assume again the ranking function to be $\pi(\cdot) = \widehat{\mu}^{Poi}(\cdot)$ in order to construct the band. 

\begin{figure}[htb!]
\begin{center}
\includegraphics[width=0.7\textwidth]{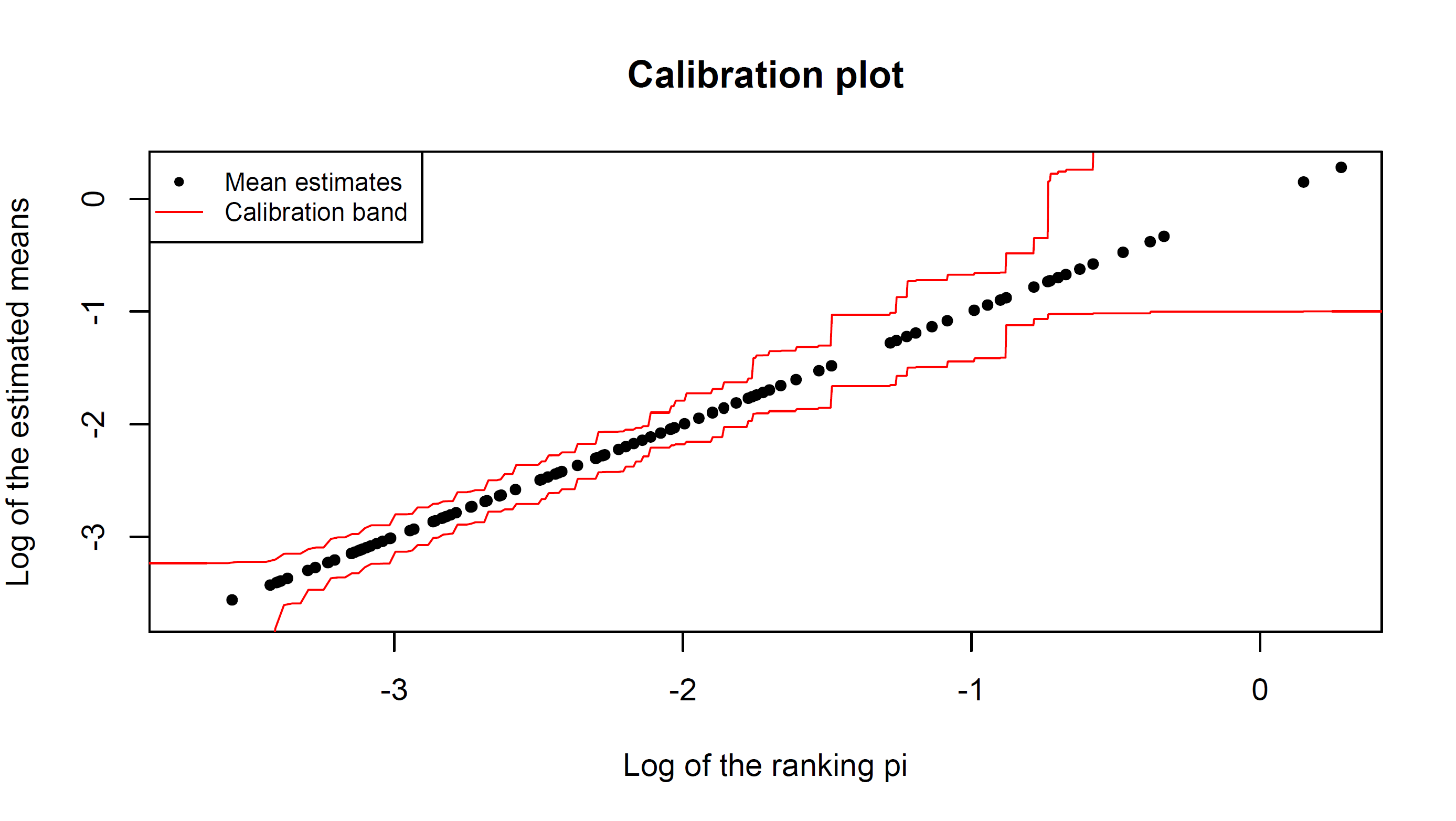}
\end{center}
\caption{Calibration plot of the regression function $\widehat{\mu}^{Poi}_{rec} : \X \to (0, \infty)$ on the log scale. The calibration band is constructed using $\pi(\cdot) = \widehat{\mu}^{Poi}_{rec}(\cdot)$ as a ranking function and is plotted in red, whereas the mean estimates $(\widehat{\mu}^{Poi}_{rec}(\bx_i))_{i=1}^n$ are drawn in black.}
\label{Fig:last_plot}
\end{figure}

This time, all the mean estimates lie at the middle of the constructed band, leading us not to reject the calibration of this model. Moreover, we point out that the regression function $\widehat{\mu}^{Poi}_{rec} : \X \to \R$ is empirically auto-calibrated, we refer to Wüthrich--Ziegel \cite{Wuethrich_Ziegel}. That is, the isotonic recalibration step provides an empirically auto-calibrated regression function, for which we do not reject the null-hypothesis of calibration. Note that one could alternatively construct a calibration band using using $\pi(\cdot) = \widehat{\mu}^{Poi}_{rec}(\cdot)$ as a ranking function in order to assess the calibration of $\widehat{\mu}^{Poi}_{rec} : \X \to \R$. The corresponding calibration plot is given in Figure \ref{Fig:last_plot} and, again, we do not reject the null-hypothesis of calibration as all the mean estimates lie within the band.

\subsection{Example 5 : power of the statistical test for calibration}

Finally, we study the power of the statistical test for calibration \eqref{classical test} derived in Section \ref{stat tests for calibration}. To do so, we consider the example of Wüthrich \cite{Wüthrich}. That is, we first simulate $n$ i.i.d.~mean values $\mu_i$ from the law described by
\begin{equation*}
    \mu_i = 
    \begin{cases}
    10, \quad \textrm{with }p = 0.1,   \\
    11, \quad \textrm{with }p = 0.15, \\
    12, \quad \textrm{with }p = 0.25, \\
    13, \quad \textrm{with }p = 0.25, \\
    14, \quad \textrm{with }p = 0.15, \\
    15, \quad \textrm{with }p = 0.1.
    \end{cases}
\end{equation*}

Then, we permute the indices of the sampled means so that \eqref{ranking no regression} holds and simulate $n$ independent responses $Y_i$ by assuming
\begin{equation*}
    Y_i \sim \Gamma(3 \mu_i, 3), \quad \textrm{for } 1 \leq i \leq n.
\end{equation*}
The above permutation ensures that the calibration bands are constructed using responses that are correctly ranked according to their means. Our aim is to evaluate the power of the statistical test in \eqref{classical test}. For this, we define the mean estimates
\begin{equation}
    \label{estimator gamma}
    \widehat{\mu}_i = \mu_i, \quad \textrm{for } 1 \leq i \leq n,
\end{equation}
and assess their calibration when the realizations $(y_i)_{i=1}^n$ of the above simulated responses are shifted. We consider here two different types of contamination. First, we contaminate the observations by a \textit{global shift} $\delta \in \{0, 0.5, 1\}$, i.e.,
\begin{equation*}
    y_i^{\delta} = y_i + \delta, \quad \textrm{for } 1 \leq i \leq n,
\end{equation*}
and we want to test for the calibration of $(\widehat{\mu}_i)_{i=1}^n$ for these new observations. Then, the same procedure is repeated but this time, only observations being associated to a given single mean are shifted, i.e., observations are transformed such that 
\begin{equation*}
    y_i^{l,\delta} = y_i + \delta \mathds{1}_{\{\mu_i = l\}},
\end{equation*}
for $l \in \{10, 13, 15 \}$. We refer to this transformation as a \textit{local shift of level $l$}.  

For both of these shifts, the mean estimates $(\widehat{\mu}_i)_{i=1}^n$ are calibrated if and only if $\delta = 0$ and we want to understand whether the statistical test in \eqref{classical test} is able to detect these violations of calibration. The results, showing the number of rejections of the calibration of $(\widehat{\mu}_i)_{i=1}^n$ at a confidence level $1-\alpha = 0.95$, are summarized in Table \ref{Tab1}. They should be compared to the plots on page 12 of Wüthrich \cite{Wüthrich} as the mean estimates $(\widehat{\mu}_i)_{i=1}^n$ are calibrated if and only if they are auto-calibrated for the shifts considered in this example. Note that we only simulate $1000$ different samples, each containing $n=1000$ responses, while Wüthrich \cite{Wüthrich} performs $10,000$ simulations of size $n=1000$. Additionally, we consider here only a limited set of contaminations $\delta \in \{0,0.5,1\}.$ 

When the shift factor $\delta$ is equal to $0$, i.e., when the mean estimates $(\widehat{\mu}_i)_{i=1}^n$ are calibrated, we see in Table \ref{Tab1} that the rejection rate is equal to $5/1000$, which is $10$ times smaller than the significance level $\alpha = 0.05$. This is not surprising as the constructed calibration bands rely on the union bound inequality in \eqref{prob bound J}, which implies that the power of the corresponding statistical test might be (much) lower than the significance level. Furthermore, we see in Table \ref{Tab1} that the constructed statistical test is not fully capable of detecting small deviations from calibration. However, the test seems to be more effective at identifying such deviations when they occur on a global scale or at the middle of the range of interest.

\begin{table}[htb!]
    \centering
    \begin{tabular}{l c c c}
      \toprule
    Contamination $\delta$ & $0$ & $0.5$ & $1$\\
      \midrule
      Global shift &$5/1000$& $291/1000$ & $1000/1000$ \\
      Local shift of level $l = 10$ & - & $13/1000$ & $612/1000$ \\
      Local shift of level $l = 13$ & - &$128/1000$& $994/1000$ \\
      Local shift of level $l = 15$ & - & $3/1000$ & $270/1000$ \\
      \bottomrule 
    \end{tabular}
    \caption{Power of the performed statistical tests with confidence level $1-\alpha = 0.95$.}
    \label{Tab1}
\end{table}

The same experiment is then repeated but this time, all the observations are binned according to their estimated means $(\widehat{\mu}_i)_{i=1}^n$ in \eqref{estimator gamma}, meaning that the calibration band is constructed using a full set of ordered pairs $\mathcal{J}^{full}$ and only six observations having large and different volumes due to aggregation. The corresponding rejection rates are provided in Table \ref{Tab2}. As expected, the power is now much better as the size of the set of ordered pairs heavily decreases, while we keep using the whole dataset to construct the calibration bands. Although the number of rejections are significantly higher, the same conclusions as for Table \ref{Tab1} hold. That is, the power of the test when the observations are not contaminated is still below the significance level $\alpha = 0.05$ and the test is more effective at identifying large deviations from calibration and contaminations that happen on a global scale or at the middle of the range of interest.
\begin{table}[htb!]
    \centering
    \begin{tabular}{l c c c}
      \toprule
    Contamination $\delta$ & $0$ & $0.5$ & $1$\\
      \midrule
      Global shift &$16/1000$& $995/1000$ & $1000/1000$ \\
      Local shift of level $l = 10$ & - & $360/1000$ & $986/1000$ \\
      Local shift of level $l = 13$ & - &$769/1000$& $1000/1000$ \\
      Local shift of level $l = 15$ & - & $203/1000$ & $917/1000$ \\
      \bottomrule 
    \end{tabular}
    \caption{Power of the performed statistical tests with confidence level $1-\alpha = 0.95$ for binned responses.}
    \label{Tab2}
\end{table}

This example shows that the calibration bands we constructed in this paper can be wide, leading to statistical tests with lower power. We mention again that the reason for this lies in the union bound inequality \eqref{prob bound J} that might not be very sharp for large sets of ordered pairs. An interesting tool to reduce this size while using all the observations is to bin those observations. As a result, the bands get narrower and, thus, allow for more powerful statistical tests. This method assumes that the observations within a given bin have approximately the same mean and, by construction, this is the case in this example.

\section{Conclusion}
\label{sec: Conclusion}

Using the stochastic ordering properties and the convolution formulas of the EDF, we extended the construction of the calibration bands on the mean from the binary case of Dimitriadis et al.~\cite{Dimitriadis} to the whole EDF. Our construction enables us to find closed form expressions for the calibration bands of independent binomial, Poisson, negative binomial, gamma and normal responses and, the bands can be computed using a root-finding algorithm in the other cases of the EDF. Interestingly, we showed that our band is narrower than Yang--Barber's \cite{Yang_Barber} band in the normal case.

As for the calibration bands derived by Dimitriadis et al.~\cite{Dimitriadis} and Yang--Barber \cite{Yang_Barber}, our construction relies on the assumption that the responses are ranked such that their true means are increasing. In a regression modeling context, we showed how this assumption can be extended by introducing a ranking function that provides the ordering of the true mean function for almost every feature in the feature space. In practice, such a ranking function is often unknown and in such cases, it has to be approximated by the regression function itself in order to construct the statistical tests for calibration or auto-calibration. Through numerical examples, we showed how these tests can be applied to detect violations of these properties, and we emphasized that in contrast to other approaches, our tests can be applied for any sample size as the construction of the calibration bands does not rely on asymptotic results. Moreover, we discussed some important factors that influence the shape of the bands 
and proposed methods to construct them for large datasets. One of these methods consists in binning the available observations, and we argued that it leads to suitable bands as it allows for using all the observations while remaining computationally efficient.

The decision of the statistical tests we propose depends on the underlying calibration band, which itself depends on the chosen ranking function.
Going forward, it will be interesting to better understand the role of the ranking function on the resulting band. This is particularly true in cases where the ranking function cannot easily be inferred from the observations, e.g., when the signal-to-noise ratio of the available observations is low. 
Another next step is to study alternative methods for binning the observations and other choices for the set of ordered pairs in order to understand the impact on the resulting band. Finally, the rate of convergence of the calibration bands for an increasing amount of observations is of interest, as well as asymptotic results.

\bigskip

{\bf Acknowledgments.} We kindly thank the referees and the editor for their useful remarks that have helped us to improve this manuscript.

\bigskip

{\small 
\renewcommand{\baselinestretch}{.51}
}

\newpage

\appendix


\newpage

\appendix


\section{Proofs}

\label{proofs}
We prove all statements in this appendix.

\bigskip

{\Beweis {\bf Proof of Proposition \ref{stochastic order EDF}.}
     Let $t$ be in the support of the density $f_{Y_1}$ (which is the same as the support of the density $f_{Y_2}$, see Section \ref{sec definitions}) and set $\theta_i = h(\mu_i)$ for $i \in \{1,2\}$. Since the canonical link $h$ is strictly increasing on $\kappa'(\interior{\Theta})$, we have $\theta_1 \leq \theta_2$. Thus, if we divide the density of $Y_2$ by the density of $Y_1$, we obtain that the function
    \begin{equation*}
        \begin{split}
            t \mapsto \frac{f_{Y_2}(t)}{f_{Y_1}(t)} &= \frac{\exp \left\{\frac{t \theta_2 - \kappa(\theta_2)}{\varphi/v} + a(t; v/\varphi)\right\}}{\exp \left\{\frac{t \theta_1 - \kappa(\theta_1)}{\varphi/v} + a(t; v/\varphi)\right\}} \\
            &= \exp \left\{\frac{t (\theta_2 - \theta_1) - \kappa(\theta_2) + \kappa(\theta_1)}{\varphi/v}\right\}
        \end{split}
    \end{equation*}
    is non-decreasing. This implies $Y_1 \leq_{lr} Y_2$ and using Theorem 1.C.1 in Shaked--Shanthikumar \cite{Stochastic Orders}, we conclude that $Y_1 \leq_{st} Y_2$.
\bigskip
\EndProof}

{\Beweis {\bf Proof of Proposition \ref{bound one EDF}.}
    The lower bound in \eqref{lower bound delta} satisfies
    \begin{equation*}
        \begin{split}
        \p \left( \E[Y] \geq l^{\delta}(Y, v, \varphi, \kappa(\cdot))  \right) &= 1 - \p \left( \E[Y] < l^{\delta}(Y, v, \varphi, \kappa(\cdot))  \right) \\
        &\geq 1 - \p \left( F^{*}(Y; h(\E[Y]), v, \varphi, \kappa(\cdot)) > 1-\delta  \right) \\
        &\geq 1- \p (U > 1- \delta) = 1-\delta,
        \end{split}
    \end{equation*}
    with $U \sim \textnormal{Unif}([0,1])$ and where we used in the first inequality that
    \begin{equation*}
        \E[Y] < l^{\delta}(Y, v, \varphi, \kappa(\cdot)) \implies F^{*}(Y; h(\E[Y]), v, \varphi, \kappa(\cdot)) > 1-\delta.
    \end{equation*}
    Similarly, we have for the upper bound in \eqref{upper bound delta} that
    \begin{equation*}
        \begin{split}
        \p \left( \E[Y] \leq u^{\delta}(Y, v, \varphi, \kappa(\cdot))  \right) &= 1 - \p \left( \E[Y] > u^{\delta}(Y, v, \varphi, \kappa(\cdot))  \right) \\
        &\geq 1 - \p \left( F(Y; h(\E[Y]), v, \varphi, \kappa(\cdot)) < \delta  \right) \\
        &\geq 1- \p (U < \delta) = 1-\delta.
        \end{split}
    \end{equation*}
\bigskip
\EndProof}

{\Beweis {\bf Proof of Lemma \ref{stochastic bound j k}.}
    The proof relies on Theorem 1.A.3 in Shaked--Shanthikumar \cite{Stochastic Orders}, which states that for any set of independent random variables $X_1, \dots, X_n$ and any another set of independent random variables $Y_1, \dots, Y_n$ satisfying $X_i \leq_{st} Y_i$ for all $i \in \{1, \dots, n\}$, we have
    \begin{equation*}
        \psi(X_1, \dots, X_n) \leq_{st} \psi(Y_1, \dots, Y_n),
    \end{equation*}
    for any non-decreasing function $\psi: \R^n \to \R$. By choosing $\psi(x_1, \dots, x_n) = \sum_{i=1}^n v_i x_i / v_{1:n}$ and using Proposition \ref{stochastic order EDF}, the stochastic ordering result holds for
    \begin{equation*}
        Z_{j:k}^{-} = \frac{1}{v_{j:k}} \sum_{i=j}^{k} v_i Y_i^{-}, \quad \textrm{with} \quad Y_i^{-} \stackrel{\mathrm{ind.}}\sim \textnormal{EDF}(\theta_j, v_i, \varphi, \kappa(\cdot)),
    \end{equation*}
    and
    \begin{equation*}
        Z_{j:k}^{+} = \frac{1}{v_{j:k}} \sum_{i=j}^{k} v_i Y_i^{+},\quad \textrm{with} \quad Y_i^{+} \stackrel{\mathrm{ind.}} \sim \textnormal{EDF}(\theta_k, v_i, \varphi, \kappa(\cdot)).
    \end{equation*}
    Note that these random variables satisfy \eqref{random variables jk} due to the convolution formula for the EDF in Corollary 2.15 of Wüthrich--Merz \cite{WM2023}. This completes the proof.

\bigskip
\EndProof}

{\Beweis {\bf Proof of Proposition \ref{last corollary}.}
    Let $1 \leq j \leq k \leq n$, it follows from Proposition \ref{bound one EDF} that
    \begin{equation*}
        \p\left(  \mu_j \leq  u^{\delta}(Z_{j:k}^{-}, v_{j:k}, \varphi, \kappa(\cdot))   \right) \geq 1- \delta \quad \textrm{and} \quad \p\left(  \mu_k \geq  l^{\delta}(Z_{j:k}^{+}, v_{j:k}, \varphi, \kappa(\cdot))   \right) \geq 1- \delta,
    \end{equation*}
    for the random variables $Z_{j:k}^{-}$ and $Z_{j:k}^{+}$ introduced in \eqref{random variables jk}. Since $Z_{j:k}^{-} \leq_{st} Z_{j:k} \leq_{st} Z_{j:k}^{+}$ due to Lemma \ref{stochastic bound j k} and since the functions 
    \begin{equation*}
        y \in \R \mapsto u^{\delta}(y, v_{j:k}, \varphi, \kappa(\cdot)),
    \end{equation*}
    and
    \begin{equation*}
        y \in \R \mapsto l^{\delta}(y, v_{j:k}, \varphi, \kappa(\cdot)),
    \end{equation*}
    are non-decreasing in $y$, we have by Theorem 1.A.3 of Shaked-Shanthikumar \cite{Stochastic Orders} that 
    \begin{equation*}
        u^{\delta}(Z_{j:k}^{-}, v_{j:k}, \varphi, \kappa(\cdot)) \leq_{st} u^{\delta}(Z_{j:k}, v_{j:k}, \varphi, \kappa(\cdot)),
    \end{equation*}
    and 
    \begin{equation*}
        l^{\delta}(Z_{j:k}, v_{j:k}, \varphi, \kappa(\cdot)) \leq_{st} l^{\delta}(Z_{j:k}^{+}, v_{j:k}, \varphi, \kappa(\cdot)).
    \end{equation*}
    The claim then follows.
\bigskip
\EndProof}

{\Beweis {\bf Proof of Theorem \ref{theo calibration bands}.}
Let $\mathcal{J}$ be any set of ordered pairs. By Proposition \ref{last corollary}, we deduce using a union bound argument that
\begin{equation*}
    \p \Bigl(\mu_j \leq  u^{\delta}(Z_{j:k}, v_{j:k}, \varphi, \kappa(\cdot) ) \textrm{ and } \mu_k \geq  l^{\delta}(Z_{j:k}, v_{j:k}, \varphi, \kappa(\cdot)) \textrm{ for all } (j,k) \in \mathcal{J} \Bigl) \geq 1-2|\mathcal{J}|\delta.
\end{equation*}
Due to the ordering assumed in \eqref{ranking no regression}, the above inequality can be rewritten as
\begin{equation*}
    \p\left(  \sup_{(j,k) \in \mathcal{J} \, : \, \theta_i \geq \theta_k} l^{\delta}(Z_{j:k}, v_{j:k}, \varphi, \kappa(\cdot)) \leq \mu_i \leq  \inf_{(j,k) \in \mathcal{J} \, : \, \theta_i \leq \theta_j}u^{\delta}(Z_{j:k}, v_{j:k}, \varphi, \kappa(\cdot)) \, \textnormal{ for all } \, i \in \{1, \dots, n\} \right) \geq 1-2|\mathcal{J}|\delta. 
\end{equation*}
Choosing $\alpha = 2|\mathcal{J}|\delta$ provides the claim.

\bigskip
\EndProof}

{\Beweis {\bf Proof of Proposition \ref{prop discrete distributions}.} The binomial, Poisson and negative binomial distributions belong to the additive form of the EDF for carefully chosen canonical parameters, volumes, dispersion parameters and cumulant functions, we refer to Table 3.3 in J{\o}rgenssen \cite{Jorgenssen2}. The lower and upper bounds defined in \eqref{L_i general}-\eqref{U_i general} can be explicitly expressed for the following three cases using the weighted partial sums $Z_{j:k}$ and aggregated volumes $v_{j:k}$ in \eqref{sum j k}-\eqref{aggregated volume}.

\medskip

\textit{Binomial case.}
    The lower bound in \eqref{L_i bin} can be expressed as
    $$L_{\bY, i}^\alpha = \sup_{(j,k) \in \mathcal{J} \, : \, \mu_i \geq \mu_k} l^{\delta}(Z_{j:k}, v_{j:k}, \varphi),$$
    with
    $$l^{\delta}(Z_{j:k}, v_{j:k},\varphi)= \inf \left\{\mu \in (0,1) \,| \, F^{*}(v_{j:k} Z_{j:k}/\varphi; v_{j:k}/\varphi, \mu) \leq 1-\delta \right\},$$
    and where $F^{*}(\cdot; m, \mu)$ denotes the left-continuous, right-limit distribution of a $\textrm{Bin}(m,\mu)$ random variable. Let $I_\mu(x,y)$ be the regularized incomplete beta function. For $l \in \{1, \dots, m\}$, we have
    \begin{equation*}
        \begin{split}
            F^{*}(l; m, \mu) \leq 1-\delta &\iff 1-I_{\mu} (1+(l-1),m-(l-1)) \leq 1-\delta, \\
            &\iff 1-G_{B}(\mu; l, 1+m-l) \leq 1-\delta, \\
            &\iff \mu \geq q_{B}(\delta; l, 1+m-l),
        \end{split}
    \end{equation*}

where $G_{B}(y; \alpha, \beta)$ and $q_{B}(\delta; \alpha, \beta)$ denote the distribution and the $\delta$-quantile of a $\textrm{Beta}(\alpha, \beta)$ random variable, respectively. This shows the claim for the lower bound. The result for the upper bound in \eqref{U_i bin} follows similarly as
\begin{equation*}
        \begin{split}
            F(l; m, \mu) \geq \delta &\iff 1-I_{\mu} (1+l,m-l) \geq \delta, \\
            &\iff 1-G_{B}(\mu; 1+l, m-l) \geq \delta, \\
            &\iff \mu \leq q_{B}(1-\delta; 1+l, m-l),
        \end{split}
    \end{equation*}
    where $l \in \{0, \dots, m-1\}$. This completes the proof.

\medskip

\textit{Poisson case.}
    The lower bound in \eqref{L_i pois} can be expressed as
    $$L_{\bY, i}^\alpha = \sup_{(j,k) \in \mathcal{J} \, : \, \mu_i \geq \mu_k} l^{\delta}(Z_{j:k}, v_{j:k}, \varphi),$$
    with
    $$l^{\delta}(Z_{j:k}, v_{j:k}, \varphi)= \inf \left\{\mu \in (0, \infty) \,| \, F^{*}(v_{j:k} Z_{j:k}/\varphi; \mu v_{j:k}/\varphi) \leq 1-\delta \right\},$$
    and where $F^{*}(\cdot; \mu v/\varphi)$ denotes the left-continuous, right-limit distribution of a $\textrm{Poi}(\mu v/\varphi)$ random variable. For $l \in \N$, we have
    \begin{equation*}
        \begin{split}
            F^{*}(l; \mu v/\varphi) \leq 1-\delta &\iff \frac{\Gamma(l, \mu v/\varphi)}{(l-1)!} \leq 1-\delta, \\
            &\iff \frac{\Gamma(l, \mu v/\varphi)}{\Gamma(l)} \leq 1-\delta,  \\
            &\iff 1-G_{\Gamma}(\mu v/\varphi; l, 1) \leq 1-\delta, \\
            &\iff \mu v/\varphi \geq\, q_{\Gamma}(\delta; l, 1), \\
            &\iff \mu \geq \, \frac{\varphi q_{\Gamma}(\delta; l, 1)}{v},
        \end{split}
    \end{equation*}

where $G_{\Gamma}(y; \gamma, c)$ and $q_{\Gamma}(\delta; \gamma, c)$ denote the distribution and the $\delta$-quantile of a $\Gamma(\gamma,c)$ random variable, respectively. This shows the claim for the lower bound. Similarly, the result for the upper bound in \eqref{U_i pois} follows from
\begin{equation*}
        \begin{split}
            F(l; \mu v/\varphi) \geq \delta &\iff \frac{\Gamma(1+l, \mu v/\varphi)}{l!} \geq \delta, \\
            &\iff \frac{\Gamma(1+l, \mu v/\varphi)}{\Gamma(1+l)} \geq \delta,  \\
            &\iff 1-G_{\Gamma}(\mu v/\varphi; 1+l, 1) \geq \delta, \\
            &\iff \mu v/\varphi \leq \, q_{\Gamma}(1-\delta; 1+l, 1), \\
            &\iff \mu \leq \frac{\varphi q_{\Gamma}(1-\delta; 1+l, 1)}{v},
        \end{split}
    \end{equation*}
    where $l \in \N_0$. 

\medskip

{\it Negative binomial case}.
    The lower bound in \eqref{L_i neg bin} can be expressed as
    $$L_{\bY, i}^\alpha = \sup_{(j,k) \in \mathcal{J} \, : \, \mu_i \geq \mu_k} l^{\delta}(Z_{j:k}, v_{j:k}, \varphi),$$
    with
    $$l^{\delta}(Z_{j:k}, v_{j:k}, \varphi)= \inf \left\{\mu \in (0, \infty) \,| \, F^{*}(v_{j:k} Z_{j:k}/\varphi; \mu, v_{j:k}/\varphi) \leq 1-\delta \right\},$$
    and where $F^{*}(\cdot; \mu, \gamma)$ denotes the left-continuous, right-limit distribution of a $\textrm{NegBin}(\mu, \gamma)$ random variable with mean parameter $\mu > 0$ and shape parameter $\gamma > 0$. Define $p= \mu/(1 + \mu)$ and let $I_p(x,y)$ be the regularized incomplete beta function. For $l \in \mathbb{N}$, we have
    \begin{equation*}
        \begin{split}
            F^{*}(l; \mu, v/\varphi) \leq 1-\delta &\iff 1-I_{p} (1+(l-1), v/\varphi) \leq 1 - \delta,\\
            &\iff 1-G_B(p; l, v/\varphi) \leq 1-\delta, \\
            &\iff p \geq q_{B}(\delta; l, v/\varphi), \\
            &\iff \frac{\mu}{1+\mu} \geq q_{B}(\delta; l, v/\varphi), \\
            &\iff \mu \geq \frac{q_{B}(\delta; l, v/\varphi)}{1-q_{B}(\delta; l, v/\varphi)},
        \end{split}
    \end{equation*}

where $G_{B}(y; \alpha, \beta)$ and $q_{B}(\delta; \alpha, \beta)$ denote the distribution and the $\delta$-quantile of a $\textrm{Beta}(\alpha, \beta)$ random variable, respectively. This shows the claim for the lower bound. Similarly, the result for the upper bound in \eqref{U_i neg bin} follows from
\begin{equation*}
        \begin{split}
            F(l; \mu, v/\varphi) \geq \delta &\iff 1-I_{p} (1+l,v/\varphi) \geq \delta, \\
            &\iff 1-G_B(p; 1+l, v/\varphi) \geq \delta, \\
            &\iff p \leq q_{B}(1-\delta; 1+l, v/\varphi), \\
            &\iff \frac{\mu}{1+\mu} \leq q_{B}(1-\delta; 1+l, v/\varphi), \\
            &\iff \mu \leq \frac{q_{B}(1-\delta; 1+l, v/\varphi)}{1-q_{B}(1-\delta; 1+l, v/\varphi)},
        \end{split}
    \end{equation*}
    where $l \in \N_0$. 
\bigskip
\EndProof}

{\Beweis {\bf Proof of Proposition \ref{prop continuous distributions}.}
The gamma and normal distributions belong to the EDF for carefully chosen canonical parameters, volumes, dispersion parameters and cumulant functions, we refer to Table 3.1 in J{\o}rgenssen \cite{Jorgenssen2}. The lower and upper bounds defined in \eqref{L_i general}-\eqref{U_i general} can be explicitly expressed for the following two cases using the weighted partial sums $Z_{j:k}$ and aggregated volumes $v_{j:k}$ in \eqref{sum j k}-\eqref{aggregated volume}.
\medskip

{\it Gamma case}. The lower bound in \eqref{L_i gamma} can be expressed as
    $$L_{\bY, i}^\alpha = \sup_{(j,k) \in \mathcal{J} \, : \, \mu_i \geq \mu_k} l^{\delta}(Z_{j:k}, v_{j:k}, \varphi),$$
    with
    $$l^{\delta}(Z_{j:k}, v_{j:k}, \varphi)= \inf \left\{\mu \in (0, \infty) \,| \, F^{*}(Z_{j:k}; v_{j:k}/\varphi, v_{j:k}/(\varphi\mu)) \leq 1-\delta \right\},$$
    and where $F^{*}(\cdot; v/\varphi, v/(\varphi\mu))$ denotes the left-continuous, right-limit distribution of a $\Gamma(v/\varphi, v/(\varphi\mu))$ random variable. Let $\mathcal{G}(x,y)$ be the lower incomplete gamma function. For $l > 0$, we have
    \begin{equation*}
        \begin{split}
            F^{*}(l; v/\varphi, v/(\varphi\mu)) \leq 1-\delta &\iff \frac{\mathcal{G}(v/\varphi,vl/(\varphi\mu))}{\Gamma(v)} \leq 1 - \delta,\\
            &\iff G_\Gamma(v/(\varphi\mu); v/\varphi, l) \leq 1-\delta, \\
            &\iff \frac{v}{\varphi\mu} \leq q_{\Gamma}(1-\delta; v/\varphi, l), \\
            &\iff v/\varphi \leq \mu \cdot q_{\Gamma}(1-\delta; v/\varphi, l),  \\
            &\iff \mu \geq \frac{v/\varphi}{q_{\Gamma}(1-\delta; v/\varphi, l)},
        \end{split}
    \end{equation*}

where $G_{\Gamma}(y; \alpha, \beta)$ and $q_{\Gamma}(\delta; \alpha, \beta)$ denote the distribution and the $\delta$-quantile of a $\Gamma(\gamma, c)$ random variable, respectively. This shows the claim for the lower bound and as the left-continuous, right-limit distribution of a gamma random variable coincides with the distribution of this random variable, the result for the upper bound in \eqref{U_i gamma} follows similarly.
\medskip

{\it Normal case.}
    The lower bound in \eqref{L_i normal} can be expressed as
    $$L_{\bY, i}^\alpha = \sup_{(j,k) \in \mathcal{J} \, : \, \mu_i \geq \mu_k} l^{\delta}(Z_{j:k}, v_{j:k}, \varphi),$$
    with
    $$l^{\delta}(Z_{j:k}, v_{j:k}, \varphi)= \inf \left\{\mu \in \R \,| \, F(Z_{j:k}; \mu, v_{j:k}/\varphi) \leq 1-\delta \right\},$$
    and where $F(\cdot; \mu, v_{j:k}/\varphi)$ denotes the distribution of a $\mathcal{N}(\mu, \varphi/v_{j:k})$ random variable. 
    This pointwise infimum is always attained as the map
    \begin{equation*}
        \mu \mapsto \p(Y_{\mu, v_{j:k}/\varphi} \leq z), \textrm{ where } Y_{\mu, v_{j:k}/\varphi} \sim \mathcal{N}(\mu, \varphi/v_{j:k}),
    \end{equation*}
    is continuous for a fixed $z \in \R$.
    The lower band $l = l^{\delta}(Z_{j:k}, v_{j:k}, \varphi)$ thus satisfies
    \begin{equation*}
        \begin{split}
            F(Z_{j:k}; l, v_{j:k}/\varphi) = 1-\delta &\iff \Phi \left( \sqrt{v_{j:k}/\varphi} \left(Z_{j:k} - l \right)  \right) = 1- \delta ,\\
            &\iff \sqrt{v_{j:k}/\varphi} \left(Z_{j:k} - l \right) = \Phi^{-1}\left(1-\delta\right),\\
            &\iff l = Z_{j:k} - \frac{\Phi^{-1}(1-\delta)}{\sqrt{v_{j:k}/\varphi}}.
        \end{split}
    \end{equation*}
    A similar derivation provides the result for the upper bound in \eqref{U_i normal}.
\bigskip
\EndProof}

{\Beweis {\bf Proof of Theorem \ref{our bands better}.}

    From Proposition B1 of Dimitriadis et al.~\cite{Dimitriadis} and by choosing $\tau = \sqrt{2 \sigma^2 \log(1/\delta)}$, we know that 
    \begin{equation*}
        U^{\alpha,YB}_{\bY ,i} = Z_{j:k}^{Iso} + \frac{\tau}{\sqrt{k-j+1}},
    \end{equation*}
    for some pair $(j,k) \in \mathcal{J}^{full}$ with $j = i$ and such that either $\widehat{\mu}^{Iso}(\bY, \bone)_k < \widehat{\mu}^{Iso}(\bY, \bone)_{k+1}$ or $k = n$. Moreover, we have
    \begin{equation*}
        Z_{j:k} \leq Z_{j:k}^{Iso}, \textrm{ whenever }\widehat{\mu}^{Iso}(\bY, \bone)_k <  \widehat{\mu}^{Iso}(\bY, \bone)_{k+1} \textrm{ or } k = n,
    \end{equation*}
    due to a property of isotonic regression (we refer, e.g., to Characterization II provided by Henzi et al.~\cite{Henzi}). By defining the new set
    \begin{equation*}
        \tilde{\mathcal{J}}_i = \left\{  \left. (i,k) \in \mathcal{J}^{full} \, \right| \, \widehat{\mu}^{Iso}(\bY, \bone)_k < \widehat{\mu}^{Iso}(\bY, \bone)_{k+1}  \textrm{ or }  k = n\right\},
    \end{equation*}
    for $1 \leq i \leq n$, we thus obtain
    \begin{equation*}
        \begin{split}
            U^{\alpha,YB}_{\bY ,i} &= \min_{(j,k) \in \tilde{\mathcal{J}}_i} Z_{j:k}^{Iso} + \frac{\sqrt{2 \sigma^2 \log(1/\delta)}}{\sqrt{k-j+1}} \\ 
            &\geq \min_{(j,k) \in \tilde{\mathcal{J}}_i} Z_{j:k} + \frac{\sqrt{2 \sigma^2 \log(1/\delta)}}{\sqrt{k-j+1}} \\
            &\geq \min_{(j,k) \in \tilde{\mathcal{J}}_i} Z_{j:k} - \frac{\sigma \Phi^{-1}(\delta)}{\sqrt{k-j+1}} \\
            &\geq \min_{(j,k) \in \mathcal{J}^{full} : \, \mu_i \leq \mu_j} Z_{j:k} - \frac{\sigma \Phi^{-1}(\delta)}{\sqrt{k-j+1}} = U^{\alpha}_{\bY ,i},
        \end{split}
    \end{equation*}
    where we used that $\sqrt{2 \log(1/\delta)} \geq - \Phi^{-1}(\delta)$ for all $\delta \in (0,1)$. A similar computation provides the result for the lower band.

\bigskip
\EndProof}

{\Beweis {\bf Proof of Theorem \ref{cor calibration bands}.}
Define the set
        \begin{equation}
        \label{set A}
        A = \left\{ L^\alpha_{\pi, (Y_i, \bX_i)_{i=1}^n}(\bX_l) \leq \mu_\pi^{*}(\bX_l) \leq  U^\alpha_{\pi, (Y_i, \bX_i)_{i=1}^n}(\bX_l)  \, \textnormal{ for all } \, l \in \{1, \dots, n\} \right\},
    \end{equation}
    that lies in $\mathcal{F}$ as all the random variables involved in its definition are measurable. Instead of proving the existence of the uniform calibration band in \eqref{prob bound Q_x} holding simultaneously for all $\bx \in \mathcal{X}$, we first prove that for a.e.~realization $(\bx_i)_{i=1}^n$ of the features $(\bX_i)_{i=1}^n$, we have
    \begin{equation}
    \label{prob bound Q_x i}
    \Q_{(\bx_i)_{i=1}^n} \Bigl(  L^\alpha_{\pi, (Y_i, \bX_i)_{i=1}^n}(\bX_l) \leq \mu_\pi^{*}(\bX_l) \leq  U^\alpha_{\pi, (Y_i, \bX_i)_{i=1}^n}(\bX_l) \, \textnormal{ for all } \, l \in \{1, \dots, n\} \Bigl) \geq 1- \alpha.
    \end{equation}
    Note that since the random variables $(Y_i, \bX_i)_{i=1}^n$ are independent, they satisfy  $$Y_{i} \, | \, \bX_1 = \bx_1, \dots, \bX_n = \bx_n \sim \textnormal{EDF}(\theta(\bx_{i}), v_{i}, \varphi, \kappa(\cdot)), \quad \textrm{for } 1 \leq i \leq n.$$ By using the permutation function $\tau_{\bx_1, \dots, \bx_n}$ introduced in \eqref{permutation p}, the indices of these random variables can be permuted such that $$
    \theta(\bx_{\tau_{
\bx_1, \dots, \bx_n}(1)}) \leq \dots \leq \theta(\bx_{\tau_{\bx_1, \dots, \bx_n}(n)}),$$
whenever the features $\bx_1, \dots, \bx_n$ satisfy $\mu_\pi^{*}(\bx_i) = \mu^{*}(\bx_i)$ for all $i \in \{1, \dots, n\}$. The latter happens for a.e.~realization of the features $(\bX_i)_{i=1}^n$ and in this case, the inequality in \eqref{prob bound Q_x i} follows from Theorem \ref{theo calibration bands}. 
    In order to show the inequality in \eqref{prob bound Q_x}, it suffices then to notice that under Assumption \ref{ranking pi}, the set $A$ in \eqref{set A}
    is equal to the set $B$ given by
    \begin{equation*}
        B = \left\{ L^\alpha_{\pi, (Y_i, \bX_i)_{i=1}^n}(\bx) \leq \mu_\pi^{*}(\bx) \leq  U^\alpha_{\pi, (Y_i, \bX_i)_{i=1}^n}(\bx)  \, \textnormal{ for all } \, \bx \in \mathcal{X} \right\}.
    \end{equation*}
    This concludes the proof.

\bigskip
\EndProof}

\end{document}